\documentclass[11pt,letterpaper]{article}

\usepackage[margin=1in]{geometry}
\usepackage[T1]{fontenc}
\usepackage{newtxtext}
\usepackage{microtype}
\usepackage{setspace}
\usepackage{amsmath,amsthm,mathtools}
\usepackage{newtxmath}
\newtheorem{theorem}{Theorem}
\newtheorem{lemma}{Lemma}
\newtheorem{proposition}{Proposition}
\newtheorem{corollary}{Corollary}

\usepackage{algorithm}
\usepackage{algpseudocode}

\makeatletter
\newenvironment{breakablealgorithm}
  {\begin{center}
   \vspace{0.75\baselineskip}
   \refstepcounter{algorithm}%
   \hrule height .8pt depth 0pt \kern 2pt%
   \renewcommand{\caption}[2][\relax]{%
     {\raggedright\textbf{\ALG@name~\thealgorithm} ##2\par}%
     \ifx\relax##1\relax
       \addcontentsline{loa}{algorithm}{\protect\numberline{\thealgorithm}##2}%
     \else
       \addcontentsline{loa}{algorithm}{\protect\numberline{\thealgorithm}##1}%
     \fi
     \kern 2pt\hrule\kern 2pt}%
  }
  {\kern 2pt\hrule\relax\end{center}\vspace{0.75\baselineskip}}
\makeatother

\usepackage{graphicx,caption,subcaption,adjustbox,float,xcolor}
\usepackage{tabularx,multirow,booktabs,longtable,array,tikz}
\usepackage[most]{tcolorbox}
\tcbuselibrary{breakable}
\usepackage{enumitem}
\newlist{steps}{enumerate}{1}
\setlist[steps,1]{wide=0pt,leftmargin=\parindent,label=Step \arabic*:,font=\bfseries}

\newcommand{\OmegaSet}{\Omega}
\usepackage[round,authoryear]{natbib}
\bibpunct[, ]{(}{)}{,}{a}{}{,}

\usepackage[colorlinks=true,linkcolor=blue,citecolor=blue,urlcolor=blue]{hyperref}

\newcommand{\Halmos}{\hfill\ensuremath{\square}}

\makeatletter
\newcommand{\arxivincludegraphics}[2][]{%
  \IfFileExists{#2}{\includegraphics[#1]{#2}}{%
    \fbox{\parbox[c][1.1in][c]{0.85\linewidth}{\centering\footnotesize Missing figure: \texttt{\detokenize{#2}}}}%
  }%
}
\makeatother

\begin{document}
\title{Recursive Pricing and Incumbent Subtree Column Generation for Solving Optimal Classification Trees}
\author{%
Tengmu Hu\\
\small Supply Chain, Operations, and Technology\\
\small Olin Business School, Washington University in St. Louis\\
\small \texttt{tengmu@wustl.edu}
\and
Chen Chen\\
\small Integrated Systems Engineering\\
\small The Ohio State University\\
\small \texttt{chen.8018@osu.edu}
\and
Theodore T. Allen\\
\small Industrial and Systems Engineering\\
\small University of Massachusetts Amherst\\
\small \texttt{tedallen@umass.edu}
}
\date{}
\maketitle

\begin{abstract}
Solving optimal classification tree (OCT) problems is computationally challenging because branching, routing, and leaf-assignment decisions must be optimized jointly. We develop an incumbent subtree column generation algorithm (ISCG-OCT) that exploits the recursive structure of classification trees. We first introduce a mixed set-packing and set-partitioning formulation (MSSF). We then construct a restricted master problem whose columns represent feasible subtree patterns. Unlike conventional dual-based column generation, our incumbent-centered pricing evaluates objective improvements relative to the current solution. We establish a recursive decomposition in which fixing a branching decision separates the remaining pricing problem into independent subtree subproblems, enabling exact dynamic programming. We also develop threshold-based variants that trade optimality guarantees for computational efficiency. Computational experiments on 27 benchmark datasets demonstrate competitive training and testing accuracy and substantial reductions in solution time on many instances. When initialized with a CART solution, ISCG-OCT can certify global optimality faster than warm-started mixed-integer programming solvers on selected instances. These results demonstrate how recursive pricing and structured subtree optimization can improve the computational efficiency of exact classification tree methods.
\end{abstract}



\section{Introduction}

Classification trees have gained popularity because of their simple and interpretable structure. Compared to black-box models such as neural networks, classification trees offer the advantage of transparency, making it easier to identify and explain the features used by the model. In scenarios where interpretability and feature-level decision transparency are essential, classification trees are a strong choice. Moreover, they are based purely on the dataset and do not require standard parametric assumptions or assumed data distributions. This flexibility allows them to be applied across a wide range of practical situations. In addition, classification trees can handle both numerical and categorical data, making them suitable for various applications.

A classification tree predicts categorical outcomes by recursively partitioning observations based on their feature values. In a binary classification tree, each branch node selects a feature and a corresponding split value to determine how observations are routed. Observations with feature values less than or equal to the selected split value are directed to the left subtree, while those with greater values are directed to the right subtree. This recursive partitioning continues until a stopping criterion is met, such as reaching a predefined maximum tree depth. Each leaf node is assigned a class label that determines the prediction for all observations reaching that leaf. An observation is correctly classified if its predicted label matches its true class. The objective of optimal classification tree (OCT) is to jointly determine the splitting decisions and leaf-class assignments that maximize the number of correctly classified observations, subject to constraints on the tree structure.

Despite their intuitive structure and interpretability, solving OCT problems is computationally challenging because split-feature selection, data-point routing, and leaf-class assignment decisions are interdependent and must be optimized jointly. Since splitting decisions at upper-level nodes affect the observations available to their descendant nodes, these decisions are closely interconnected, resulting in a challenging optimization problem. Consequently, widely used classification tree algorithms typically rely on greedy heuristics that make sequential, locally optimal splitting decisions without guaranteeing global optimality. Although several exact optimization formulations have been proposed to overcome this limitation, their computational requirements can be substantial, particularly as the number of observations, features, or tree depth increases. These challenges motivate the development of more efficient exact algorithms.

In this paper, we develop a recursive-pricing framework for solving OCT problems. We first introduce a mixed set-packing and set-partitioning formulation (MSSF) that integrates branching and assignment decisions while preserving the underlying tree structure. In MSSF, routing is enforced by set-packing compatibility constraints between branching and leaf-assignment decisions rather than by big-$M$ constraints, so each feasible tree corresponds to a selection of branching and assignment decisions that is consistent along every root-to-leaf path. We propose a data-dependent finite reduction of candidate split values based on recursive class-boundary information across feature orderings, thereby reducing the branching search space while retaining class-boundary information identified across recursively generated observation subsets. The reduction preserves optimality for depth-one trees and at every node whose children are leaves, and we show that it can exclude optimal trees at higher levels, so that the optimality certificates of ISCG-OCT are stated relative to the reduced candidate sets.

To overcome the computational limitations of existing exact OCT formulations, we develop a column generation framework that exploits the recursive structure of classification trees. Specifically, we aggregate the interdependent branching and assignment decisions into closed feasible subtree columns and construct a restricted master problem (RMP) over these columns. This formulation gives rise to a structured subtree-pricing problem that jointly determines the branching, deactivation, routing, and leaf-assignment decisions within each candidate subtree. Unlike conventional column generation, which relies on dual-based reduced costs, our pricing procedure is centered on the current integer-feasible tree and evaluates candidate columns based on their objective improvement relative to the incumbent solution. As a result, the algorithm maintains an integer-feasible tree whose objective never deteriorates, so a usable tree is available whenever the procedure is stopped. Its progress does not depend on a linear programming relaxation. And, because the incumbent defines a zero objective change, every partial subtree configuration can be compared with a known cutoff before any complete column is generated.

We show that the subtree-pricing problem arising from the proposed column generation framework can be solved exactly using dynamic programming (DP) by exploiting the recursive structure of classification trees. Specifically, once a branching or deactivation decision is fixed at an internal node, the observations reaching its child nodes are determined, allowing the remaining pricing problem to decompose into independent left- and right-subtree subproblems. This decomposition enables us to derive an exact Bellman recursion that jointly optimizes branching, routing, and leaf-class assignment decisions within each feasible subtree column. To improve computational efficiency, we further incorporate a lower-bound-guided search procedure that eliminates partial subtree configurations that cannot yield sufficiently improving columns. The assignment component of this bound is computed exactly by a hierarchical procedure that respects the common-class requirement at each leaf, so the bound remains valid while the search is screened. These components form the basis of our proposed Incumbent Subtree Column Generation for Optimal Classification Trees (ISCG-OCT) algorithm. Optimality over the candidate split values is certified when exact pricing over the complete tree establishes that no improving feasible tree exists.

The way DP is used distinguishes ISCG-OCT from existing DP methods for OCT. Methods such as DL8.5 \citep{aglin2020learning} and MurTree \citep{demirovic2022murtree} apply DP to the complete OCT problem. They search over the subproblems defined by the possible paths of a tree, typically on binarized features, and cache the optimal subtree of every subproblem they solve. ISCG-OCT applies the recursion only to the subtree being priced, with the observations entering that subtree and all decisions outside it fixed by the incumbent. Because the incumbent supplies an objective-change cutoff from the start, a partial configuration whose lower bound cannot improve on the incumbent is discarded before its descendant states are created. Only values that have been established exactly are kept for reuse, and they are released when the incumbent changes. The memory used by recursive pricing is therefore governed by the states of the current pricing problem that survive this screening, rather than by all subproblems of the complete tree. 

We evaluate the proposed framework on 27 benchmark datasets from the UCI Machine Learning Repository and compare its computational performance with existing exact OCT methods, particularly mixed-integer programming (MIP) solvers equipped with warm starts. Our computational study focuses on the ability of ISCG-OCT to certify optimal solutions over its candidate split values more efficiently than conventional MIP approaches, even when these approaches are provided with initial feasible solutions. The results demonstrate that ISCG-OCT can achieve faster optimality certification on certain instances while maintaining competitive training and testing accuracy. We further investigate threshold-based variants that trade exact optimality certification for reduced computational effort, providing a flexible mechanism for balancing solution quality and computational efficiency.

\subsection{Literature Review}

The broad applicability of classification trees has led to their use across a wide range of domains and case studies. We focus on methodological developments most closely related to optimal classification trees. Recent surveys provide broader perspectives on decision tree learning and interpretable machine learning. \citet{costa2023recent} review developments in decision tree research over the past decade, with particular attention to model fitting, generalization, and interpretability. \citet{ntoutsi2020bias} discuss ethical concerns associated with biased artificial intelligence decisions, particularly when sensitive demographic attributes are involved. \citet{rudin2022interpretable} summarize key principles of interpretable machine learning and identify important open technical challenges. Consistent with the focus of this paper, we restrict attention to single-tree methods and exclude ensemble approaches \citep{mivsic2020optimization,mistry2021mixed}, which combine multiple trees. For example, XGBoost \citep{chen2016xgboost} sequentially constructs multiple trees to improve predictive performance.

Constructing an optimal classification tree is computationally challenging because the problem is NP-hard \citep{laurent1976constructing}. Consequently, a substantial portion of the classification-tree literature has focused on heuristic methods that sacrifice global optimality in exchange for computational efficiency. One of the most widely used approaches is Classification and Regression Trees (CART), introduced by \citet{breiman1984classification}, which employs a greedy procedure that sequentially selects splitting features and thresholds at individual nodes. Other classical greedy methods, including ID3 \citep{quinlan1986induction} and C4.5 \citep{quinlan1993c4}, recursively construct trees using criteria based on information gain or entropy. More recent heuristic approaches incorporate metaheuristics and column-generation-based procedures \citep{santos2022decision,col1}. \citet{bennett1996optimal} combine linear programming with tabu search. Although these methods can often construct high-quality trees efficiently, they do not generally provide guarantees of global optimality.

In parallel with heuristic approaches, considerable attention has been devoted to exact optimization methods for obtaining globally optimal classification trees. Among mathematical programming approaches, \citet{bertsimas2017optimal} propose a mixed-integer programming formulation for optimal classification trees. Subsequent research has developed alternative formulations and solution strategies intended to improve computational performance and robustness. \citet{aghaei2025strong} formulate OCT through a maximum-flow representation and develop a Benders decomposition approach. \citet{verwer2019learning} propose a binary linear programming formulation, while \citet{liu2024optimal} introduce additional constraints to improve the stability of OCT optimization. Constraint programming and SAT-based methods provide further exact approaches \citep{verhaeghe2020learning,narodytska2018learning}. These studies demonstrate the importance of formulation structure and specialized solution procedures in overcoming the computational limitations of direct mixed-integer optimization. Nevertheless, existing methods remain computationally demanding when seeking globally optimal solutions to the OCT problem.

Column generation provides another decomposition-based approach for solving large-scale OCT formulations \citep{barnhart1998branch,lubbecke2005selected}. Previous work has constructed restricted master problems based on root-to-leaf paths. In \citet{col1}, each column represents a path followed by observations from the root to a leaf, requiring the split decisions along that path to satisfy the corresponding routing conditions. An advantage of this representation is the relative simplicity of the restricted master problem because it avoids explicitly modeling some of the classification relationships associated with individual feature values. However, the resulting pricing problem requires verification of the feasibility of generated paths and is computationally challenging. Consequently, heuristic procedures are used to solve the subproblems. Building on this approach, \citet{col2} develop improved heuristic procedures for the associated pricing problem. 

Column generation has also been used to construct interpretable Boolean decision rules \citep{dash2018boolean}. In contrast, the framework proposed in this paper represents branching and assignment decisions jointly within a single family of closed feasible subtree columns and develops exact pricing procedures that exploit the recursive structure of the underlying tree. Because its pricing is based on objective changes relative to an integer-feasible incumbent rather than LP dual reduced costs, the proposed approach is also related to exact large-neighborhood search \citep{ahuja2002very,pisinger2010large}.

DP has emerged as an important exact solution approach for OCT because it can exploit the recursive structure of a decision tree. Rather than relying entirely on a general-purpose mixed-integer solver, DP decomposes the tree-construction problem into smaller subproblems associated with subtrees. Early work by \citet{nijssen2007mining} and subsequent caching-based search \citep{aglin2020learning} illustrate this line of research. \citet{demirovic2022murtree} show that this structure can be exploited to efficiently search the space of possible trees while preserving global optimality under suitable separability conditions. \citet{van2023necessary} further characterize conditions under which such decomposition is valid for more general objectives and constraints. DP has also been combined with bounding techniques to incorporate global requirements \citep{van2022fair}. More recently, DP and branch-and-bound strategies have been integrated to address continuous feature spaces directly, reducing dependence on coarse discretization while improving computational performance \citep{brita2025optimal}. Collectively, these developments demonstrate the value of exploiting recursive tree structure in exact OCT algorithms.

Our framework builds on this insight in a different way. Rather than applying DP directly to the complete OCT problem, we use DP as the exact solution mechanism for the subtree pricing problems induced by the column-generation reformulation. This decomposition restricts each pricing problem to a candidate replacement subtree with a fixed entering observation set, while decisions outside that subtree remain fixed. Consequently, the DP explores only the branching and assignment configurations relevant to the current pricing problem. In addition, the incumbent-centered pricing representation provides immediately available objective-change bounds that can be used to screen partial subtree configurations before complete patterns are generated. Related work on numerical split selection \citep{fayyad1992handling,fayyad1993multi,elomaa1999general} provides context for the recursive class-boundary candidate reduction developed in Section~\ref{sec:mssf}. The distinction is that our reduction considers boundaries arising within recursively generated observation subsets.

\subsubsection{OCT Variants, Special Settings, and Extensions}

Beyond general OCT formulations, several studies develop specialized algorithms for restricted data settings. \citet{lin2020generalized}, \citet{hu2019optimal}, and \citet{xin2022exploring} exploit analytical bounds and specialized data structures to reduce the search space. Other developments include SPO Trees for decision-focused learning \citep{elmachtoub2020decision}, optimization approaches incorporating cost sensitivity, fairness, and interpretability \citep{carrizosa2021mathematical}, sparse oblique trees \citep{carreira2018alternating}, black-box-guided tree construction \citep{mctavish2022fast}, and mixed-integer formulations for categorical and numerical features \citep{gunluk2021optimal}.

Extensions also include randomized decision trees \citep{blanquero2021optimal} and Quant-BnB for continuous features \citep{mazumder2022quant}. The proposed ISCG-OCT framework accommodates numerical, ordinal, and suitably encoded binary features through ordered threshold decisions. Nominal categorical features may be incorporated through appropriate encoding or extensions to the candidate split representation and pricing structure. Such modifications also provide opportunities to exploit preprocessing and structural restrictions to reduce the search space.

\subsection{Paper Organization}

This paper is organized as follows. Section~\ref{sec:mssf} presents the MSSF formulation and proposes a finite reduction of the candidate split values. Section~\ref{sec:column_generation} develops the incumbent-centered column-generation reformulation, including the restricted master problem and the closed feasible subtree-pricing problem that jointly generates branching and assignment decisions. Section~\ref{sec:dynamic_pricing} establishes the recursive pricing structure and develops the DP procedures for solving the pricing subproblems. Section~\ref{subsec:iscg} presents the complete ISCG-OCT algorithm. Section~\ref{cresults} evaluates the proposed framework computationally on benchmark datasets and compares its performance with existing exact OCT methods. Finally, Section~\ref{sec:conclusions} concludes the paper. Proofs, algorithms, and additional computational results are provided in the appendix.

\section{Mixed Set-Packing and Set-Partitioning Formulation}
\label{sec:mssf}

The purpose of the MSSF formulation is not only to provide an exact representation of the OCT problem, but also to expose the structural decomposition required by the column-generation framework developed in Section~\ref{sec:column_generation}. Existing compact OCT formulations typically couple split-feature selection, data-point routing, and leaf-class assignment within a single system of variables and constraints. In contrast, MSSF provides an explicit representation of branching and leaf-assignment decisions together with the compatibility conditions that link them. This structure enables these decisions to be organized into closed feasible subtree patterns, where each pattern jointly specifies the branching and assignment decisions within the corresponding subtree. The resulting patterns are then represented as columns in the RMP, thereby preserving routing and assignment compatibility and avoiding incompatible combinations that could arise if branching and assignment patterns were generated independently.

The formulation combines set-partitioning and set-packing structures. The set-partitioning constraints ensure that each branch node selects exactly one branching decision and that each observation is assigned exactly once. The set-packing constraints exclude combinations of branching and assignment decisions that are incompatible with the routing structure induced by the tree. This separation provides the structural foundation for the pricing and DP procedures developed in the subsequent sections.

In addition to the general MSSF formulation, we establish two structural results. First, Section~\ref{subsec:class_boundary_reduction} introduces a data-dependent reduction of the candidate split values, which reduces the number of branching decisions that need to be considered while preserving the relevant class-boundary information across features. Second, for the special case of a depth-one tree with a single feature, we show in Appendix~\ref{convexhull} that the complete combined $C_p$ set-packing inequalities yield a convex-hull description of the depth-one, single-feature MSSF relaxation under the stated threshold assumptions.

\subsection{Notation}

We consider a classification dataset indexed by $I=\{1,\ldots,N_I\}$. Each observation $i\in I$ is represented by a feature vector $x_i=(x_{ij})_{j\in J}$, where $J=\{1,\ldots,N_J\}$ is the set of features and $x_{ij}$ denotes the value of feature $j$ for observation $i$. Each observation $i$ has a class label $y_i\in K,K=\{1,\ldots,N_K\}$. We use $N_I=|I|$, $N_J=|J|$, and $N_K=|K|$ to denote the corresponding set cardinalities.

A classification tree is represented by a rooted binary tree of fixed maximum depth $N_D$. Let $T=T_B\cup T_L$ denote the set of nodes, where $T_B$ is the set of branch nodes and $T_L$ is the set of leaf nodes. The root node is indexed by $1$ and has depth zero. Each node $v\in T$ has depth $d(v)\in\{0,1,\ldots,N_D\}$, and every leaf node $t\in T_L$ satisfies $d(t)=N_D$. Each branch node $v\in T_B$ has a left child and a right child at depth $d(v)+1$, and every node other than the root has a unique parent. The underlying tree structure is fixed a priori and is therefore part of the problem instance.

\begin{figure}[h]
    \centering
    \includegraphics[width=0.7\columnwidth]{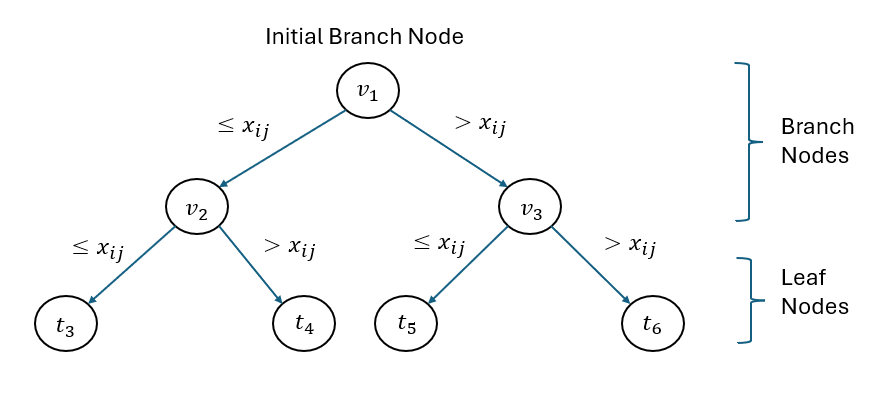}
    \caption{Example of the structure of a classification tree. Node $v_2$ is the parent of leaf nodes $t_3$ and $t_4$, where $t_3$ is the left child of $v_2$ and $t_4$ is the right child of $v_2$.}
    \label{FIG:tree_structure}
\end{figure}

For each leaf node $t\in T_L$, let $P_L(t)\subseteq T_B$ and $P_R(t)\subseteq T_B$ denote the branch nodes on the unique root-to-leaf path to $t$ at which the path proceeds through the left and right child, respectively. Thus, an observation can reach leaf $t$ only if it satisfies the left-routing conditions associated with every node in $P_L(t)$ and the right-routing conditions associated with every node in $P_R(t)$. For the example in Figure~\ref{FIG:tree_structure}, $P_L(t_4)=\{v_1\}$ and $P_R(t_4)=\{v_2\}$.

At each active branch node $v\in T_B$, a feature $j\in J$ and an observed feature value $x_{ij}$ are selected to define the split. Let $b_{ijv}\in\{0,1\}$ equal $1$ if feature $j$ with split value $x_{ij}$ is selected at branch node $v$, and $0$ otherwise. Following \citet{bertsimas2017optimal}, a branch node may also be deactivated. A deactivated node does not perform a split, and every observation reaching that node is routed directly to its right child. We represent this decision using $b_{0v}\in\{0,1\}, v\in T_B$. When $b_{0v}=1$, branch node $v$ is deactivated. When $b_{0v}=0$, the node remains active and a split is selected. Equivalently, deactivation has the routing effect of selecting a threshold below the minimum observed value of any feature, but it is represented by the single alternative $b_{0v}$ rather than by separate feature-specific variables.

The assignment decisions are represented by $z_{ikt}\in\{0,1\},i\in I,\ k\in K,\ t\in T_L$, where $z_{ikt}=1$ if observation $i$ is assigned to leaf $t$ with class label $k$, and $0$ otherwise. The MSSF compatibility constraints ensure that assignments are consistent with the selected branching decisions and that observations assigned to the same nonempty leaf satisfy the required leaf-label consistency.

\subsubsection{Objective Coefficients}

Let $c^z_{ik}$ denote the objective coefficient associated with assigning observation $i\in I$ to class $k\in K$. We define
\[
c^z_{ik}
=
\begin{cases}
-1, & \text{if } k=y_i,\\
 1, & \text{if } k\neq y_i.
\end{cases}
\]
Thus, a correct classification contributes negatively to the objective, whereas an incorrect classification contributes positively. Minimizing the assignment component of the objective therefore favors correct classifications.

Let $c^b_{ij}$ denote the objective coefficient associated with selecting feature value $x_{ij}$ as an active split, and let $c^b_0$ denote the objective coefficient associated with branch-node deactivation. We define $c^b_{ij}=1, \forall i\in I,\ j\in J$, and $c^b_0=-1$. Thus, an active branch node contributes a positive coefficient to the objective through its selected split, whereas a deactivated branch node contributes a negative coefficient. Since the objective is minimized, this construction favors branch-node deactivation and therefore encourages simpler trees with fewer active branch nodes.

Finally, let $\alpha\geq 0$ and $\beta\geq 0$ denote weighting parameters controlling the relative importance of classification performance and tree complexity, respectively. The parameter $\alpha$ determines the weight placed on classification performance, while $\beta$ controls the preference for simpler trees.

\subsubsection{Routing and Compatibility Sets}

We next define the routing and compatibility sets used to construct the set-packing constraints of MSSF. The purpose of these sets is to identify combinations of branching and assignment decisions that cannot simultaneously occur in a feasible classification tree.

For a candidate split $(i,j)$, let $S^l_{ij}\subseteq I$ denote the set of observations that are infeasible for the left-routing condition induced by split value $x_{ij}$, and let $S^r_{ij}\subseteq I$ denote the set of observations that are infeasible for the corresponding right-routing condition. The left- and right-routing conventions are defined so that every observation is associated with exactly one side of the split. Consequently, $S^l_{ij}\cup S^r_{ij}=I$ and $S^l_{ij}\cap S^r_{ij}=\emptyset$. Thus, $S^l_{ij}$ identifies observations that cannot be assigned to leaves reached through the left child when split $(i,j)$ is selected, while $S^r_{ij}$ identifies observations that cannot be assigned to leaves reached through the right child.

Branch-node deactivation can be represented within the same routing framework. Because $b_{0v}=1$ routes every observation reaching node $v$ to its right child, we define $S^l_0=I$ and $S^r_0=\emptyset$. Thus, every observation is infeasible for a left descendant of a deactivated branch node, while deactivation creates no additional exclusion for the right descendant. This is exactly the routing effect of the lower artificial threshold $u_j^{(0)}=\min_{i\in I}x_{ij}-0.001$ for any feature $j$.

To treat deactivation and active splits consistently, let $\mathcal Q^a$ denote the set of admissible active split candidates, indexed by $(i,j)$, and let $\mathcal Q=\{0\}\cup\mathcal Q^a$. For $q=(i,j)\in\mathcal Q^a$, write $b_{qv}=b_{ijv}$, $c_q^b=c_{ij}^b$, and $S_q^d=S_{ij}^d$ for $d\in\{l,r\}$. For $q=0$, retain $b_{0v}$ and $c_0^b=-1$. The lower artificial thresholds are equivalent routing representations of $q=0$ and do not create additional branching alternatives. Define the dominance sets by $B_q^d=\{q'\in\mathcal Q\setminus\{q\}:S_q^d\subseteq S_{q'}^d\}$, $q\in\mathcal Q$,\ $d\in\{l,r\}$.

Thus, a dominance set includes all other alternatives whose routing-exclusion sets contain that of $q$, including alternatives with identical exclusion sets. This convention also handles repeated observed thresholds and any active candidate equivalent to deactivation. In particular, deactivation is included in $B_q^l$ for every active candidate with $S_q^l\subseteq I$. It is not included in the right-routing dominance set of any candidate with a nonempty right-exclusion set.

\subsubsection{Assignment-Index Collections}

We now construct the assignment-index collections used in the set-packing inequalities. Fix a leaf node $t\in T_L$ and a branching alternative $q\in\mathcal Q$. For $d\in\{l,r\}$, define the assignment-index ground set
\[
\bar C^d(t;q)=\{(u,k,t):u\in S_q^d,\ k\in K\}.
\]
When $S_q^d\neq\emptyset$, define
\[
\mathcal C^d(t;q)
=\left\{C\subseteq\bar C^d(t;q):
\forall k\in K,\ \exists!\,u\in S_q^d
\text{ such that }(u,k,t)\in C\right\}.
\]
When $S_q^d=\emptyset$, set $\mathcal C^d(t;q)=\emptyset$. For active alternatives $q=(i,j)$, this is the original construction $\mathcal C^d(t;i,j)$. In particular, $S_0^l=I$ and $S_0^r=\emptyset$ give $\bar C^l(t;0)=\{(u,k,t):u\in I,\ k\in K\}$ and $\mathcal C^r(t;0)=\emptyset$.

Each element of $\mathcal C^l(t;q)$ or $\mathcal C^r(t;q)$ contains exactly one assignment index for each class label, selected from observations incompatible with the corresponding routing direction. The same observation may appear with different class labels, but each class label appears exactly once in a set. In particular, the same observation may be selected for every class label, so the collections also capture individual-observation routing restrictions. For notational convenience, enumerate the elements as $C_p^d(t;q)$, $p=1,\ldots,|\mathcal C^d(t;q)|$. These collections provide the assignment-index sets for the unified MSSF inequalities.

\subsection{Proposed Formulation}
\label{subsec:mssf_formulation}

Using the notation and compatibility sets defined above, we now present the proposed MSSF. The formulation separates branching decisions from assignment decisions while enforcing their compatibility through set-packing inequalities. The set-partitioning constraints ensure that exactly one decision is selected at each branch node and that each observation is assigned exactly once.

The branching alternatives, dominance sets, and assignment-index collections are indexed by $q\in\mathcal Q$ as defined above. The deactivation alternative $q=0$ is included in the left-routing family. Since $\mathcal C^r(t;0)=\emptyset$, it generates no right-routing inequalities.

\begingroup
\allowdisplaybreaks
\begin{align}
\textnormal{(MSSF)}\quad \min\quad &
\alpha\sum_{i\in I}\sum_{k\in K}\sum_{t\in T_L}c^z_{ik}z_{ikt}
+\beta\sum_{v\in T_B}\sum_{q\in\mathcal Q}c_q^b b_{qv}
\label{eq:mssf_obj}
\\
\textnormal{s.t.}\quad &
\sum_{q\in\mathcal Q}b_{qv}=1,
\quad\forall v\in T_B,
\label{eq:mssf_branch_partition}
\\
&\sum_{k\in K}\sum_{t\in T_L}z_{ikt}=1,
\quad\forall i\in I,
\label{eq:mssf_assignment_partition}
\\
&b_{qv}+\sum_{q'\in B_q^l}b_{q'v}
+\sum_{(u,k,t)\in C_p^l(t;q)}z_{ukt}\leq1,
\notag\\[-0.2em]
&\qquad\forall q\in\mathcal Q,\ \forall t\in T_L,\
\forall v\in P_L(t),\
\forall p=1,\ldots,|\mathcal C^l(t;q)|,
\label{eq:mssf_left}
\\
&b_{qv}+\sum_{q'\in B_q^r}b_{q'v}
+\sum_{(u,k,t)\in C_p^r(t;q)}z_{ukt}\leq1,
\notag\\[-0.2em]
&\qquad\forall q\in\mathcal Q,\ \forall t\in T_L,\
\forall v\in P_R(t),\
\forall p=1,\ldots,|\mathcal C^r(t;q)|,
\label{eq:mssf_right}
\\
&z_{ikt}\in\{0,1\},
\quad\forall i\in I,\ \forall k\in K,\ \forall t\in T_L,
\label{eq:mssf_z_binary}
\\
&b_{qv}\in\{0,1\},
\quad\forall q\in\mathcal Q,\ \forall v\in T_B.
\label{eq:mssf_b_binary}
\end{align}
\endgroup

The objective~\eqref{eq:mssf_obj} weighs classification performance, through the coefficients $c^z_{ik}$ and the weight $\alpha$, against tree complexity, through the coefficients $c_q^b$ and the weight $\beta$. For a given classification performance, it therefore favors solutions with fewer active branch nodes.

Constraint~\eqref{eq:mssf_branch_partition} is a set-partitioning constraint over the branching decisions at each branch node. It requires exactly one branching decision at every node $v\in T_B$. Thus, either one candidate split $(i,j)$ is selected or the node is deactivated through $b_{0v}=1$. Constraint~\eqref{eq:mssf_assignment_partition} is a set-partitioning constraint over the assignment decisions and requires each observation to be assigned to exactly one class--leaf pair.

Constraints~\eqref{eq:mssf_left} and~\eqref{eq:mssf_right} enforce routing consistency for all branching alternatives, including deactivation. For an alternative $q$, the sets $S_q^l$ and $S_q^r$ contain the observations infeasible for the corresponding routing directions, and the collections $\mathcal C^l(t;q)$ and $\mathcal C^r(t;q)$ translate these exclusions into assignment-index sets. When $q$ is selected at a node on the path to leaf $t$, the inequalities prevent observations that violate the required routing condition from being assigned to $t$. Through the dominance sets $B_q^l$ and $B_q^r$, each inequality also excludes the assignments incompatible with every alternative whose exclusion set contains that of $q$, including identical sets. In particular, when $b_{0v}=1$ and $v\in P_L(t)$, the $q=0$ rows prevent any assignment to leaf $t$, because deactivation sends all observations to the right child; deactivation produces no right-routing rows and does not dominate any nonempty right-exclusion set.

The assignment-index collections also encode the class structure needed by the set-packing inequalities. Each element contains exactly one assignment index for each class label, and observations may repeat across labels. Iterating over all elements produces inequalities for individual-observation restrictions as well as mixed-observation class compatibility. Together with constraint~\eqref{eq:mssf_assignment_partition}, these enforce routing and leaf-label consistency for integer-feasible solutions.

The set-partitioning constraints~\eqref{eq:mssf_branch_partition} and~\eqref{eq:mssf_assignment_partition} contribute $|T_B|$ and $N_I$ constraints, respectively. The total number of unified set-packing inequalities is $\sum_{t\in T_L}\sum_{q\in\mathcal Q}\left(|P_L(t)|\,|\mathcal C^l(t;q)|+|P_R(t)|\,|\mathcal C^r(t;q)|\right)$. For a nonempty routing-exclusion set, the construction selects one observation index independently for each class label. Thus $|\mathcal C^d(t;q)|=|S_q^d|^{|K|}$ whenever $S_q^d\neq\emptyset$, and the corresponding collection is empty otherwise. The number of inequalities can therefore grow rapidly with the number of class labels, candidate splits, and tree depth, motivating the decomposition developed in the next section.

\subsection{Class-Boundary Candidate Split Reduction}
\label{subsec:class_boundary_reduction}

To improve the computational efficiency of the proposed ISCG-OCT framework, we develop a recursive class-boundary reduction heuristic that systematically reduces the number of candidate splitting thresholds considered during optimization. The motivation is that evaluating every distinct observed feature value as a potential split can introduce substantial computational redundancy, particularly for large datasets with continuous features. The proposed heuristic exploits the relationship between feature values and class labels to identify potentially informative split thresholds. Unlike conventional boundary identification based solely on the complete training set, our approach recursively partitions observations at identified class boundaries and reexamines the resulting subsets across all features. 

When fitting an optimal classification tree, only a finite number of split values need to be considered for each feature. For each feature $j\in J$, let $U_j=\left\{u_j^{(1)}<u_j^{(2)}<\cdots<u_j^{(M_j)}\right\}$ denote the ordered set of distinct observed values of feature $j$, where $M_j\leq N_I$. Under the routing convention that observations satisfying $x_{ij}\leq u_j^{(r)}$ are sent to the left child, define the admissible nontrivial split set of the unreduced observed-value formulation as $U_j^{\mathrm{split}}=\{u_j^{(1)},\ldots,u_j^{(M_j-1)}\}$. In addition, define two artificial boundary values $u_j^{(0)}=u_j^{(1)}-0.001$ and $u_j^{(M_j+1)}=u_j^{(M_j)}+0.001$. The lower artificial boundary is equivalent in routing to the single deactivation alternative $q=0$ (with coefficient $c_0^b=-1$), while the upper artificial boundary is retained to preserve the complete MSSF compatibility structure, including leaf-label consistency, after candidate split reduction. These artificial values serve structural purposes and are not additional nontrivial branching thresholds; in particular, the lower boundary does not introduce a separate deactivation variable for each feature. For each distinct observed value $u_j^{(r)}$, define the corresponding value group $G_j^{(r)}=\left\{i\in I:x_{ij}=u_j^{(r)}\right\}$. Observations belonging to the same value group cannot be separated by a threshold on feature $j$ and are therefore treated jointly when identifying class boundaries.

The proposed reduction is based on class-boundary split values identified recursively over subsets generated by class-boundary partitions. The central idea is that a class boundary detected under one feature ordering defines a potentially relevant partition of the observations in the current subset. The observations in each resulting child subset are then reordered by every feature, and new class boundaries are identified within those restricted subsets. Repeating this procedure for all features captures class boundaries that may not be visible in the ordering of the complete training set but become visible after earlier class-boundary partitions.

For any nonempty subset $\Omega\subseteq I$ and feature $j\in J$, let $U_j(\Omega)=\left\{u_{j,\Omega}^{(1)}<\cdots<u_{j,\Omega}^{(M_j(\Omega))}\right\}$ be the ordered distinct feature-$j$ values among observations in $\Omega$, and let $G_{j,\Omega}^{(r)}=\left\{i\in \Omega:x_{ij}=u_{j,\Omega}^{(r)}\right\}$ be the corresponding value groups. Two consecutive groups $G_{j,\Omega}^{(r)}$ and $G_{j,\Omega}^{(r+1)}$ define a class boundary in $\Omega$ if there exist observations $i\in G_{j,\Omega}^{(r)}$ and $i'\in G_{j,\Omega}^{(r+1)}$ such that $y_i\neq y_{i'}$. Whenever such a boundary occurs, both values immediately before and after the class transition, $u_{j,\Omega}^{(r)}$ and $u_{j,\Omega}^{(r+1)}$, are retained as class-boundary split values for feature $j$ in subset $\Omega$, provided they induce nontrivial splits of $\Omega$. Denote this set by $\mathcal Q_j(\Omega)=\left\{u_{j,\Omega}^{(q)}:q\in\{r,r+1\},\ r\in\mathcal B_j(\Omega),\ 1\leq q<M_j(\Omega)\right\}$, where $\mathcal B_j(\Omega)=\left\{r:G_{j,\Omega}^{(r)}\text{ and }G_{j,\Omega}^{(r+1)}\text{ define a class boundary in }\Omega\right\}$. The restriction $q<M_j(\Omega)$ excludes the largest observed value in $\Omega$, which would not induce a nontrivial split of that subset.

Each retained class-boundary split value also generates the two subsets that would result from splitting $\Omega$ at that value. For a retained value $u_{j,\Omega}^{(q)}\in\mathcal Q_j(\Omega)$, define $\Omega_{j,q}^{L}=\left\{i\in \Omega:x_{ij}\leq u_{j,\Omega}^{(q)}\right\}$, $\Omega_{j,q}^{R}=\Omega\setminus \Omega_{j,q}^{L}$. The same class-boundary search is then applied to each nonempty child subset over all features. Thus, after a boundary is identified using feature $j$, all observations in the resulting subset are retained when that subset is reordered by every feature $j'\in J$. This is important because the recursive search is not restricted to observations adjacent to the preceding boundary. It searches for new class boundaries among all observations that survive the preceding class-boundary partitions.

Let $\OmegaSet_0=\{I\}$. Recursively, for $d=0,1,\ldots$, define $\OmegaSet_{d+1}$ as the collection of all nonempty child subsets generated by retained class-boundary splits of every $\Omega\in\OmegaSet_d$ over all features:

\begin{equation}
\OmegaSet_{d+1}
=
\bigcup_{\Omega\in\OmegaSet_d}
\bigcup_{j\in J}
\bigcup_{q\in\mathcal C_j(\Omega)}
\left\{\Omega_{j,q}^{L},\Omega_{j,q}^{R}\right\},
\end{equation}
where
\begin{equation}
\mathcal C_j(\Omega)
=
\left\{
q\in\{1,\ldots,M_j(\Omega)-1\}:
q\in\mathcal B_j(\Omega)
\text{ or }
q-1\in\mathcal B_j(\Omega)
\right\}.
\end{equation}
Duplicate subsets are retained only once. The procedure is repeated over all features for every newly generated subset and continues until no previously unseen subset is generated. Since $I$ is finite, only finitely many nonempty subsets of $I$ exist, and therefore the recursion terminates after finitely many iterations. Let $\OmegaSet^{\star}=\bigcup_{d\geq 0}\OmegaSet_d$ denote the resulting collection of recursively generated subsets. The global recursive class-boundary split set for feature $j\in J$ is obtained by taking the unique union of all feature-$j$ boundary values identified over the recursively generated subsets. Accordingly, define $\mathcal R_j=\left[\left(\bigcup_{\Omega\in\OmegaSet^{\star}}\mathcal Q_j(\Omega)\right)\cap U_j^{\mathrm{split}}\right]\cup\left\{u_j^{(0)},u_j^{(M_j+1)}\right\}$.

Hence, an observed value of feature $j$ is retained when it lies immediately before or after a class transition in the feature-$j$ ordering of the complete training set or of any recursively generated subset, provided it induces a nontrivial split of that subset. A value that is not adjacent to a class boundary in the complete ordering may therefore be retained once preceding class-boundary partitions restrict the observations under consideration. The construction is independent of the maximum depth of the tree subsequently optimized and avoids retaining every observed value a priori. The sets $\mathcal R_j$ specify the globally retained candidate thresholds for the reduced MSSF formulation. In pricing, candidate generation must cover the feasible thresholds from these sets for each state; using class changes within a current state to discard further thresholds would require a separate justification to preserve exactness over the reduced formulation.
Let \(N_I=|I|\), \(N_J=|J|\), and \(N_{\Omega}=|\OmegaSet^{\star}|\) denote the number of unique nonempty observation subsets generated by the recursive class-boundary procedure. For each subset \(\Omega\in\OmegaSet^{\star}\), sorting the observations according to one feature requires \(O(|\Omega|\log|\Omega|)\) time, while scanning consecutive value groups to identify class boundaries requires \(O(|\Omega|)\) time. Since at most \(2(|\Omega|-1)\) class-boundary split values can be identified for each feature, and constructing the corresponding child subsets requires \(O(|\Omega|)\) time per retained split value, the total computational complexity of a direct implementation is \(O\!\left(N_J\sum_{\Omega\in\OmegaSet^{\star}}\left(|\Omega|\log|\Omega|+|\Omega|^2\right)\right)\). Using \(|\Omega|\leq N_I\), this expression is bounded by \(O(N_JN_{\Omega}N_I^2)\).

\section{Incumbent Subtree Column Generation}
\label{sec:column_generation}

Column generation is a decomposition framework for optimization problems whose reformulations contain a potentially very large number of decision patterns. Rather than enumerating all such patterns in advance, column generation maintains a restricted subset of columns in a RMP. The master problem determines how the currently available columns are combined, while a pricing problem searches for omitted columns that can improve the current solution. This decomposition is particularly attractive when the original formulation contains a large number of local decisions but these decisions can be organized into structured patterns that admit a more efficient pricing procedure.

The MSSF formulation developed in Section~\ref{sec:mssf} has precisely this structure. A feasible OCT solution is determined jointly by the branching decisions $b_{ijv}$ and $b_{0v}$ and the observation assignment decisions $z_{ikt}$. The exact-one constraints require one branching or deactivation decision at every branch node and one leaf--class assignment for every observation, while the MSSF compatibility constraints ensure that these decisions collectively describe a valid classification tree. Directly solving MSSF requires all of these decisions to be optimized simultaneously.

We therefore develop an incumbent subtree column generation decomposition in which each generated column represents a closed feasible subtree replacement. Starting from a current integer-feasible MSSF tree, a column replaces the decisions associated with a subtree while retaining the incumbent decisions outside that subtree. Closure requires the column to contain all branching and assignment decisions whose feasibility can be affected by the replacement. Branching and assignment decisions are therefore generated jointly within a single subtree column and satisfy the MSSF routing structure by construction. Columns associated with compatible nonoverlapping subtrees may be selected simultaneously in the RMP.

This construction yields a set-partitioning representation over structured subtree columns. Each decision owned by the decomposition is covered exactly once, either by retaining its incumbent realization or by selecting a compatible subtree replacement column. The RMP therefore determines which generated closed subtree patterns are combined to form the next integer-feasible MSSF tree, while the pricing problem searches the implicitly defined family of feasible subtree patterns for improving columns.

The proposed framework differs from conventional LP-based column generation in an important respect. Standard column generation typically solves an LP restricted master problem and prices new columns using reduced costs obtained from the dual multipliers of that master problem. In ISCG-OCT, the master is maintained over integer-feasible subtree patterns and pricing is recentered at the current integer-feasible tree. The coefficient of a candidate subtree pattern is therefore expressed directly as its objective change relative to the incumbent. If $p$ denotes a closed feasible subtree pattern, its pricing value can be written as $\Delta C_p=C_p-C_p^*$, where $C_p^*$ is the contribution of the incumbent decisions owned by the same subtree. A pattern satisfying $\Delta C_p<0$ is consequently an improving integer-feasible replacement.

Finally, when the priced subtree is rooted at the root of the complete tree, its closed-pattern family spans the complete integer-feasible MSSF tree space. Hence, exhaustive exact root pricing provides the global certificate $\min_{p\in\mathcal P_{\mathrm{root}}}\Delta C_p\ge 0$, which implies that no integer-feasible MSSF tree has an objective value strictly better than the current incumbent. Thus, although intermediate subtree pricing can be interpreted as structured large-neighborhood improvement, exact root pricing provides the connection between the column-generation representation and global integer optimality.

\subsection{Overview of the Decomposition}
\label{subsec:cg_overview}

Column generation maintains a restricted subset of decision patterns, or columns, in an RMP, while a pricing problem searches for omitted columns that can improve the current solution. It is attractive when a formulation contains many local decisions that can be organized into structured patterns admitting an efficient pricing procedure.

The MSSF formulation developed in Section~\ref{sec:mssf} has precisely this structure. A feasible OCT solution is determined jointly by the branching decisions $b_{ijv}$ and $b_{0v}$ and the observation assignment decisions $z_{ikt}$. The exact-one constraints require one branching or deactivation decision at every branch node and one leaf--class assignment for every observation, while the MSSF compatibility constraints ensure that these decisions collectively describe a valid classification tree. Directly solving MSSF requires all of these decisions to be optimized simultaneously.

We therefore develop an incumbent subtree column generation decomposition in which each generated column represents a closed feasible subtree replacement. Starting from a current integer-feasible MSSF tree, a column replaces the decisions associated with a subtree while retaining the incumbent decisions outside that subtree. Closure requires the column to contain all branching and assignment decisions whose feasibility can be affected by the replacement. Branching and assignment decisions are therefore generated jointly within a single subtree column and satisfy the MSSF routing structure by construction. Columns associated with compatible nonoverlapping subtrees may be selected simultaneously in the RMP.

This construction yields a set-partitioning representation over structured subtree columns. Each decision owned by the decomposition is covered exactly once, either by retaining its incumbent realization or by selecting a compatible subtree replacement column. The RMP therefore determines which generated closed subtree patterns are combined to form the next integer-feasible MSSF tree, while the pricing problem searches the implicitly defined family of feasible subtree patterns for improving columns.

The proposed framework differs from conventional LP-based column generation in an important respect. Standard column generation typically solves an LP restricted master problem and prices new columns using reduced costs obtained from the dual multipliers of that master problem. In ISCG-OCT, the master is maintained over integer-feasible subtree patterns and pricing is recentered at the current integer-feasible tree. The coefficient of a candidate subtree pattern is therefore expressed directly as its objective change relative to the incumbent. If $p$ denotes a closed feasible subtree pattern, its pricing value can be written as $\Delta C_p=C_p-C_p^*$, where $C_p^*$ is the contribution of the incumbent decisions owned by the same subtree. A pattern satisfying $\Delta C_p<0$ is consequently an improving integer-feasible replacement.

Finally, when the priced subtree is rooted at the root of the complete tree, its closed-pattern family spans the complete integer-feasible MSSF tree space. Hence, exhaustive exact root pricing provides the global certificate $\min_{p\in\mathcal P_{\mathrm{root}}}\Delta C_p\ge 0$, which implies that no integer-feasible MSSF tree has an objective value strictly better than the current incumbent. Thus, although intermediate subtree pricing can be interpreted as structured large-neighborhood improvement, exact root pricing provides the connection between the column-generation representation and global integer optimality.

\subsection{Incumbent Decisions and Subtree Patterns}
\label{subsec:incumbent_change_costs}

For every observation $i\in I$, let $(k_i^*,t_i^*)$ denote its unique incumbent assignment $z^*_{ik_i^*t_i^*}=1$. Define $c_i^{z,*}=c^z_{ik_i^*}$. Since $c^z_{ik}\in\{-1,1\}$, $c_i^{z,*}=-1$ if observation $ i$ is correctly classified, and $ 1$ otherwise. For every branch node $v\in T_B$, define the coefficient of its incumbent branching decision by $c_v^{b,*}=c_0^b b^*_{0v}+\sum_{i\in I}\sum_{j\in J}c^b_{ij}b^*_{ijv}$. By the exact-one branching constraint, $c_v^{b,*}=1$ if v is active in the incumbent tree, and $-1$ if v is deactivated. The objective-change coefficient of an assignment decision is $\Delta c^z_{ikt}=\alpha\left(c^z_{ik}-c_i^{z,*}\right),i\in I,\ k\in K,\ t\in T_L$. Thus,

\begin{equation}
\Delta c^z_{ikt}
=
\begin{cases}
0,
& \text{current correct, candidate correct},\\
2\alpha,
& \text{current correct, candidate incorrect},\\
-2\alpha,
& \text{current incorrect, candidate correct},\\
0,
& \text{current incorrect, candidate incorrect}.
\end{cases}
\label{eq:delta_z_cases}
\end{equation}

Similarly, for an active branching decision, $\Delta c^b_{ijv}=\beta\left(c^b_{ij}-c_v^{b,*}\right)$, and for deactivation, $\Delta c^b_{0v}=\beta\left(c^b_0-c_v^{b,*}\right)$. Because $c^b_{ij}=1$ and $c^b_0=-1$,
\begin{equation}
\Delta c^b_{ijv}
=
\begin{cases}
0, & v\text{ currently active},\\
2\beta, & v\text{ currently deactivated},
\end{cases}
\qquad
\Delta c^b_{0v}
=
\begin{cases}
-2\beta, & v\text{ currently active},\\
0, & v\text{ currently deactivated}.
\end{cases}
\label{eq:delta_b_cases}
\end{equation}

For $v\in T_B$, let $T(v)$ denote the nodes in the subtree rooted at $v$, and let $T_B(v)=T(v)\cap T_B$, $T_L(v)=T(v)\cap T_L$. Let $I_v^*\subseteq I$ denote the observations that reach node $v$ under the incumbent branching decisions on the path from the root to $v$. For the root node, $I_1^*=I$.

A subtree pattern rooted at $v$ is denoted by $p\in\mathcal P_v$. It contains $b^p_{ijw}, b^p_{0w}, w\in T_B(v)$, and $z^p_{ikt}, i\in I_v^*, t\in T_L(v)$. For each active leaf $t\in T_L(v)$, all observations assigned to leaf $t$ use a common class label $k_t^p\in K$. A pattern $p\in\mathcal P_v$ is admissible if it satisfies the MSSF branching partition equations for every $w\in T_B(v)$, $\sum_{q\in\mathcal Q}b^p_{qw}=1,\forall w\in T_B(v)$, the assignment partition equations $\sum_{k\in K}\sum_{t\in T_L(v)}z^p_{ikt}=1,\forall i\in I_v^*$, and all unified MSSF routing constraints \eqref{eq:mssf_left}--\eqref{eq:mssf_right} whose branching node and affected assignment decisions lie in the subtree rooted at $v$. The branching decisions above $v$ remain fixed at their incumbent values. Therefore, the boundary path by which observations in $I_v^*$ enter the subtree is unchanged. Closure additionally requires that the replacement and all retained incumbent decisions satisfy every unified MSSF compatibility inequality crossing the boundary of the replaced subtree.

For $p\in\mathcal P_v$, the branch-row incidence is $a^b_{wp}=1$ if $w\in T_B(v)$ and $0$ otherwise, and the assignment-row incidence is $a^z_{ip}=1$ if $i\in I_v^*$ and $0$ otherwise.

The objective-change coefficient of pattern $p$ is
\begin{align}
\Delta C_p
={}&
\sum_{i\in I_v^*}
\sum_{k\in K}
\sum_{t\in T_L(v)}
\Delta c^z_{ikt}z^p_{ikt}
+
\sum_{w\in T_B(v)}
\left(
\sum_{i\in I}\sum_{j\in J}
\Delta c^b_{ijw}b^p_{ijw}
+
\Delta c^b_{0w}b^p_{0w}
\right).
\label{eq:pattern_delta_cost}
\end{align}

\subsection{Restricted Master Problem}
\label{subsec:subtree_rmp}

Let $\bar{\mathcal P}$ be the set of generated closed subtree replacement patterns. Let $\lambda_p$ denote the master variable for pattern $p\in\bar{\mathcal P}$. Let $\lambda_v^0$ denote the incumbent keep variable for branch row $v$, and let $\lambda_i^0$ denote the incumbent keep variable for assignment row $i$. The incumbent-centered restricted master problem is

\begin{align}
\textnormal{(RMP)}\quad
\min\quad&
\sum_{p\in\bar{\mathcal P}}\Delta C_p\lambda_p
\label{eq:subtree_rmp_obj}
\\
\textnormal{s.t.}\quad
&
\lambda_v^0
+
\sum_{p\in\bar{\mathcal P}}
a^b_{vp}\lambda_p
=
1,
&&
\forall v\in T_B,
\label{eq:subtree_rmp_branch}
\\
&
\lambda_i^0
+
\sum_{p\in\bar{\mathcal P}}
a^z_{ip}\lambda_p
=
1,
&&
\forall i\in I,
\label{eq:subtree_rmp_assignment}
\\
&
\lambda_v^0,\ \lambda_i^0,\ \lambda_p\ge0,
&&
\forall v\in T_B,\ i\in I,\ p\in\bar{\mathcal P}.
\label{eq:subtree_rmp_nonneg}
\end{align}

For the integer master, the variables are restricted to
\begin{equation}
\lambda_v^0,\lambda_i^0,\lambda_p\in\{0,1\}.
\label{eq:subtree_integer_master}
\end{equation}

All master constraints therefore have right-hand side one and no slack variables are required. A replacement pattern covering a branch node or observation excludes the corresponding keep column and excludes every other replacement pattern covering the same row.

\begin{proposition}[Feasible Assembly of Closed Subtree Patterns]
\label{prop:feasible_assembly}
Let $\lambda$ be an integer-feasible solution of \eqref{eq:subtree_rmp_branch}--\eqref{eq:subtree_integer_master}. Suppose every generated pattern is closed and internally MSSF-feasible, and that closure ensures compatibility with all retained incumbent decisions across the boundaries of the selected patterns. Then the selected patterns, together with the incumbent decisions retained by the keep columns, define an integer-feasible MSSF solution.
\end{proposition}

The objective-change formulation provides the desired pricing structure when the master is recentered at the current incumbent. At the beginning of a pricing round, consider the initial restricted master containing only
the keep columns:
\begin{align}
\min\quad &0
\notag\\
\textnormal{s.t.}\quad
&\lambda_v^0=1,
&&\forall v\in T_B,
\label{eq:keep_only_b}
\\
&\lambda_i^0=1,
&&\forall i\in I,
\label{eq:keep_only_z}
\\
&\lambda_v^0,\lambda_i^0\geq0.
\notag
\end{align}

Let $\pi_v$ be the dual variable of~\eqref{eq:keep_only_b} and let $\sigma_i$ be the dual variable of~\eqref{eq:keep_only_z}. Its dual maximizes $\sum_{v\in T_B}\pi_v+\sum_{i\in I}\sigma_i$ subject to $\pi_v\le0$ for all $v\in T_B$ and $\sigma_i\le0$ for all $i\in I$.

\begin{proposition}[Zero Dual Prices of the Re-Centered Master]
\label{prop:zero_initial_duals}
The unique optimal dual solution of the keep-only re-centered master is $\pi_v^*=0,\forall v\in T_B$, and $\sigma_i^*=0,\forall i\in I$. Consequently, the reduced cost of any candidate subtree pattern $p$ at the beginning of the pricing round is exactly $\bar C_p=\Delta C_p$.

\end{proposition}

Proposition~\ref{prop:zero_initial_duals} is a property of the recentered keep-only master. If generated replacement columns are retained and the LP master is reoptimized in the conventional column-generation manner, the master dual variables need not remain zero. Accordingly, the identity $\bar C_p=\Delta C_p$ is guaranteed at each pricing round by recentering the master at the accepted incumbent. This distinguishes the proposed incumbent-centered integer column-generation scheme from conventional LP column generation with a permanently accumulating RMP.

\subsection{Subtree Pricing Problem}
\label{subsec:subtree_pricing}

For a candidate replacement root $v\in T_B$, pricing searches only over the subtree rooted at $v$. The branching decisions above $v$ are fixed at their incumbent values, and the observations entering the subtree are fixed to $I_v^*$. Thus, the pricing problem does not optimize a complete tree unless $v$ is the root.

Let $\mathcal F_v$ denote the set of binary subtree decisions satisfying the exact-one branching and assignment conditions and the restrictions of the unified MSSF routing constraints \eqref{eq:mssf_left}--\eqref{eq:mssf_right} to the subtree rooted at $v$, together with the closure conditions, including compatibility with retained incumbent decisions across the subtree boundary, defining an admissible subtree
pattern. The pricing problem associated with root $v$ is based on this feasible set.

\begin{align}
\textnormal{(PP$_v$)}\qquad
\min\quad&
\sum_{i\in I_v^*}
\sum_{k\in K}
\sum_{t\in T_L(v)}
\Delta c^z_{ikt}z_{ikt}
\notag\\
&+
\sum_{w\in T_B(v)}
\left(
\sum_{i\in I}\sum_{j\in J}
\Delta c^b_{ijw}b_{ijw}
+
\Delta c^b_{0w}b_{0w}
\right)
\label{eq:subtree_pricing_obj}
\\
\textnormal{s.t.}\quad
&(b,z)\in\mathcal F_v.
\label{eq:subtree_pricing_feasible}
\end{align}

By Proposition~\ref{prop:zero_initial_duals}, the objective of PP$_v$ is exactly the reduced cost of the generated subtree column in the recentered master. The pricing problem has a direct feasibility and improvement interpretation. Let $p_v^*$ denote an optimal solution of $\mathrm{PP}_v$ with objective value $\Delta C_{p_v^*}$. By the definition of $\mathcal F_v$, $p_v^*$ is a closed MSSF-feasible subtree replacement. Replacing the incumbent subtree by $p_v^*$ while retaining the incumbent decisions outside the covered rows therefore produces an integer-feasible MSSF tree by Proposition~\ref{prop:feasible_assembly}. 

Moreover, $\Delta C_{p_v^*}$ is exactly the change in the original MSSF objective produced by this replacement. Hence, if $\Delta C_{p_v^*}<0$, the replacement yields a strictly improving feasible tree. Conversely, if $\Delta C_{p_v^*}\geq0$, then no closed feasible replacement of the subtree rooted at $v$ can improve the incumbent objective while the branching decisions on the path above $v$ remain fixed, since $(\mathrm{PP}_v)$ minimizes the exact objective change over all such replacements.

Let $\mathcal P=\bigcup_{v\in T_B}\mathcal P_v$ denote the complete family of admissible closed subtree patterns. When the root is included as an admissible replacement root, $\mathcal P_1$ contains a pattern corresponding to every integer-feasible MSSF tree. Indeed, the subtree rooted at the root contains the entire fixed tree topology and $I_1^*=I$. Hence, an admissible root pattern specifies all branching, deactivation, and assignment decisions and satisfies the complete MSSF feasibility structure. Conversely, every integer-feasible MSSF solution defines an admissible pattern in $\mathcal P_1$.

The recursive pricing procedure is embedded within an incumbent-centered integer column generation framework. Starting from the current incumbent tree, recursive pricing generates one or more closed subtree replacement patterns together with their exact objective changes. These patterns are added to the restricted master problem, which is solved as an integer optimization problem to determine the best compatible collection of subtree replacements. If the resulting solution improves the incumbent objective, the selected patterns are assembled into a new integer-feasible tree, which becomes the incumbent for the next iteration. The objective change coefficients, boundary observation sets, and keep columns are then redefined relative to the new incumbent, and the pricing procedure is repeated until no improving replacement pattern can be generated.

\section{Dynamic Programming for Recursive Pricing}
\label{sec:dynamic_pricing}

This section develops a DP framework for solving the subtree-pricing problems recursively and relates the exact Bellman recursion to the lower-bound-guided search used in ISCG-OCT. The central observation is that, once the decision at a branch node is fixed, the observations reaching that node are routed to its descendant subtrees according to the MSSF routing rules. For an active split, the observations are partitioned between the left and right child subtrees. For a deactivated node, the left child receives no observations and all observations continue to the right child. Because the pricing objective is additive over the branching and assignment decisions contained in a closed subtree, the continuation value below a node can be written as the local branching contribution plus the values of its two descendant subproblems.

The DP uses the objective-change coefficients derived in the previous section. Hence, every state and transition is evaluated on the same incumbent-centered objective scale as the pricing problem. A negative value represents an improvement relative to the corresponding incumbent decisions, whereas a nonnegative value indicates that the represented replacement does not improve the incumbent. In addition to the exact continuation value, the recursive pricing procedure associates each unresolved state with the lower bound $LB_v(I_v)$ developed in Appendix~\ref{subsec:lower_bound_computation}. This lower bound is used to screen partial states before their complete descendant continuations are evaluated.

For every branch node $v\in T_B$, let $\ell(v)$ and $r(v)$ denote its left and right children. A pricing state is represented by $(v,I_v)$, where $I_v\subseteq I$ is the set of observations reaching node $v$ under the decisions already fixed on the path from the root of the priced subtree to $v$. At the root of a subtree replacement, the initial observation set is the boundary set $I_v^*$. Once the decision at $v$ is selected, the MSSF routing structure uniquely determines the observation sets entering its two children.

For an active split represented by $(i,j)$ at node $v$, define $I^l_{ijv}(I_v)$ and $I^r_{ijv}(I_v)$ as the observations in $I_v$ routed to the left and right child, respectively. These sets satisfy $I^l_{ijv}(I_v)\cap I^r_{ijv}(I_v)=\emptyset$ and $I^l_{ijv}(I_v)\cup I^r_{ijv}(I_v)=I_v$. Let $\mathcal B_v(I_v)$ denote the set of split candidates $(i,j)$ that are feasible for state $(v,I_v)$ under the restrictions defining the subtree-pricing feasible set $\mathcal F_v$. In particular, the induced child observation sets must satisfy the MSSF routing restrictions associated with the selected split. No additional feasibility condition is imposed by the DP beyond those defining $\mathcal F_v$. The admissible splits are drawn from predefined candidate threshold sets $\{\mathcal R_j\}_{j\in J}$, which may be constructed by any prescribed procedure. The DP is exact relative to the resulting candidate-restricted feasible family, provided all its admissible alternatives are considered.

Let $V_v(I_v)$ denote the minimum objective change achievable by a closed feasible completion of the subtree rooted at $v$, according to $\mathcal F_v$, when the observations entering $v$ are $I_v$. Thus, $V_v(I_v)$ is the exact continuation value associated with the subtree-pricing problem. At a leaf $t\in T_L$, no further branching decision is required. The remaining decision is the common class assigned to the observations reaching that leaf. For $I_t\subseteq I$, define $V_t(I_t)=\min_{k\in K}\sum_{i\in I_t}\alpha\left(c^z_{ik}-c_i^{z,*}\right)$.

The minimizing class determines the assignment variables of the generated pattern. Hence, $V_t(I_t)$ is the exact change in classification cost obtained by assigning the observations reaching $t$ to their best common class relative to their incumbent classification contributions. For an empty leaf state, $V_t(\emptyset)=0$.

A branch node may alternatively be deactivated by selecting $b_{0v}=1$, which is routing-equivalent to selecting the lower artificial threshold introduced in Section~\ref{subsec:class_boundary_reduction}. Under the MSSF deactivation convention, deactivating node $v$ does not terminate the subtree rooted at $v$. Instead, no observation is routed to the left child and every observation in $I_v$ continues to the right child. Because a closed pattern rooted at $v$ contains the descendant branching and assignment decisions, these decisions must also be included in the value of the deactivation alternative. Accordingly, define $V_v^0(I_v)=\Delta c^b_{0v}+V_{\ell(v)}(\emptyset)+V_{r(v)}(I_v)$.

The first term is the local objective change associated with deactivating $v$, while the remaining terms account for the optimal closed completions of the left and right descendant subtrees. More generally, an empty branch state still requires a feasible closed completion of the descendant branch rows. Therefore, $V_v(\emptyset)$ is evaluated by the same Bellman recursion as any other branch state, with all induced child observation sets empty. In particular, it is not assumed to be zero unless $v$ is a leaf. For an admissible active split $(i,j)\in\mathcal B_v(I_v)$, define $Q_{ijv}(I_v)=\Delta c^b_{ijv}+V_{\ell(v)}\left(I^l_{ijv}(I_v)\right)+V_{r(v)}\left(I^r_{ijv}(I_v)\right)$.

The first term is the local branching contribution, and the remaining terms are the optimal objective changes of the two descendant subproblems induced by the selected split. For every branch state $(v,I_v)$, the exact value function is $V_v(I_v)=\min\left\{V_v^0(I_v),\;\min_{(i,j)\in\mathcal B_v(I_v)}Q_{ijv}(I_v)\right\}$. Equivalently,
\begin{equation}
V_v(I_v)
=
\min
\left\{
\Delta c^b_{0v}
+
V_{\ell(v)}(\emptyset)
+
V_{r(v)}(I_v),
\;
\min_{(i,j)\in\mathcal B_v(I_v)}
\left[
\Delta c^b_{ijv}
+
V_{\ell(v)}\left(I^l_{ijv}(I_v)\right)
+
V_{r(v)}\left(I^r_{ijv}(I_v)\right)
\right]
\right\}.
\label{eq:dp_bellman_expanded}
\end{equation}
Thus, the two branch alternatives differ only in the observation sets passed to the descendants. The exact recursion is evaluated from the leaves toward the root of the priced subtree. For each exactly solved state, the DP stores both its minimum continuation value and a minimizing branching or deactivation decision. Following these stored decisions recursively recovers the minimum-cost complete subtree pattern.

The Bellman recursion above identifies the exact continuation value when every feasible decision and reachable state is evaluated. The implemented recursive pricing procedure uses the lower bound $LB_v(I_v)$ to reduce this search. While a subtree state remains unresolved, its feasible branching and deactivation extensions are considered and their lower bounds are used to screen extensions before evaluating their complete descendant continuations. After a decision is fixed, the resulting child extension defines a partially fixed state. When applied to such an extension, $LB_v(I_v)$ denotes a lower bound conditional on the decisions already fixed in that extension. Its lower bound combines the objective contribution of decisions already fixed with the assignment lower bound under the remaining routing flexibility. Therefore, for any closed feasible completion $p$ extending that partial state, $LB_v(I_v)\leq\Delta C_p$, where the notation includes the contribution of the fixed decisions when compared with the global cutoff. 

A partial extension is eligible for further expansion only when its lower bound satisfies $LB_v(I_v)<UB-\Gamma(\rho)$, where $UB$ is the current upper bound on the objective change and $\Gamma(\rho)=\rho N_I\alpha$ for an improvement-threshold parameter $\rho\ge0$. Once a state $(v,I_v)$ is solved, its complete subtree is represented by the minimum-cost feasible pattern attaining $V_v(I_v)$. The state and its minimizing subtree decisions are stored and need not be recomputed. Screening a partial extension does not change the Bellman minimum of any state that is solved exactly. Only values whose minimum over all admissible continuations has been established are stored as exact; cutoff-dependent results are not reused as exact Bellman values. Screening against a global cutoff includes the contributions of fixed decisions outside the local continuation. Exact memoized values are cleared when the incumbent-relative coefficients change. A positive improvement threshold may exclude feasible completions that do not meet the cutoff; termination of thresholded pricing alone therefore does not certify a zero global optimality gap.

A DP state is completely characterized by the current node and the set of observations reaching that node, provided the feasible decisions and objective-change coefficients depend only on this pair. States having the same pair $(v,I_v)$ have identical continuation problems regardless of the sequence of ancestor decisions by which they were reached. Consequently, their exact values and minimizing subtree patterns need to be computed only once. Let $\mathcal S=\left\{(v,I_v):v\in T_B\cup T_L,\;I_v\subseteq I\right\}$ denote the set of reachable pricing states. For each solved state, the DP stores its optimal value and the branching or deactivation decision attaining that value, together with pointers to the minimizing child states. Repeated solved states are retrieved from memory rather than recomputed. In the lower-bound-guided implementation, unresolved partial extensions are screened before evaluation, while an exactly solved subtree is replaced by its stored minimum-cost complete pattern whenever the same state is encountered.

Let $B=\max_{(v,I_v)}|\mathcal B_v(I_v)|$ denote the maximum number of admissible active split candidates at any reachable branch state. Including deactivation, each branch state has at most $B+1$ admissible decisions. Let $\mathcal S_B$ and $\mathcal S_L$ denote the distinct branch and leaf states actually solved, respectively. With observation subsets represented so that routing an active split costs at most $O(N_I)$, evaluating each branch state requires at most $O((B+1)N_I)$ operations, apart from retrieving its already computed child values. Evaluating each leaf state by comparing all classes costs at most $O(N_KN_I)$. Since each node can be reached with at most $2^{N_I}$ distinct observation subsets, $|\mathcal S_B|\leq |T_B(v)|2^{N_I}$ and $|\mathcal S_L|\leq |T_L(v)|2^{N_I}$. Consequently, a valid worst-case time bound for exact memoized DP is $O\!\left(2^{N_I}N_I\left[(B+1)|T_B(v)|+N_K|T_L(v)|\right]\right)$.

\section{ISCG-OCT}
\label{subsec:iscg}

ISCG-OCT addresses the OCT problem by combining incumbent subtree column generation with a recursive DP search. The incumbent tree provides a feasible reference configuration and induces the objective-change costs used in pricing. The dynamic program then constructs a candidate subtree recursively, one branch-node decision at a time. Each DP state is characterized by the current node and the observations reaching that node, while the partially fixed decisions determine the routing restrictions used to compute the associated lower bound. The lower bound associated with a partially fixed state provides a lower bound for all complete subtree columns that can be generated from that state. Consequently, states that cannot produce a sufficiently improving column can be discarded before the remaining subtree is enumerated.

The initial incumbent tree is obtained from CART. Let $\bar p$ denote the subtree configuration inherited from this incumbent. Pricing is expressed in terms of objective change relative to $\bar p$. Thus, keeping the incumbent subtree unchanged has value $\Delta C(\bar p)=0$. Since $\bar p$ is feasible, it immediately supplies the baseline upper bound $UB=0$. Accordingly, a newly generated subtree column is improving only if its objective-change value is negative. This incumbent-centered representation is the connection between the column-generation interpretation and the recursive pricing procedure. The DP does not search for an arbitrary classification tree, but for a subtree configuration that improves upon the current incumbent configuration.

If another complete feasible subtree with value $\widehat{UB}<0$ is available, then $\widehat{UB}$ is a tighter upper bound and may replace $UB=0$. More generally, any valid feasible completion can be used to tighten the upper bound. We use $UB=0$ in the baseline implementation because it is available immediately from the incumbent CART solution.

As in Section~\ref{sec:dynamic_pricing}, the state $(v,I_v)$ contains exactly the information required to continue the recursion: fixing the decision at $v$ decomposes the remaining pricing problem into its two descendant states, and the candidate decisions at a descendant are updated as earlier decisions change the observations reaching it.

In the implemented recursive search, the assignment component of this lower bound is computed by the hierarchical assignment lower-bound procedure. After the current branch-node decision is fixed, its routing restrictions determine the descendant leaves that remain feasible for each observation. The appendix procedure optimizes the remaining observation--leaf--class assignments over these feasible-leaf sets while retaining the common-class requirement at each leaf. Under the Subtree Feasible-Leaf Property, the value returned by the procedure is the optimal assignment objective value under the current partial routing restrictions. Proposition~\ref{prop:exact_assignment} establishes this exactness result. The resulting assignment value is then used in the computation of $LB_v(I_v)$. Consequently, fixing one branch-node decision produces both the descendant DP states and the assignment value used to construct the lower bound for all complete subtree patterns extending that partial fixing.

\begin{proposition}[Exact Assignment Evaluation]
\label{prop:exact_assignment}
Under the Subtree Feasible-Leaf Property established in Lemma~\ref{lem:sp2_subtree_feasible_leaf}, Algorithm~\ref{alg:assignment_lower_bound} computes $L_v^z=\min_{\emptyset\neq C\subseteq K}H_v(C)$, which is the exact optimal assignment value for a partially fixed pricing state $(v,I_v)$.
\end{proposition}

This lower bound calculation prevents the recursive pricing procedure from enumerating every complete subtree pattern. If the lower bound associated with a partially fixed state cannot satisfy the required improvement criterion, none of its completions needs to be generated. The particular lower bound described above is the one used in the baseline implementation. If a tighter valid lower bound is available for the same partially fixed state, it can be used without changing the recursive pricing logic.

With the baseline upper bound $UB=0$, the expansion rule of Section~\ref{sec:dynamic_pricing} becomes $LB_v(I_v)<-\rho N_I\alpha$, and the recursion retains all qualifying extensions of each partially fixed parent state. For $\rho>0$, a state is expanded only if its lower bound indicates a possible improvement of at least $\rho N_I\alpha$, so increasing $\rho$ reduces the number of expanded states but does not necessarily preserve the exact pricing optimum. Global optimality certification therefore requires valid bounds from exact pricing or an equivalent certified procedure, together with the applicable master-problem optimality conditions.

\section{Computational Results}
\label{cresults}

In this section, we evaluate the proposed ISCG-OCT framework on 27 publicly available benchmark datasets from the UCI Machine Learning Repository. The proposed method is compared with several state-of-the-art optimal decision tree algorithms: the MIP-based formulations OCT \citep{bertsimas2017optimal}, Benders OCT \citep{aghaei2025strong}, and BNP-OCT \citep{liu2024optimal}, and the dynamic-programming method STreeD \citep{van2023necessary}. The additional numerical results are provided in Appendix~\ref{computational}.

\subsection{Standard Data Sets}

All datasets are scaled to the interval $[0,1]$ prior to model training. Table~\ref{table:data-sets-info} summarizes the characteristics of the 27 benchmark datasets used in the computational study. The datasets range from 100 to 12,960 observations, contain between 4 and 72 features, and include 2 to 4 class labels. The CART benchmark is implemented using the \texttt{scikit-learn} library, while all mathematical programming models are solved with Gurobi~12.0.3. The experiments are conducted on a computer equipped with a 2.9~GHz Intel Core i7-8650U processor and 16~GB of RAM. The proposed ISCG-OCT algorithm is implemented in \texttt{C++} on Windows using the MSYS2 MinGW-w64 toolchain (GCC~15.2.0) and compiled under the \texttt{C++17} standard. Reported times include preprocessing and model construction.

STreeD is run with the same weights and time limit as the other methods. Because STreeD is exact for a given binarization of the features, every distinct training value of a continuous feature is supplied as a candidate threshold and every category of a categorical feature as a binary indicator, so a tree that STreeD certifies is optimal over all univariate splits. Its times are measured from the binarization of the features through the end of the search. The STreeD runs were performed on a different 70\%/30\% partition of each dataset. Their results are therefore reported for reference. Appendix~\ref{computational} gives further details.

The CART solution is used as the warm start for ISCG-OCT. We report two variants, denoted by ISCG $U=0\%$ and ISCG $U=5\%$, which use the improvement-threshold parameters $\rho=0$ (exact pricing) and $\rho=0.05$, respectively. $U$ is only a reporting label, and both variants use the upper bound $UB=0$. Under $\rho=0.05$, a partially fixed pricing state is expanded only if $LB_v(I_v)<-0.05N_I\alpha$, which reduces the number of states explored and can decrease computation time.

\begin{table}[!htbp]
\centering
\caption{Benchmark Datasets}
\label{table:data-sets-info}
\begin{adjustbox}{width=0.8\columnwidth}
\begin{tabular*}{\textwidth}{@{\extracolsep{\fill}}lccc}
\toprule
Dataset & Observations & Features & Classes \\
\midrule
Fertility & 100 & 9 & 2 \\
Thoracic Surgery& 470 &16 & 2 \\
Echocardiogram & 132 & 9 & 2 \\
Iris & 150 & 4 & 3 \\
Hepatitis & 155 & 19 & 2 \\
Planning Relax & 182 & 12 & 2 \\
Wine & 178 & 13 & 3 \\
House-vote & 232 & 16 & 2 \\
Seeds & 210 & 7 & 3 \\
Connectionist Bench (Sonar) & 208 & 60 & 2 \\
Ionosphere & 351 & 34 & 2 \\
Cylinder Bands & 512 & 39 & 2 \\
Breast Cancer Wisconsin (Diagnostic) & 569 & 30 & 2 \\
Balance-scale & 625 & 4 & 3 \\
Credit Approval & 690 & 15 & 2 \\
Blood Transfusion Service Center & 748 & 4 & 2 \\
Tic-tac-toe & 958 & 18 & 2 \\
German & 1000 & 20 & 2 \\
Biodeg & 1055 & 41 & 2 \\
Banknote Authentication & 1372 & 4 & 2 \\
Car-evaluation & 1728 & 6 & 4 \\
Ozone Level Detection & 2536 & 72 & 2 \\
Seismic-bumps & 2584 & 18 & 2 \\
Chess & 3196 & 36 & 2 \\
Spambase & 4601 & 57 & 2 \\
Mushroom & 5644 & 22 & 2 \\
Nursery & 12960 & 8 & 3 \\
\bottomrule
\end{tabular*}
\end{adjustbox}
\end{table}

Model performance is evaluated using training and testing accuracy. For each dataset, the observations are randomly partitioned once into 70\% training and 30\% testing sets. The maximum tree depth is varied from $D=2$ to $D=5$ and a one-hour time limit is imposed for all optimization runs. To maintain a common weighting scheme across methods, the objective weights are fixed at $\alpha=0.99$ and $\beta=0.01$, placing substantially greater weight on classification performance than on tree complexity. These weights are applied where the competing formulations contain corresponding loss and complexity terms. Differences in objective scaling and normalization may affect the comparability of reported objective values and relative optimality gaps.

Following the benchmarking methodology of \citet{dolan2002benchmarking}, we use performance profiles to compare the relative computational performance of competing algorithms over the collection of benchmark instances. For problem $p\in P$ and solver $s\in S$, the performance ratio is $r_{p,s}=t_{p,s}/\min_{s'\in S}t_{p,s'}$ for measures in which smaller values are preferred, such as computational times, where $t_{p,s}$ is the wall-clock time required by solver $s$ to terminate on problem $p$ under its specified stopping criterion, and $r_{p,s}=\max_{s'\in S}a_{p,s'}/a_{p,s}$ for measures in which larger values are preferred, such as classification accuracy, where $a_{p,s}$ is the classification accuracy achieved by solver $s$ on problem $p$. Under both definitions, the best-performing solver has a performance ratio of one. The computational-time profiles compare wall-clock times to termination, which may result from optimality certification, threshold-based stopping, or reaching the time limit, rather than times to a common optimality certificate. Runs that terminate without a feasible tree, such as STreeD runs that exceed the memory limit, are treated as failures and are not counted at any $\tau$. For a performance factor $\tau\ge1$, the performance profile $\phi_s(\tau)=|\{p\in P:r_{p,s}\le\tau\}|/|P|$ is the fraction of benchmark instances on which solver $s$ performs within a factor $\tau$ of the best-performing solver. A curve that rises more rapidly and remains higher indicates superior overall performance across the benchmark set.

\subsection{Comparison Results on Model Accuracy}

For tree depth $D=2$, the methods exhibit noticeable differences in the frequency with which they attain the best reported training accuracy, particularly at $\tau=1$. Exact ISCG-OCT attains the best reported accuracy, including ties, on approximately $85\%$ of the datasets, compared with approximately $78\%$ for Benders OCT, $70\%$ for BNP-OCT, $67\%$ for STreeD, and $41\%$ for OCT. The thresholded ISCG-OCT variant attains the best reported accuracy on approximately $33\%$ of the datasets at $\tau=1$. These percentages measure the frequency with which a method attains the best reported accuracy and should not be interpreted as the magnitude of the accuracy difference between methods. In particular, a method is not counted at a given value of $\tau$ even when its training accuracy differs only slightly from the best-performing method.

The lower best-accuracy frequency of the thresholded variant may be associated with its early stopping. If no subtree replacement improves the objective by more than the prescribed threshold, the algorithm returns the current incumbent, which may be the CART warm start or a previously identified feasible tree. The observed differences may, however, also depend on the pricing search and other termination conditions. As $\tau$ increases, the profiles also capture near-best performance. At $\tau=1.05$, exact ISCG-OCT and BNP-OCT cover all datasets, compared with approximately $81\%$ for ISCG $U=5\%$ and STreeD.

For tree depth $D=3$, STreeD attains the best reported training accuracy on approximately $74\%$ of the datasets at $\tau=1$, followed by exact ISCG-OCT with $52\%$, Benders OCT with $37\%$, and BNP-OCT, OCT, and ISCG $U=5\%$ with $30\%$ each. STreeD certifies optimality on 21 datasets at this depth, compared with 17 for exact ISCG-OCT, and its trees attain a higher training accuracy than those of exact ISCG-OCT on 13 of the 23 datasets on which it returns a tree, including several datasets for which ISCG-OCT also reports a zero gap. ISCG-OCT certifies optimality over its candidate-restricted MSSF feasible region on its training partition, whereas STreeD certifies optimality over all thresholds of the binarized features on a different partition. For $\tau\ge1.05$, however, exact ISCG-OCT has the highest profile (approximately $93\%$ at $\tau=1.05$ and $100\%$ at $\tau=1.1$), and both ISCG-OCT variants remain competitive with the MIP-based formulations.

\begin{figure}[htb]
    \centering
    \begin{subfigure}[t]{0.48\textwidth}
        \centering
        \includegraphics[width=\linewidth]{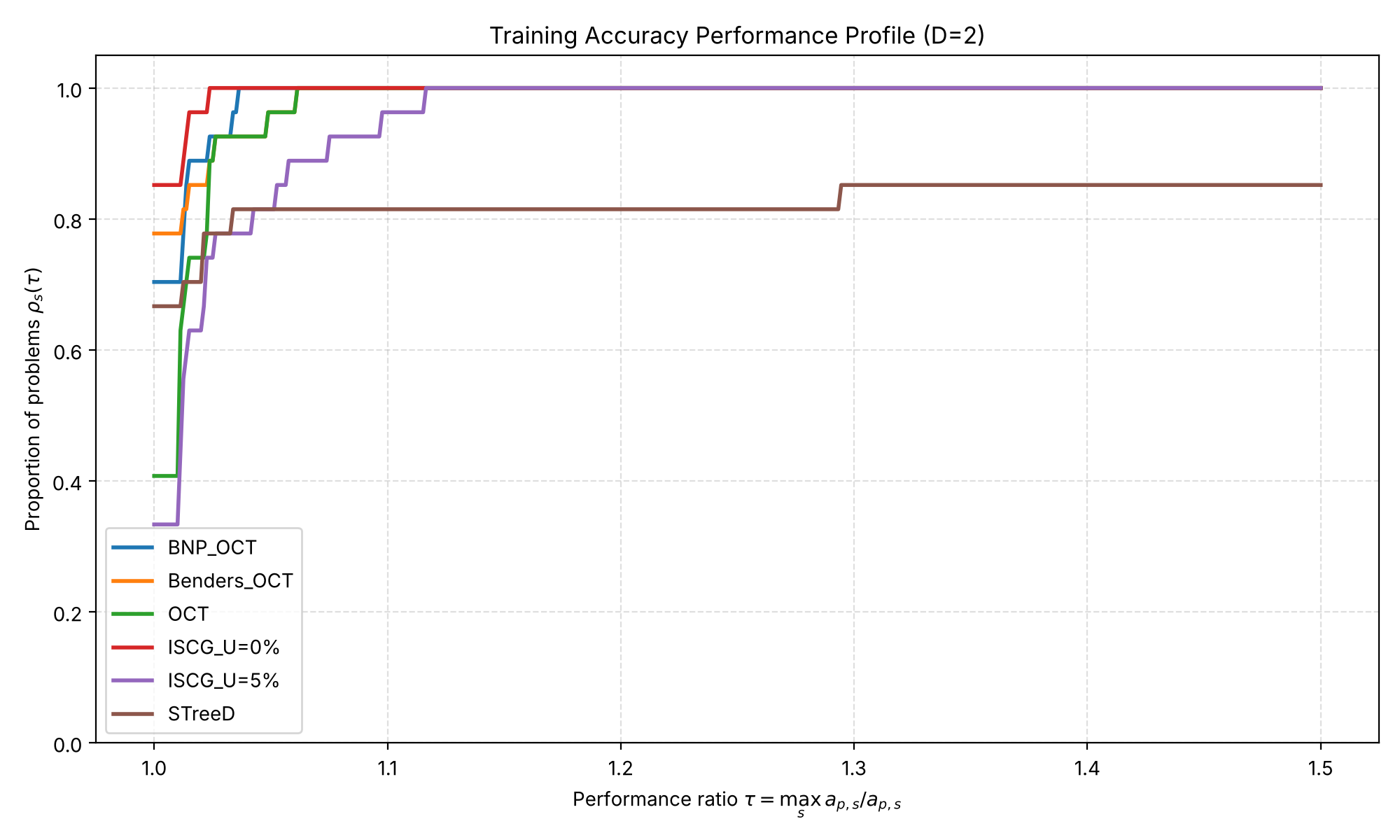}
        \caption{Training Accuracy Profile for Tree Depth 2}
        \label{fig:train_depth2}
    \end{subfigure}
    \hfill
    \begin{subfigure}[t]{0.48\textwidth}
        \centering
        \includegraphics[width=\linewidth]{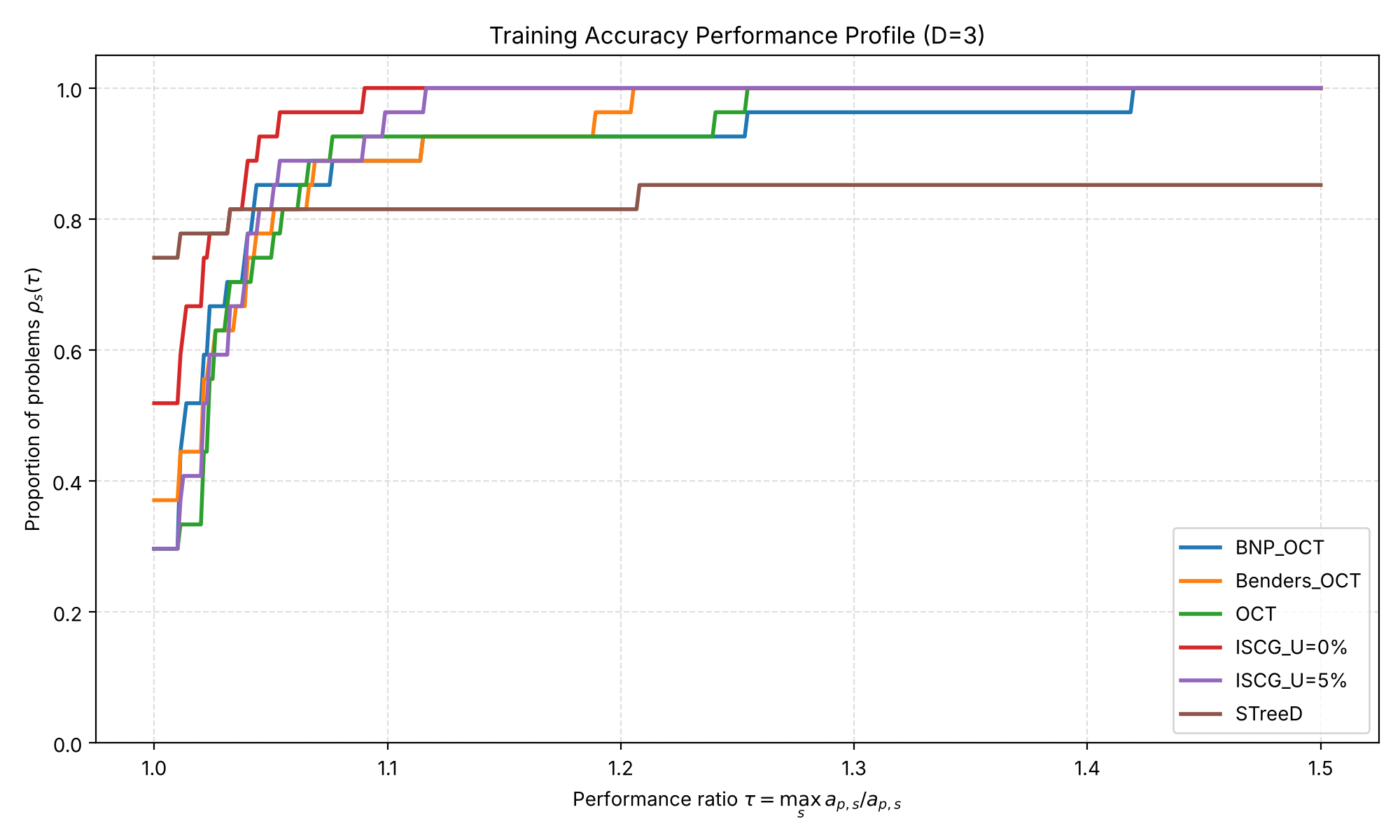}
        \caption{Training Accuracy Profile for Tree Depth 3}
        \label{fig:train_depth3}
    \end{subfigure}

    \vspace{0.3cm} 

    \begin{subfigure}[t]{0.48\textwidth}
        \centering
        \includegraphics[width=\linewidth]{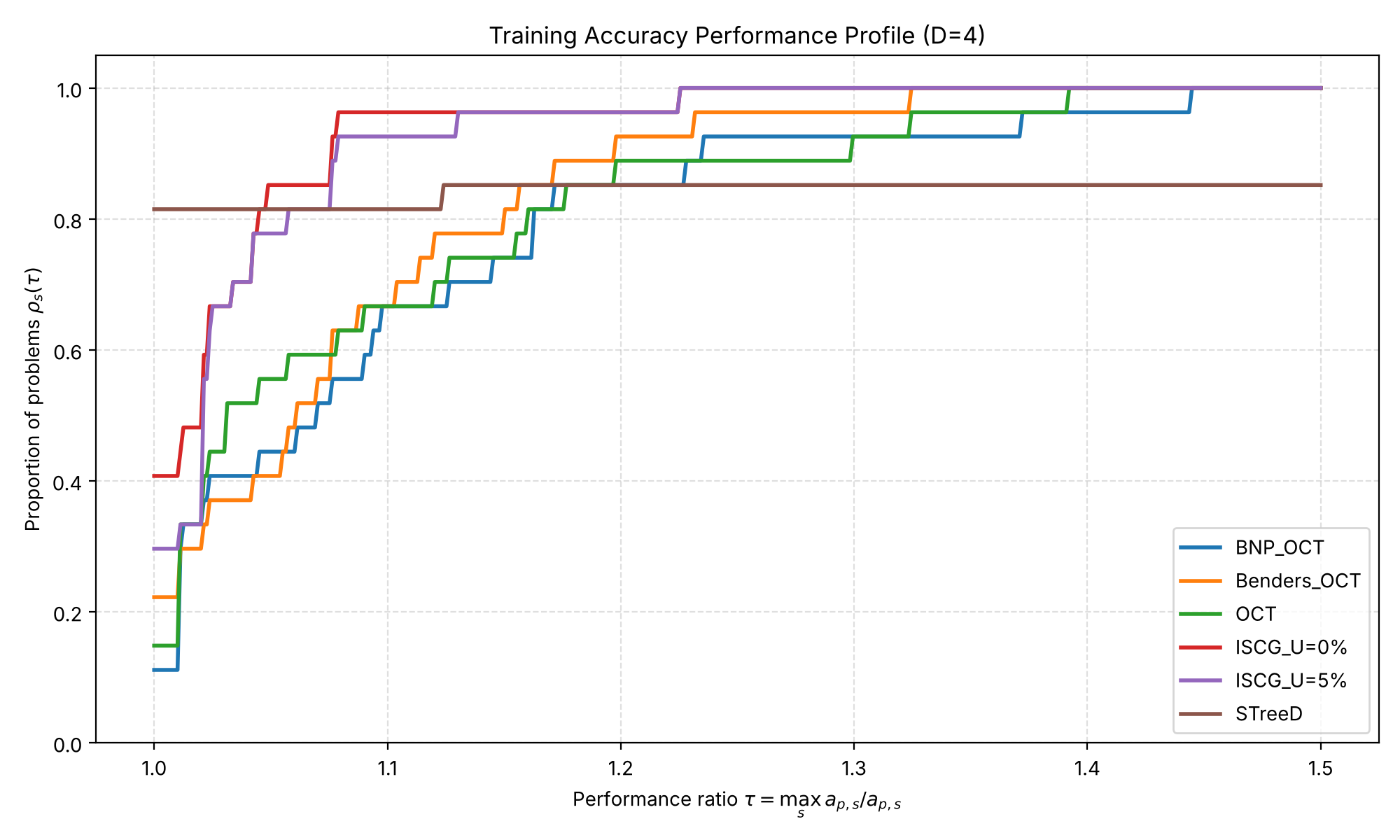}
        \caption{Training Accuracy Profile for Tree Depth 4}
        \label{fig:train_depth4}
    \end{subfigure}
    \hfill
    \begin{subfigure}[t]{0.48\textwidth}
        \centering
        \includegraphics[width=\linewidth]{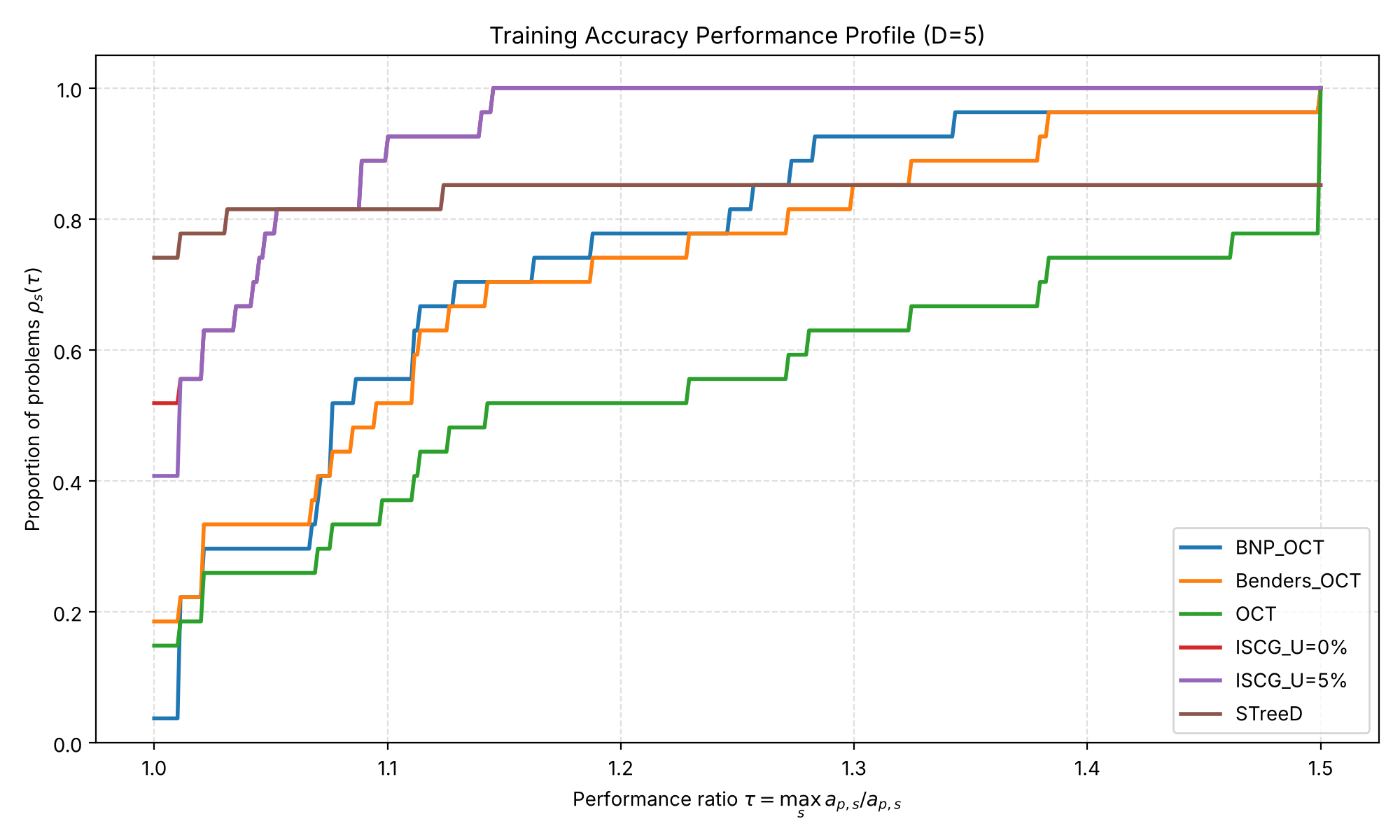}
        \caption{Training Accuracy Profile for Tree Depth 5}
        \label{fig:train_depth5}
    \end{subfigure}

    \caption{Training Accuracy Performance Profiles for Tree Depths $D=2,3,4,5$ (Time Limit of 3,600 Seconds)}
    \label{fig:train_accuracy_profiles}
\end{figure}

For tree depths $D=4$ and $D=5$, many benchmark instances cannot be solved to proven optimality within the one-hour time limit. STreeD attains the best reported training accuracy on approximately $81\%$ and $74\%$ of the datasets at $\tau=1$, compared with $41\%$ and $52\%$ for exact ISCG-OCT and at most $22\%$ for the MIP-based formulations, and its trees attain a higher training accuracy than those of exact ISCG-OCT on 16 and 13 of the 23 datasets on which it returns a tree, respectively. Both ISCG-OCT profiles, however, overtake that of STreeD before $\tau=1.1$, where exact ISCG-OCT covers approximately $96\%$ and $93\%$ of the datasets, respectively, compared with $81\%$ for STreeD and at most $67\%$ for the MIP-based formulations. The results are shown in Figure~\ref{fig:train_accuracy_profiles}.

\begin{figure}[h]
    \centering
    \begin{subfigure}[t]{0.48\textwidth}
        \centering
        \includegraphics[width=\linewidth]{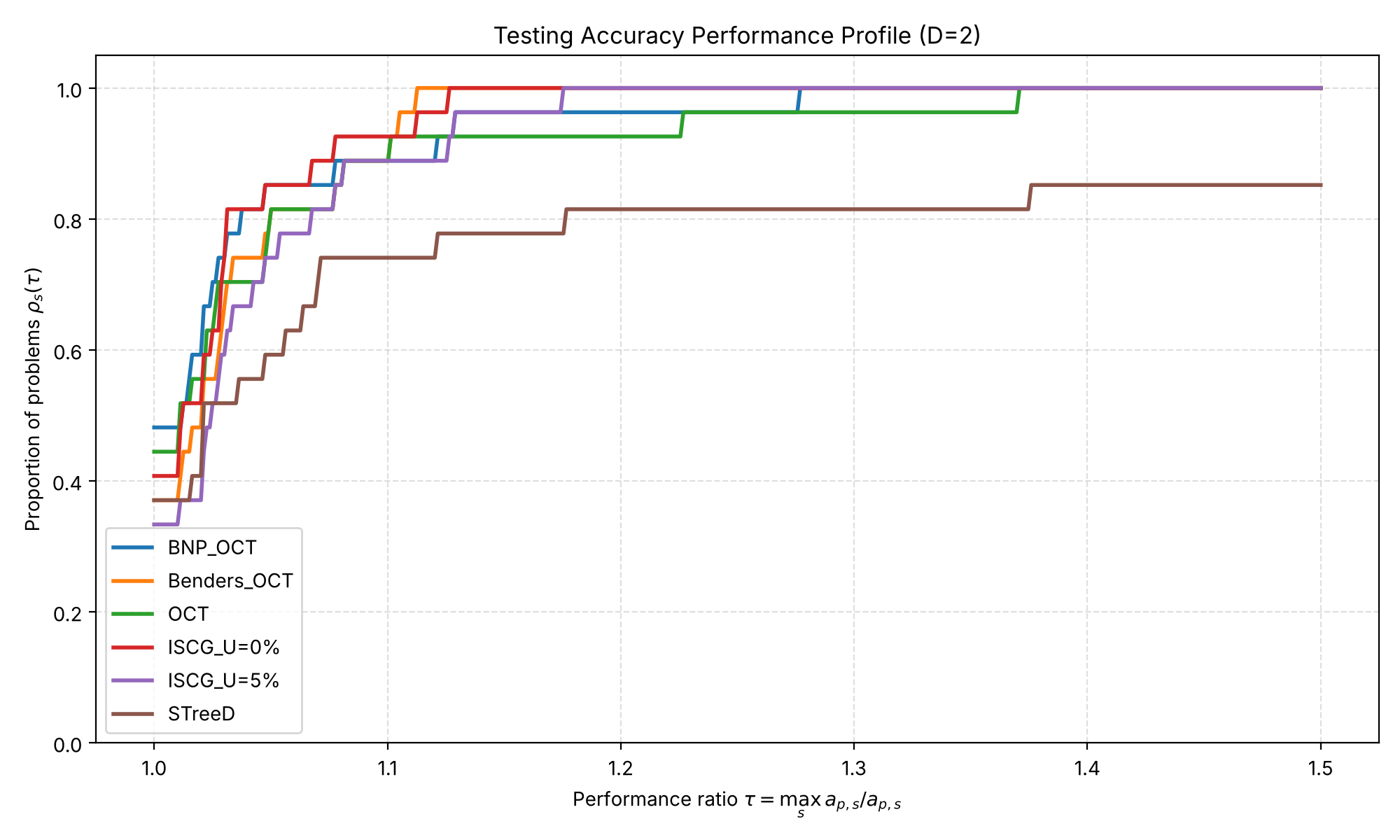}
        \caption{Testing Accuracy Profile for Tree Depth 2}
        \label{fig:test_depth2}
    \end{subfigure}
    \hfill
    \begin{subfigure}[t]{0.48\textwidth}
        \centering
        \includegraphics[width=\linewidth]{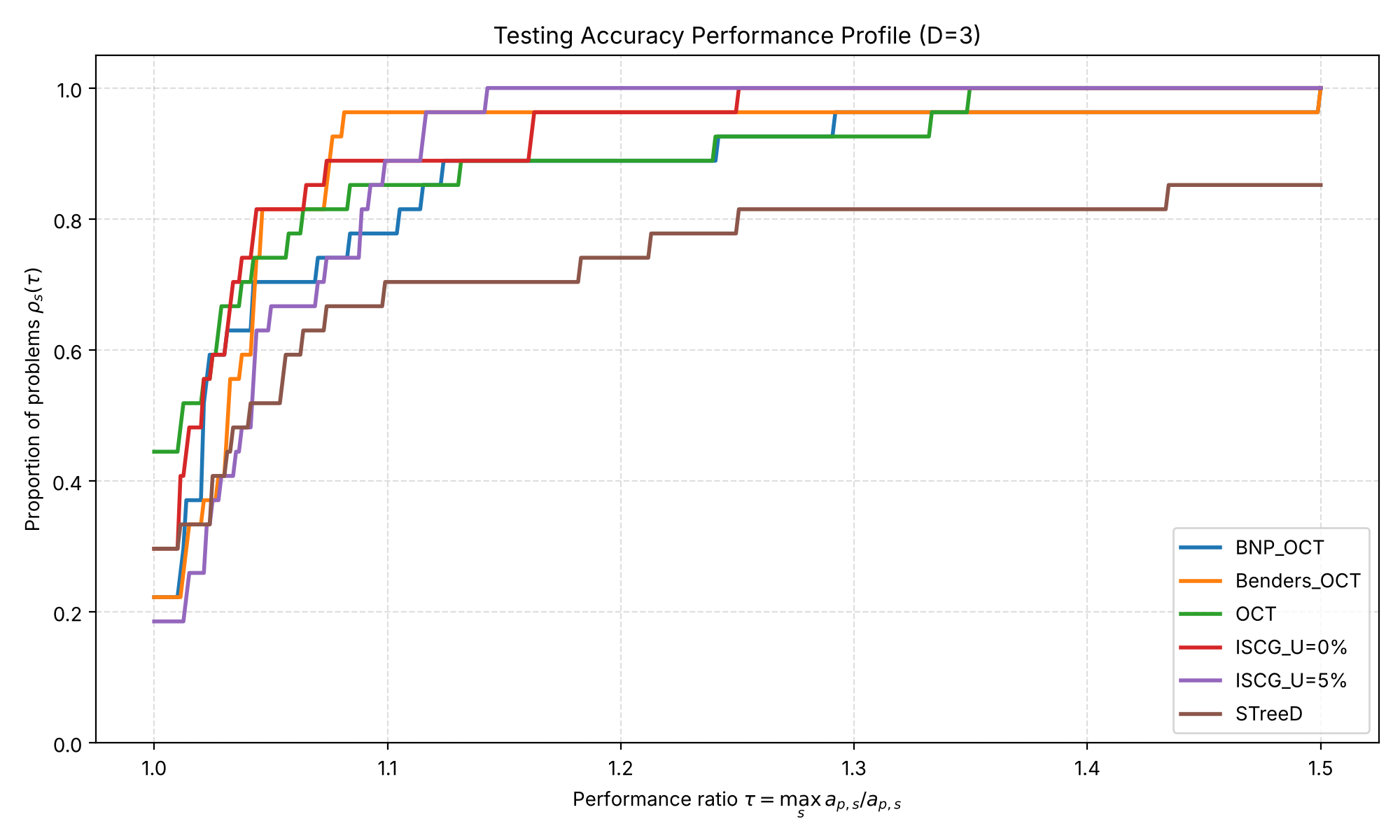}
        \caption{Testing Accuracy Profile for Tree Depth 3}
        \label{fig:test_depth3}
    \end{subfigure}

    \vspace{0.3cm} 

    \begin{subfigure}[t]{0.48\textwidth}
        \centering
        \includegraphics[width=\linewidth]{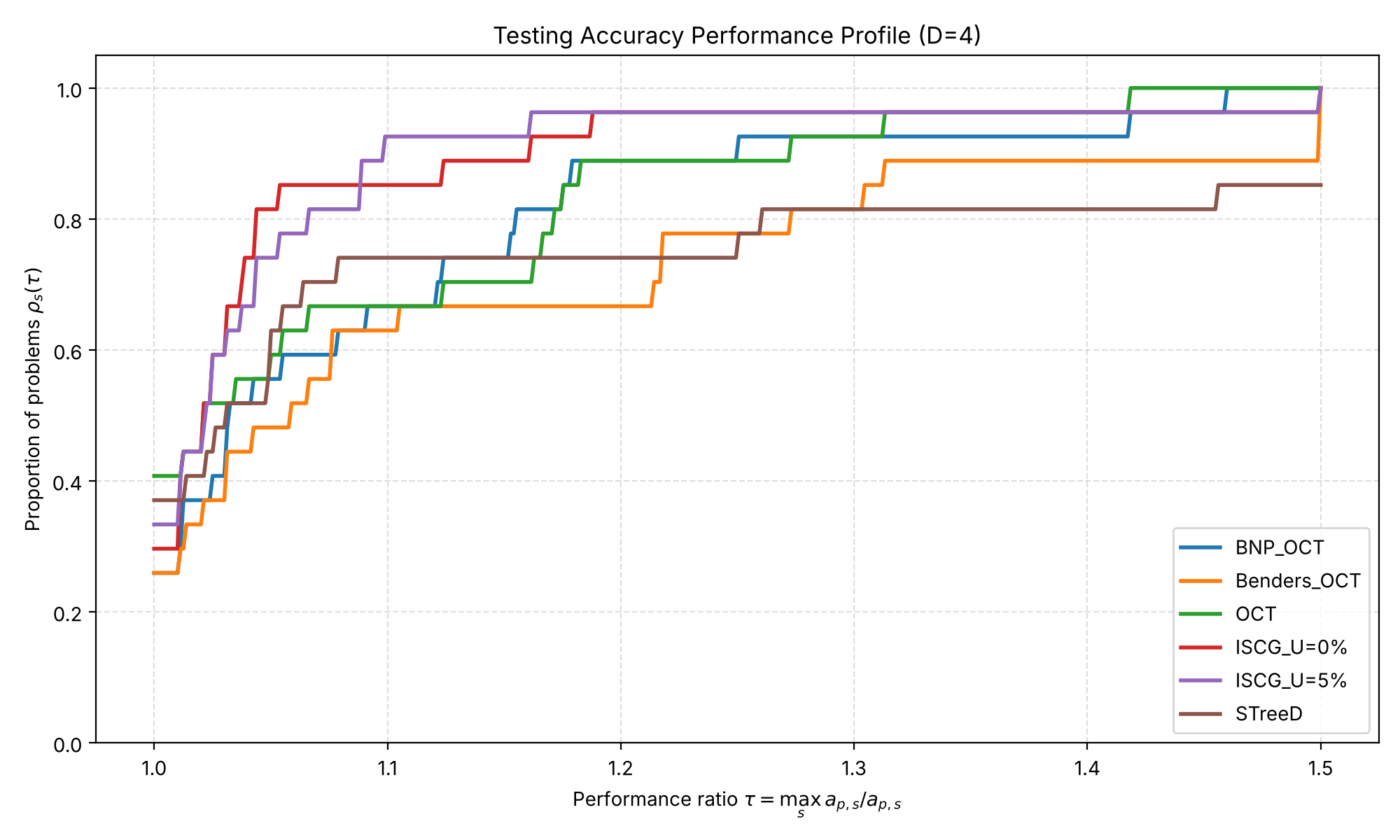}
        \caption{Testing Accuracy Profile for Tree Depth 4}
        \label{fig:test_depth4}
    \end{subfigure}
    \hfill
    \begin{subfigure}[t]{0.48\textwidth}
        \centering
        \includegraphics[width=\linewidth]{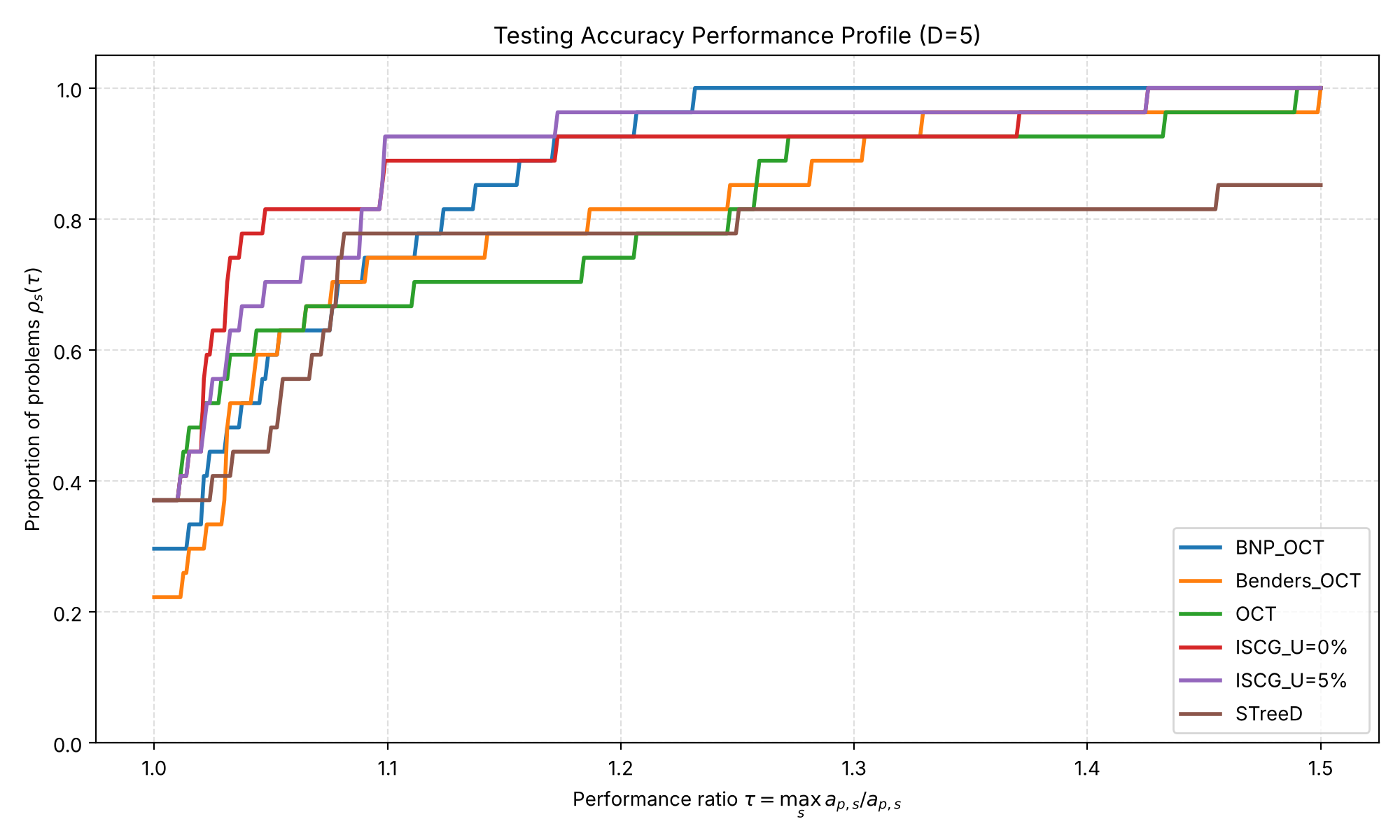}
        \caption{Testing Accuracy Profile for Tree Depth 5}
        \label{fig:test_depth5}
    \end{subfigure}

    \caption{Testing Accuracy Performance Profiles for Tree Depths $D=2,3,4,5$ (Time Limit of 3,600 Seconds)}
    \label{fig:profile_test}
\end{figure}

Figure~\ref{fig:profile_test} presents the testing accuracy performance profiles. Unlike the training accuracy results, no single algorithm consistently dominates across all benchmark instances and tree depths, reflecting the inherent trade-off between fitting the training data and generalization performance. Consequently, the testing accuracy results are more evenly distributed among the competing methods. At $\tau=1$, BNP-OCT attains the highest testing accuracy most often at $D=2$ (approximately $48\%$ of the datasets), OCT at $D=3$ and $D=4$ (approximately $44\%$ and $41\%$), and OCT, both ISCG-OCT variants, and STreeD at $D=5$ (approximately $37\%$ each). STreeD attains the highest testing accuracy on approximately $30\%$ to $37\%$ of the datasets at each depth, although its testing accuracy is measured on a different partition. At $\tau=1.05$, exact ISCG-OCT has the highest or tied-highest profile at every depth (approximately $81\%$ to $85\%$), compared with approximately $48\%$ to $63\%$ for STreeD. Overall, the results indicate that the proposed ISCG-OCT framework provides testing performance that is competitive with existing state-of-the-art exact OCT algorithms while offering the computational advantages over the MIP-based formulations reported in the next subsection.

\subsection{Comparison Results on Speed}

For tree depth $D=2$, the methods exhibit differences in their ability to certify optimality within the 3,600-second time limit. The proposed ISCG-OCT algorithm with a $0\%$ pruning threshold proves optimality for every benchmark instance except Spambase, for which none of the competing methods establishes optimality within the allotted time. In comparison, BNP-OCT and OCT terminate with nonzero optimality gaps on eight of the 27 benchmark instances, while Benders OCT also has instances with unresolved optimality gaps. STreeD certifies optimality on all 23 instances on which it stays within its memory limit, in at most 10.4 seconds, but exceeds the memory limit on Spambase, Breast Cancer Wisconsin (Diagnostic), Biodeg, and Ozone Level Detection, for which the exhaustive binarization produces between 6,386 and 19,835 binary features. Exact ISCG-OCT certifies optimality on three of these four datasets, because its recursive pricing considers only the candidate split values. These results show that ISCG-OCT certifies optimality on the largest number of the depth-two benchmark instances within the imposed limits.

The thresholded ISCG-OCT variant $U=5\%$ further reduces computational effort by pruning recursive states that cannot improve the incumbent objective by at least the prescribed threshold. Consequently, it typically terminates earlier than the exact variant while maintaining high-quality feasible solutions. For $D=2$, the reduction in computational times is achieved with only a negligible loss in solution quality because the optimal solutions are generally close to the CART warm-start on most benchmark instances.

\begin{figure}[h]
    \centering
    \begin{subfigure}[t]{0.48\textwidth}
        \centering
        \includegraphics[width=\linewidth]{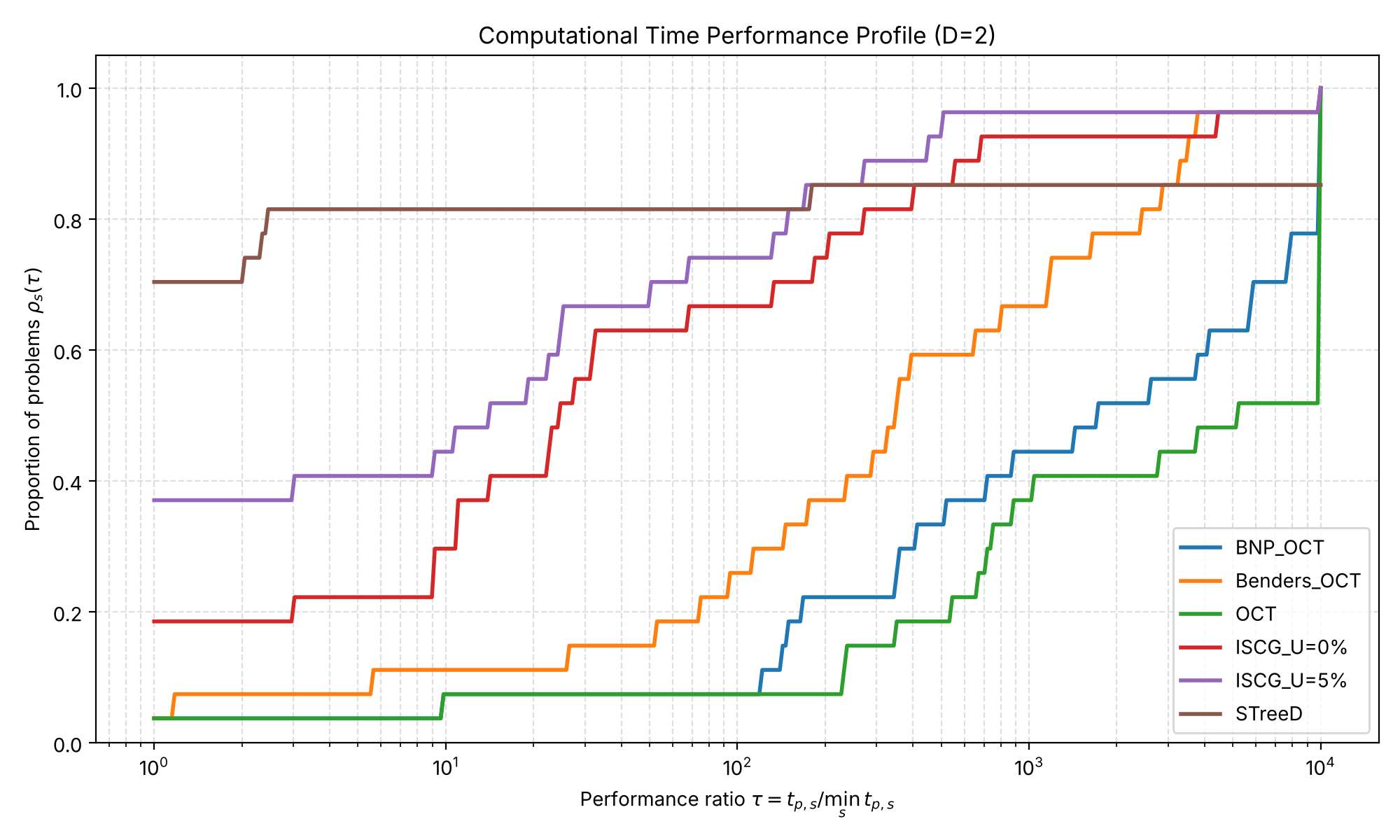}
        \caption{computational times Profile for Tree Depth 2}
        \label{fig:time_depth2}
    \end{subfigure}
    \hfill
    \begin{subfigure}[t]{0.48\textwidth}
        \centering
        \includegraphics[width=\linewidth]{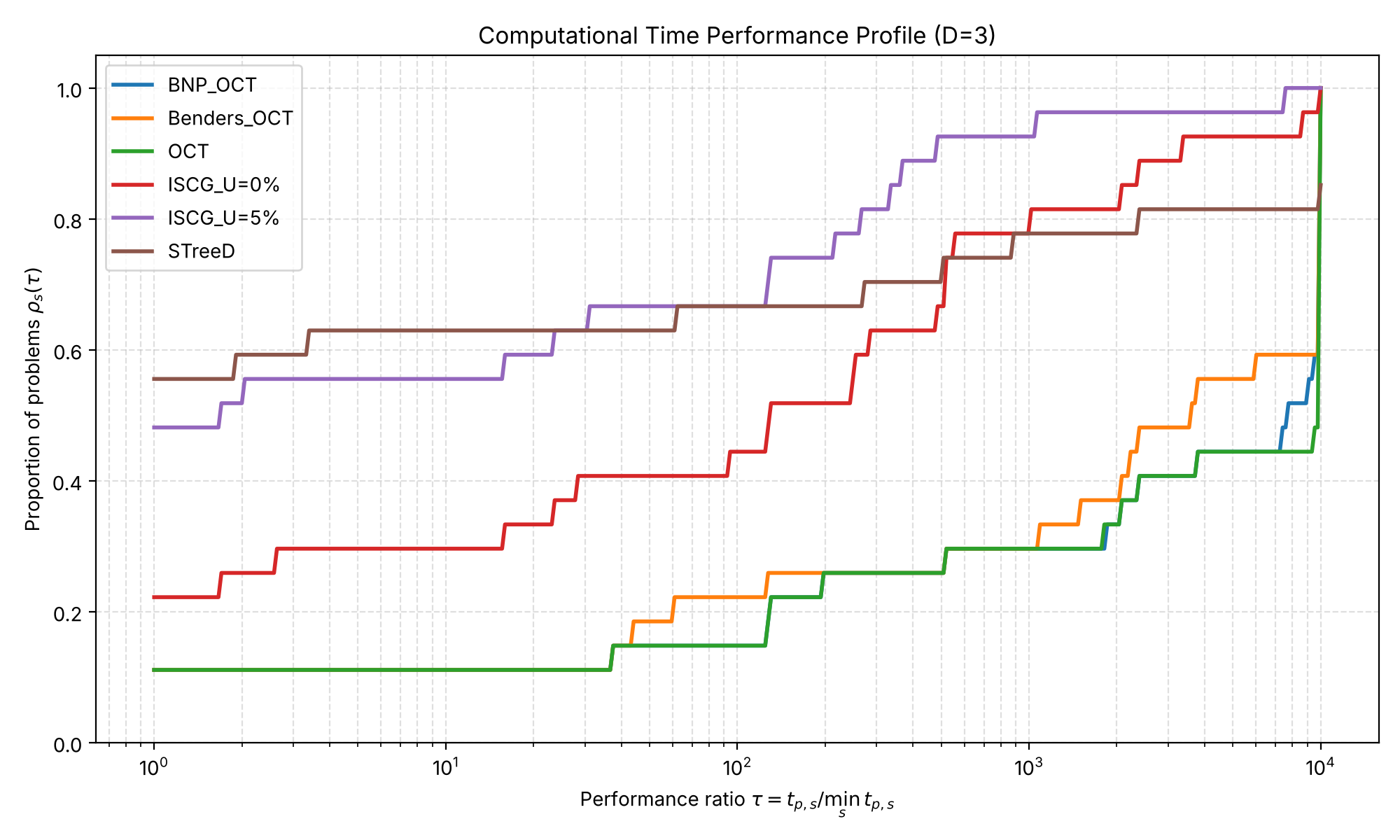}
        \caption{computational times Profile for Tree Depth 3}
        \label{fig:time_depth3}
    \end{subfigure}

    \vspace{0.3cm} 

    \begin{subfigure}[t]{0.48\textwidth}
        \centering
        \includegraphics[width=\linewidth]{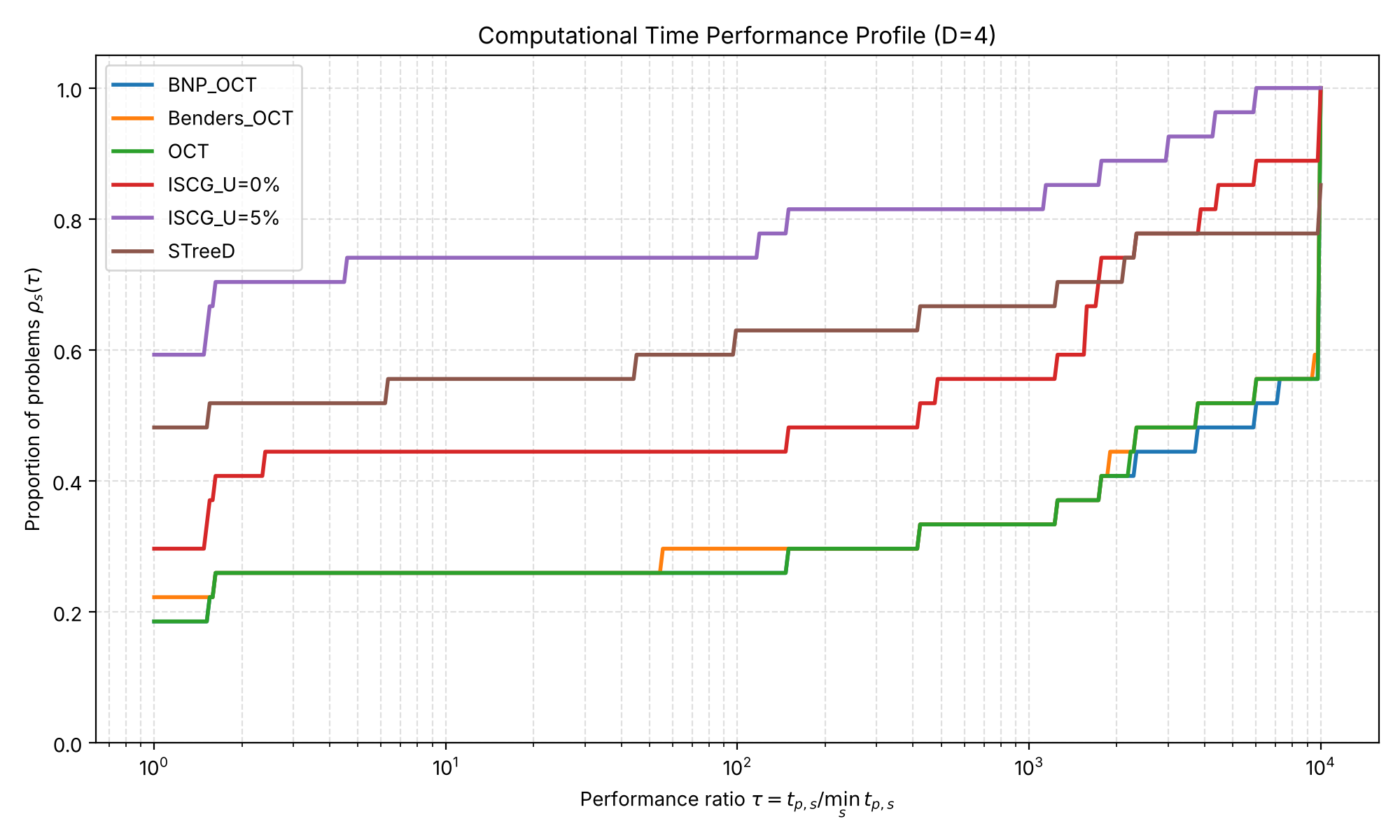}
        \caption{computational times Profile for Tree Depth 4}
        \label{fig:time_depth4}
    \end{subfigure}
    \hfill
    \begin{subfigure}[t]{0.48\textwidth}
        \centering
        \includegraphics[width=\linewidth]{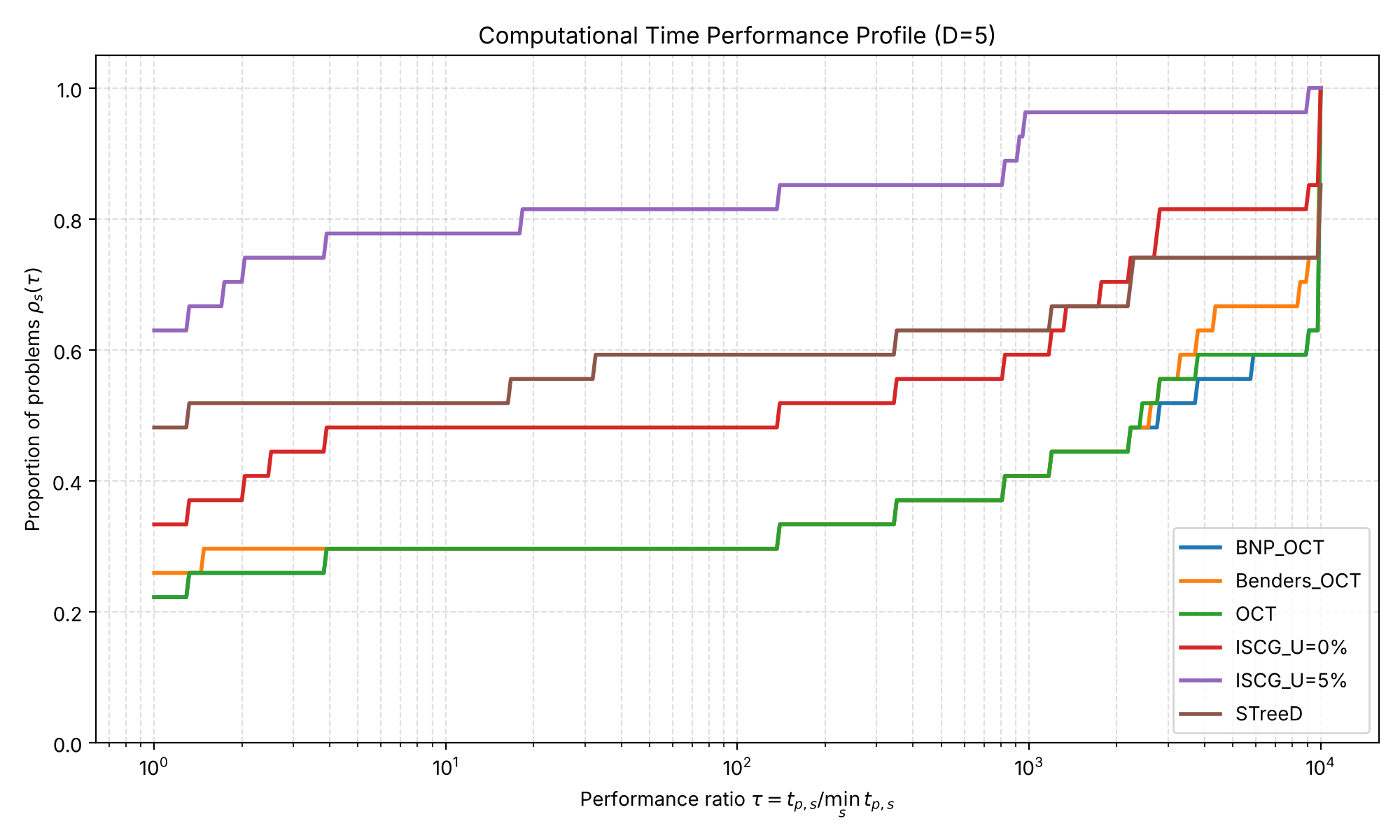}
        \caption{computational times Profile for Tree Depth 5}
        \label{fig:time_depth5}
    \end{subfigure}

    \caption{computational times Performance Profiles for Tree Depths $D=2,3,4,5$ (Time Limit of 3,600 Seconds)}
    \label{profile_computationaltime}
\end{figure}

Figure~\ref{profile_computationaltime} shows that both ISCG-OCT variants generally achieve shorter computational times than the MIP-based formulations under their respective termination criteria, with the thresholded variant $U=5\%$ providing an additional speedup over the exact variant $U=0\%$. STreeD has the highest profile at $D=2$ and $D=3$ up to ratios of about 170 and 30, respectively (approximately $70\%$ and $56\%$ at $\tau=1$), whereas ISCG $U=5\%$ has the highest profile beyond these ratios and at every ratio for $D=4$ and $D=5$ (approximately $59\%$ and $63\%$ at $\tau=1$, compared with $48\%$ for STreeD). The STreeD profile lies above that of exact ISCG-OCT up to ratios of several hundred at $D=2$ and $D=3$ and above 1,000 at $D=4$ and $D=5$. Among the instances on which both methods certify optimality, STreeD terminates sooner on 20 of 23, 13 of 16, 6 of 10, and 4 of 7 at $D=2,3,4,5$. The MIP-based profiles are the lowest at almost all ratios.

ISCG-OCT demonstrates substantial computational advantages over the MIP-based formulations on several benchmark datasets. For example, on the Thoracic Surgery dataset with $D=2$, ISCG-OCT solves the instance in less than one second. Benders OCT and BNP-OCT also certify optimality, requiring 72.32 and 328.28 seconds, respectively, while OCT reaches the 3,600-second limit with a reported $0\%$ optimality gap. This corresponds to observed speedups of at least 4,390, 87, and 390 relative to OCT, Benders OCT, and BNP-OCT, respectively; STreeD certifies the same instance in 0.03 seconds. Similarly, for the Tic-tac-toe dataset, ISCG-OCT obtains the optimal solution in approximately 10 seconds for $D=3$, compared with 0.02 seconds for STreeD. On the Mushroom dataset with $D=4$, which contains relatively few candidate split values, ISCG-OCT solves the OCT instance in approximately 15 seconds. Even at tree depth $D=5$, the solution time increases to only about 190 seconds, whereas none of the MIP-based methods proves optimality within the one-hour time limit. STreeD certifies both Mushroom instances in under two seconds, with trees that classify all training observations correctly.

\subsection{Comparison Results on Optimality Gap}

Unlike conventional branch-and-bound algorithms, ISCG-OCT obtains its optimality bound through the recursive pricing procedure rather than the primal--dual bounds maintained by a MIP solver. To provide a numerical comparison, we calculate the relative optimality gap for each algorithm using its best feasible objective value and a valid global lower bound for its respective formulation. For the MIP-based benchmark methods, these bounds are obtained from the solver. For ISCG-OCT, the upper bound is given by the current incumbent tree, while the global lower bound is constructed by combining the objective-based lower bound implied by MSSF with the corrected lower bounds of the recursive-pricing states. The resulting bound provides a global optimality certificate only when the maintained states collectively cover the complete candidate-restricted MSSF feasible region and the combined bound is valid for all unexplored completions. STreeD reports no lower bound when it stops at the time limit. Its gap is zero when it certifies optimality and is otherwise computed from a valid but weak lower bound described in Appendix~\ref{gapcomputation}.

For ISCG $U=0\%$, the reported gap is zero when exhaustive recursive pricing with valid bounds certifies global optimality within the 3,600-second limit. Otherwise, it is calculated from the best feasible objective and the valid global lower bound at termination, so a positive gap indicates that optimality has not been certified. For ISCG $U=5\%$, the positive improvement threshold controls the recursive search rather than imposing a fixed gap: the reported gap can be zero when global optimality is certified and is otherwise calculated in the same way, rather than being set to $5\%$.

\begin{figure}[h]
    \centering
    \begin{subfigure}[t]{0.48\textwidth}
        \centering
        \includegraphics[width=\linewidth]{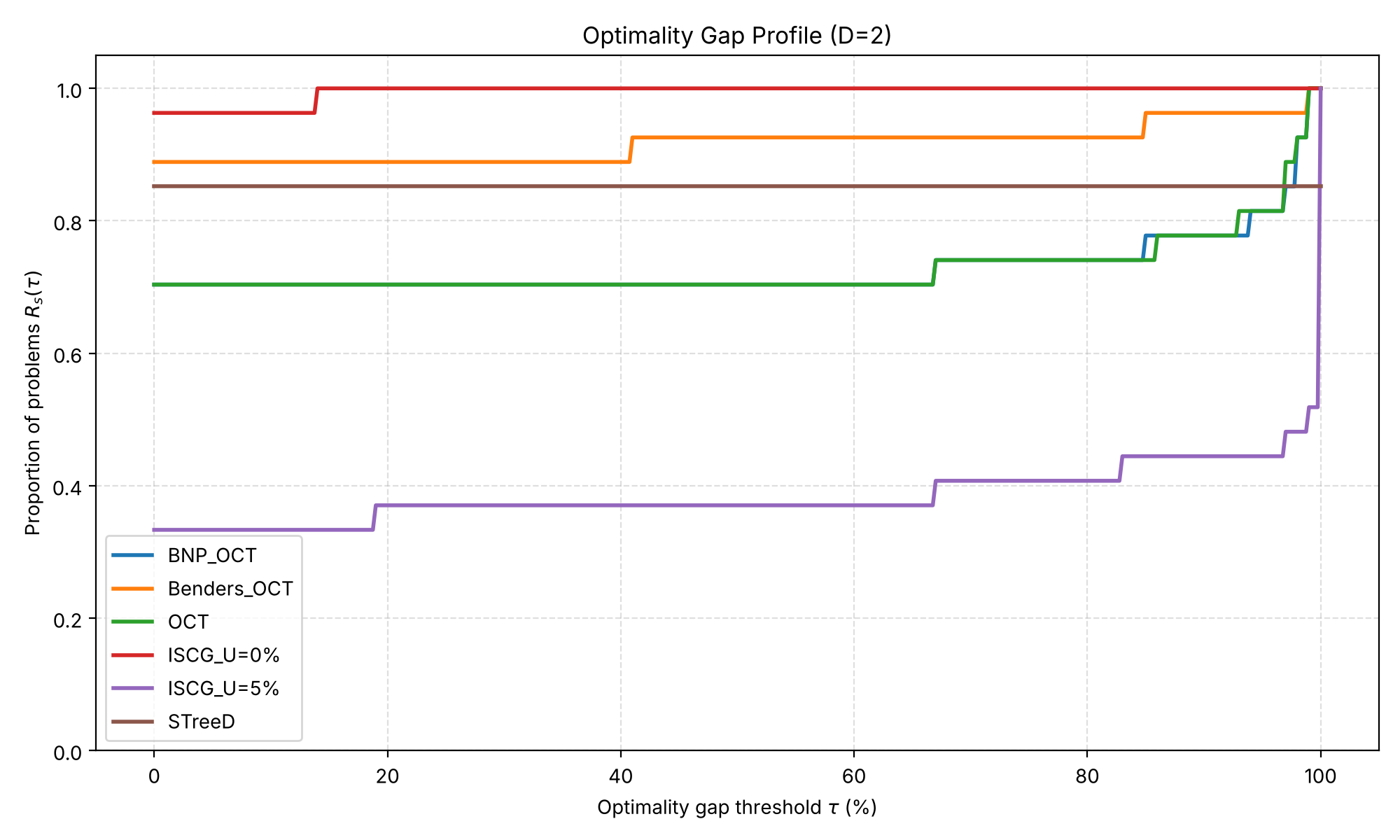}
        \caption{Optimality Gap Profile for Tree Depth 2}
        \label{fig:gap_depth2}
    \end{subfigure}
    \hfill
    \begin{subfigure}[t]{0.48\textwidth}
        \centering
        \includegraphics[width=\linewidth]{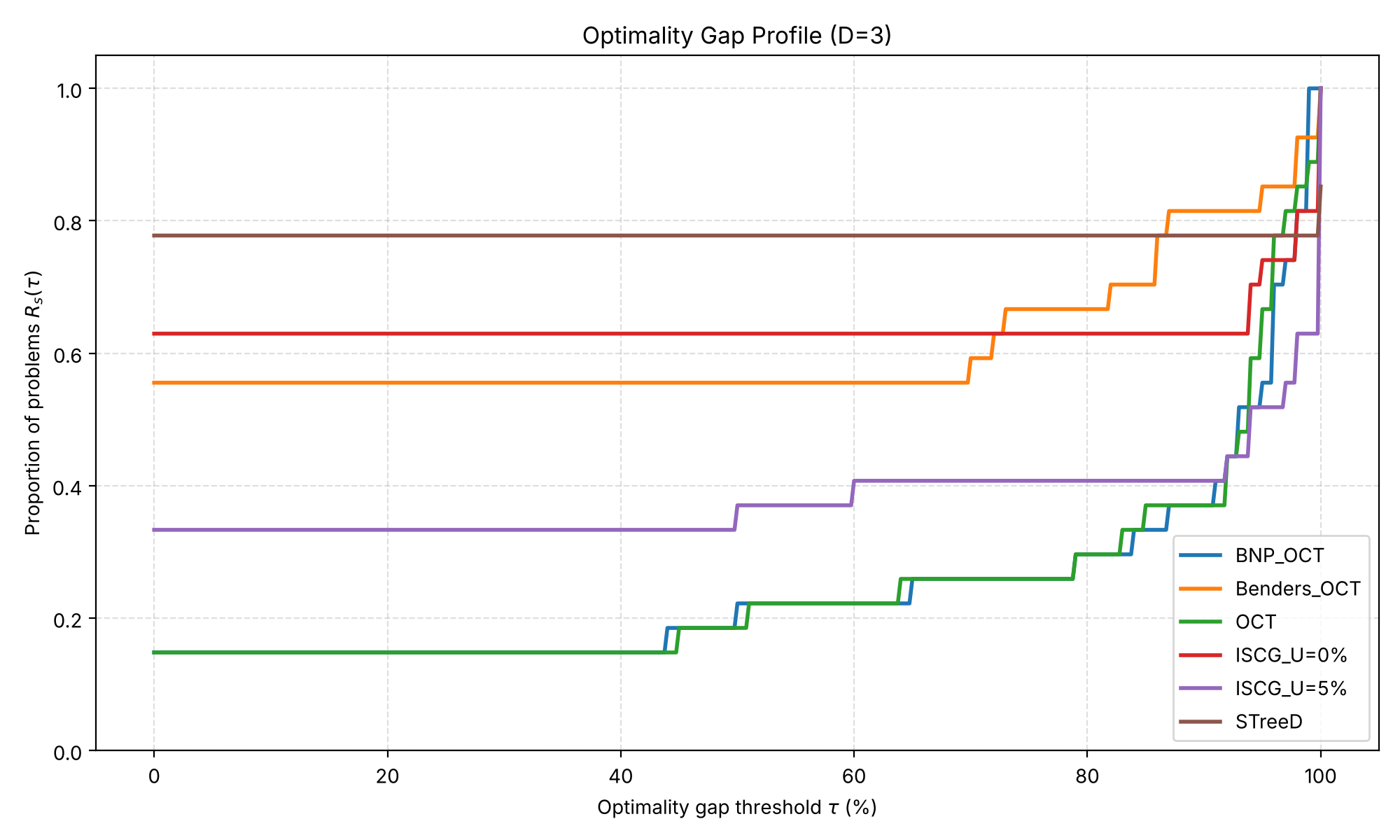}
        \caption{Optimality Gap Profile for Tree Depth 3}
        \label{fig:gap_depth3}
    \end{subfigure}

    \vspace{0.3cm} 

    \begin{subfigure}[t]{0.48\textwidth}
        \centering
        \includegraphics[width=\linewidth]{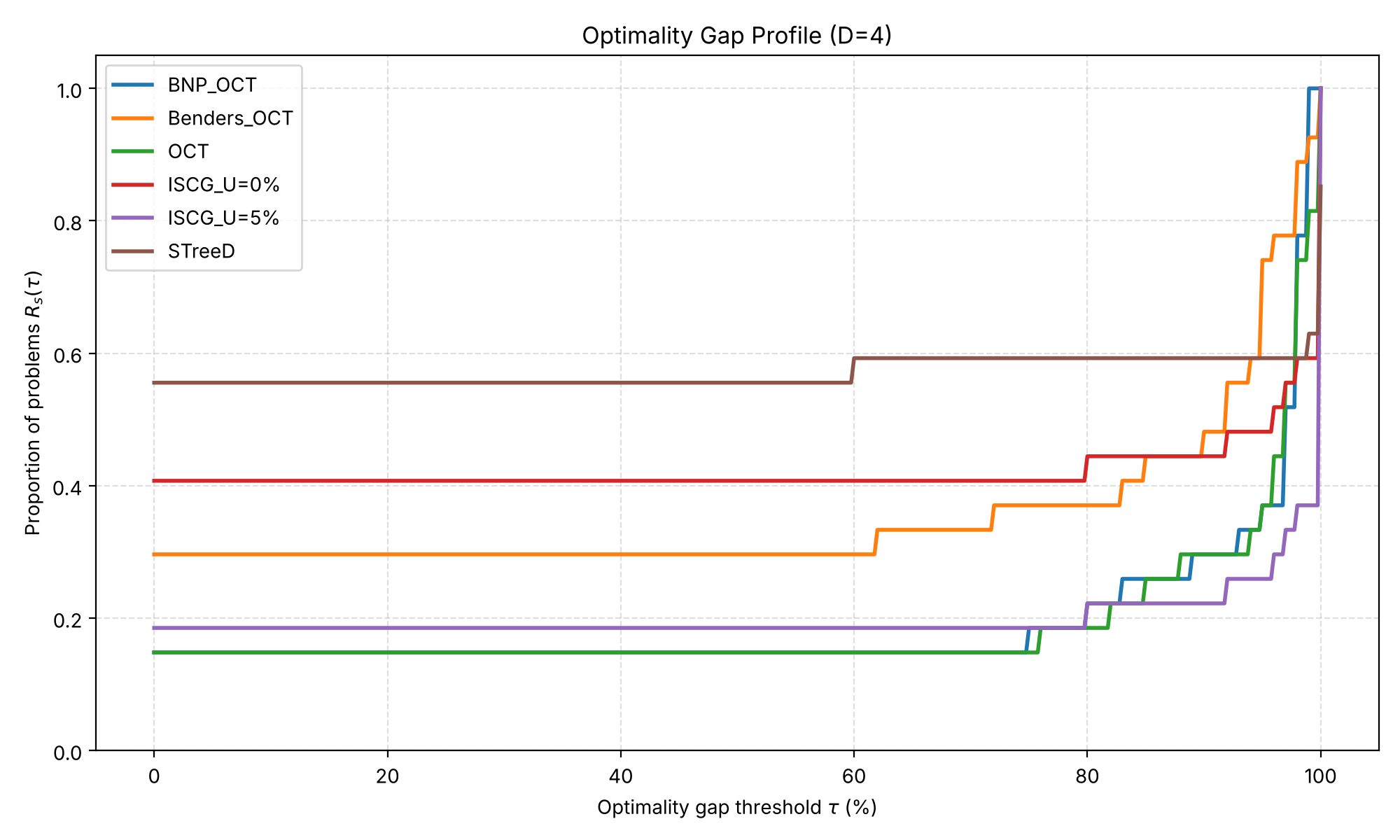}
        \caption{Optimality Gap Profile for Tree Depth 4}
        \label{fig:gap_depth4}
    \end{subfigure}
    \hfill
    \begin{subfigure}[t]{0.48\textwidth}
        \centering
        \includegraphics[width=\linewidth]{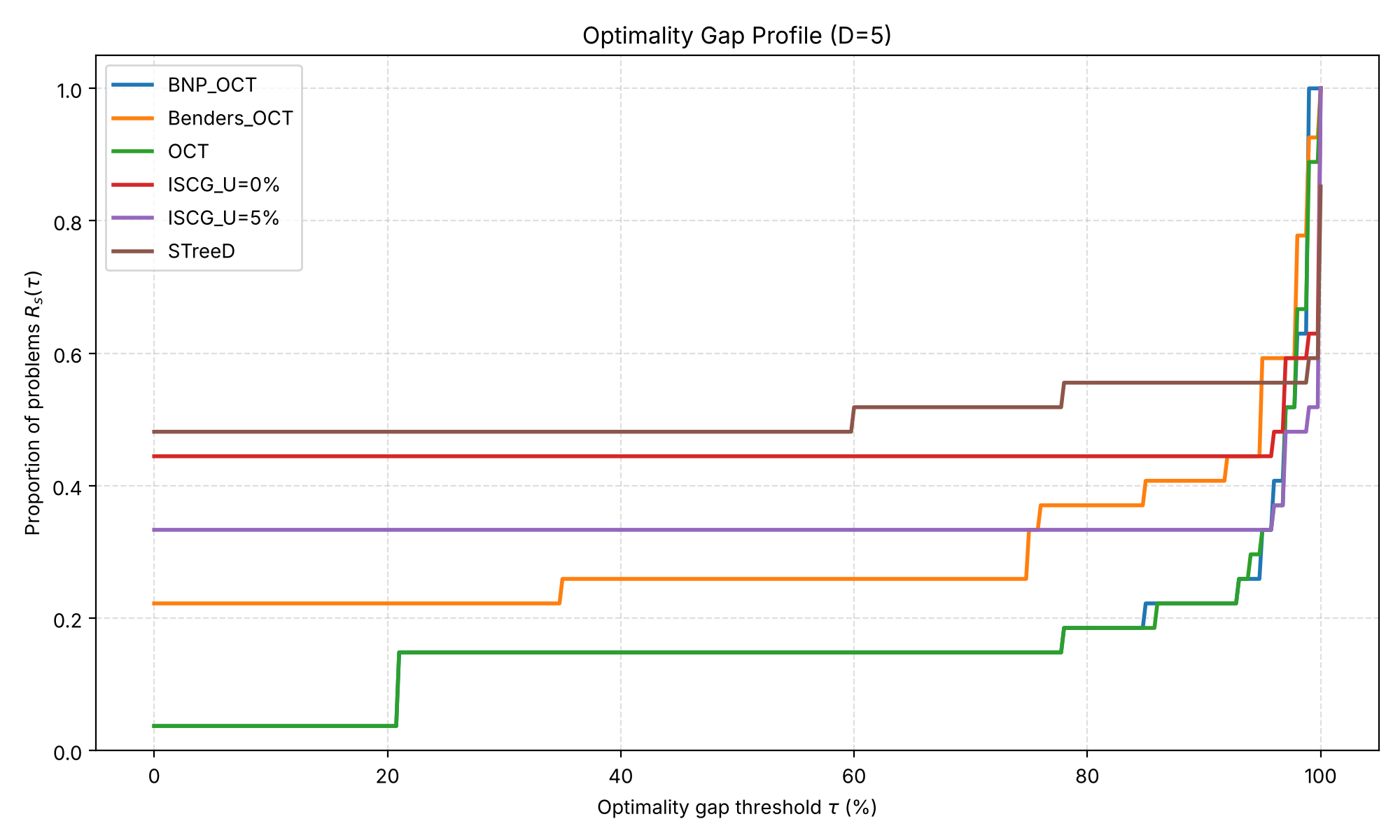}
        \caption{Optimality Gap Profile for Tree Depth 5}
        \label{fig:gap_depth5}
    \end{subfigure}

    \caption{Optimality Gap Profiles for Tree Depths $D=2,3,4,5$ (Gaps Relative to Each Method's Formulation)}
    \label{fig:gap_profiles}
\end{figure}

Figure~\ref{fig:gap_profiles} presents the optimality gap profiles for tree depths $D=2,3,4,5$. For each method, the profile $R_s(\tau)$ represents the proportion of datasets for which the reported optimality gap is no greater than the threshold $\tau$. Thus, a higher curve indicates better performance in terms of the reported gaps. At $D=2$, ISCG $U=0\%$ achieves a zero gap on approximately $96\%$ of the datasets and a gap of $14\%$ on the remaining dataset, while the other methods require substantially larger thresholds to attain comparable coverage; STreeD reports a zero gap on approximately $85\%$ of the datasets and no gap on the four datasets on which it exceeds the memory limit. As tree depth increases, the profiles generally shift downward, reflecting the greater computational difficulty of solving deeper trees. At $D=3,4,5$, STreeD reports zero gaps on approximately $78\%$, $56\%$, and $48\%$ of the datasets, compared with $63\%$, $41\%$, and $44\%$ for ISCG $U=0\%$ and $56\%$, $30\%$, and $22\%$ for Benders OCT. At these depths, the profiles of STreeD and both ISCG-OCT variants change little up to a threshold of $90\%$, because almost all of their unresolved runs report gaps above $90\%$, most of them above $99\%$ (for STreeD, based on a weak lower bound), whereas most unresolved gaps of the MIP-based formulations lie between $90\%$ and $99\%$, so that their profiles rise steeply and overtake those of ISCG-OCT and STreeD at large thresholds. ISCG $U=5\%$ also performs competitively at deeper tree depths, although its performance is generally weaker than ISCG $U=0\%$. Benders OCT exhibits intermediate performance across the four depths. Overall, these results suggest that ISCG $U=0\%$ provides strong optimality-gap performance for its formulation, with higher profiles than the MIP-based formulations for thresholds up to approximately $70\%$ at every depth, whereas STreeD certifies optimality on more instances at depths three to five despite exceeding its memory limit on four datasets. Differences in objective definitions, feasible regions, and split candidate sets should be considered when interpreting the comparisons across methods.

\section{Conclusions}\label{sec:conclusions}

This paper develops ISCG-OCT, an incumbent subtree column-generation framework for optimal classification trees. MSSF separates branching and assignment decisions while preserving the routing and leaf-label structure of the tree. Building on MSSF, we construct an RMP whose columns are closed feasible subtree replacement patterns and derive the corresponding subtree-pricing problems, in which objective changes relative to the integer-feasible incumbent replace dual-based reduced costs. Because fixing a branching or deactivation decision separates the remaining pricing problem into independent left- and right-subtree subproblems, these problems can be solved exactly by DP. ISCG-OCT combines the RMP with this recursive pricing to generate improving integer-feasible subtree patterns while maintaining global feasibility. When the complete root-pricing problem is solved exactly, a nonnegative pricing value certifies that the incumbent is globally optimal over the candidate-restricted MSSF feasible region.

To improve computational efficiency, a lower-bound-guided search prunes partial subtree configurations that cannot improve the incumbent, and a thresholded variant expands a partial state only if it may improve the incumbent by at least a prescribed amount and terminates when no such replacement is identified. Unlike the exact procedure, the thresholded variant trades optimality certification for reduced computational effort. On 27 benchmark datasets, ISCG-OCT substantially reduces solution times on many instances while achieving training and testing accuracy competitive with existing exact OCT methods. Future research may develop stronger bounds for deeper trees and extend recursive pricing to other interpretable tree models, including regression, prescriptive, and fairness-constrained trees.

\bibliographystyle{plainnat}
\bibliography{reference}

\newpage

\appendix
\phantomsection
\label{app:proofs}

\small
\setstretch{1.05}


\setcounter{figure}{0}
\renewcommand{\thefigure}{A.\arabic{figure}}

\setcounter{table}{0}
\renewcommand{\thetable}{A.\arabic{table}}

\setcounter{lemma}{0}
\setcounter{proposition}{3}
\section{Proofs of Section~\ref{sec:mssf} Results}
\label{app:proofs_mssf}

This section proves the routing and leaf-label consistency of MSSF, gives the class-boundary candidate-split algorithm, and shows that the linear relaxation of a depth-one, single-feature MSSF describes the convex hull of its integer-feasible solutions.

\subsection{Routing and Leaf-Label Consistency}

\begin{lemma}[Routing Consistency of MSSF]
\label{lem:mssf_routing_consistency}
Let $(b,z)$ be an integer-feasible solution of MSSF. If a branching alternative $q\in\mathcal Q$ selected at node $v\in T_B$ makes observation $i\in I$ infeasible for reaching leaf $t\in T_L$, then $z_{ikt}=0$ for all $k\in K$. This includes deactivation ($q=0$), which makes every observation infeasible for reaching a leaf in the left subtree of the deactivated node.
\end{lemma}
\proof{Proof.}
Let $b_{qv}=1$. If $v\in P_L(t)$ and $i\in S_q^l$, then, for every $k\in K$, some element of $\mathcal C^l(t;q)$ contains $(i,k,t)$, and the corresponding inequality~\eqref{eq:mssf_left}, together with nonnegativity, forces $z_{ikt}=0$. If $v\in P_R(t)$ and $i\in S_q^r$, inequality~\eqref{eq:mssf_right} gives $z_{ikt}=0$ in the same way. For $q=0$, $S_0^l=I$ and $S_0^r=\emptyset$, so deactivation excludes every observation from the left subtree and none from the right subtree. Hence every assignment is consistent with the selected branching and deactivation decisions. \Halmos
\endproof

\begin{lemma}[Leaf-Label Consistency of MSSF]
\label{lem:mssf_leaf_label_consistency}
Let $t^{\mathrm{right}}$ be the rightmost leaf, which corresponds to the compatibility-only artificial upper threshold with right-routing exclusion set $I$, and suppose that MSSF includes the class-consistency inequalities
\begin{equation}
\sum_{k\in K}z_{i_kkt^{\mathrm{right}}}\le1,\qquad\forall(i_k)_{k\in K}\in I^{|K|}.
\label{eq:leaf_label_rightmost}
\end{equation}
If $(b,z)$ is an integer-feasible solution of MSSF, then every nonempty leaf $t\in T_L$ has exactly one class $k\in K$ such that $z_{ikt}=1$ for some $i\in I$; equivalently, $z_{ikt}=z_{i'k't}=1$ implies $k=k'$.
\end{lemma}
\proof{Proof.}
Fix $t\in T_L$. We first show that
\begin{equation}
\sum_{k\in K}z_{i_kkt}\le1,\qquad\forall(i_k)_{k\in K}\in I^{|K|}.
\label{eq:leaf_label_cp_implied}
\end{equation}
If $P_L(t)\neq\emptyset$, choose $v\in P_L(t)$. Since $S^l_0=I$, the collection $\mathcal C^l(t;0)$ contains $\{(i_k,k,t):k\in K\}$ for every tuple $(i_k)_{k\in K}\in I^{|K|}$, and the corresponding inequality~\eqref{eq:mssf_left}, $b_{0v}+\sum_{q'\in B_0^l}b_{q'v}+\sum_{k\in K}z_{i_kkt}\le1$, gives~\eqref{eq:leaf_label_cp_implied} because the branching variables are nonnegative. If $P_L(t)=\emptyset$, then $t=t^{\mathrm{right}}$ because the tree is binary; the artificial upper threshold is not a selectable alternative, and its exclusion set $I$ is represented by~\eqref{eq:leaf_label_rightmost}, which is~\eqref{eq:leaf_label_cp_implied}. Now suppose that $z_{ikt}=z_{i'k't}=1$ with $k\neq k'$. Choosing $i_k=i$, $i_{k'}=i'$, and arbitrary indices for the remaining classes in~\eqref{eq:leaf_label_cp_implied} gives a left-hand side of at least two, a contradiction. Hence every nonempty leaf, which contains at least one assigned observation, has exactly one assigned class. \Halmos
\endproof

\subsection{Class-Boundary Candidate Split Reduction Algorithm}
\label{class_boundary}

Algorithm~\ref{alg:recursive_class_boundary} constructs the candidate sets before MSSF is solved by recursively enumerating boundary-induced observation subsets until no previously unseen subset is generated.

\begin{breakablealgorithm}
\caption{Recursive Class-Boundary Split Values}
\label{alg:recursive_class_boundary}
\begin{algorithmic}[1]
\Require Observations $I$, features $J$, values $x_{ij}$, labels $y_i$, and nontrivial observed split sets $U_j^{\mathrm{split}}=\{u_j^{(1)},\ldots,u_j^{(M_j-1)}\}$, $j\in J$.
\Ensure Candidate split sets $\mathcal R_j$, $j\in J$.
\State Set $\OmegaSet_0\gets\{I\}$, $\OmegaSet^{\mathrm{seen}}\gets\{I\}$, $\mathcal R_j\gets\emptyset$ for all $j\in J$, and $d\gets0$.
\While{$\OmegaSet_d\neq\emptyset$}
  \State Set $\OmegaSet_{d+1}\gets\emptyset$.
  \For{each unique subset $\Omega\in\OmegaSet_d$ and each feature $j\in J$}
    \State Sort the distinct feature-$j$ value groups $G_{j,\Omega}^{(1)},\ldots,G_{j,\Omega}^{(M_j(\Omega))}$ of $\Omega$.
    \State $\mathcal B_j(\Omega)\gets\{r:\exists\,i\in G_{j,\Omega}^{(r)},\ i'\in G_{j,\Omega}^{(r+1)}\text{ with }y_i\neq y_{i'}\}$.
    \State $\mathcal C_j(\Omega)\gets\{q\in\{1,\ldots,M_j(\Omega)-1\}:q\in\mathcal B_j(\Omega)\text{ or }q-1\in\mathcal B_j(\Omega)\}$.
    \For{each $q\in\mathcal C_j(\Omega)$}
      \State Add $u_{j,\Omega}^{(q)}$ to $\mathcal R_j$; let $\Omega_{j,q}^{L}=\{i\in\Omega:x_{ij}\le u_{j,\Omega}^{(q)}\}$ and $\Omega_{j,q}^{R}=\Omega\setminus\Omega_{j,q}^{L}$.
      \State Add each nonempty subset in $\{\Omega_{j,q}^{L},\Omega_{j,q}^{R}\}$ not in $\OmegaSet^{\mathrm{seen}}$ to $\OmegaSet_{d+1}$ and $\OmegaSet^{\mathrm{seen}}$.
    \EndFor
  \EndFor
  \State Remove duplicate subsets from $\OmegaSet_{d+1}$ and set $d\gets d+1$.
\EndWhile
\State $\mathcal R_j\gets(\mathcal R_j\cap U_j^{\mathrm{split}})\cup\{\min_{i\in I}x_{ij}-0.001,\ \max_{i\in I}x_{ij}+0.001\}$ for all $j\in J$.
\State \Return $\{\mathcal R_j:j\in J\}$.
\end{algorithmic}
\end{breakablealgorithm}

At each level, every current subset is ordered by every feature. Whenever two consecutive value groups have a class transition, the values immediately before and after it are retained if they induce nontrivial splits of the subset, and each retained value generates the left and right subsets examined at the next level. The final sets are intersected with $U_j^{\mathrm{split}}$, and the two artificial values are added for structural compatibility rather than as nontrivial thresholds: deactivation is represented explicitly by $q=0$, and the upper value motivates the rightmost-leaf inequalities~\eqref{eq:leaf_label_rightmost} and is not selectable.

\begin{figure}[htbp]
\centering
\includegraphics[width=0.6\textwidth]{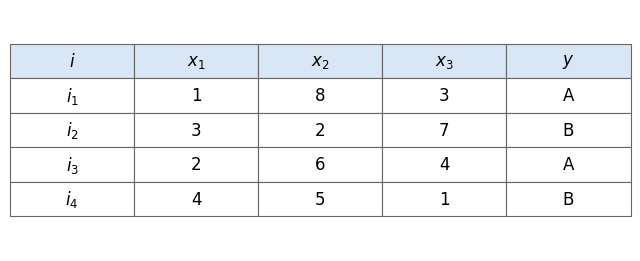}
\caption{Recursive Class-Boundary Candidate Split Construction with Three Features (Both Values Adjacent to Each Class Transition Are Retained When They Induce Nontrivial Splits, with the Artificial Lower and Upper Boundary Values)}
\label{fig:class_boundary_example}
\end{figure}

Figure~\ref{fig:class_boundary_example} illustrates the construction with four observations of classes A and B and three features. Sorting all observations by $x_1$ gives $i_1(A),i_3(A),i_2(B),i_4(B)$. The only transition lies between $i_3$ and $i_2$, so $x_{31}=2$ and $x_{21}=3$ are retained, and, with the artificial values $0.999$ and $4.001$ (the minimum and maximum of $x_1$ are $1$ and $4$), this ordering contributes $\{0.999,2,3,4.001\}$. Sorting by $x_2$ gives $i_2(B),i_4(B),i_3(A),i_1(A)$. The transition between $i_4$ and $i_3$ retains $x_{42}=5$ and $x_{32}=6$, and, with $1.999$ and $8.001$, this ordering contributes $\{1.999,5,6,8.001\}$. Sorting by $x_3$ gives $i_4(B),i_1(A),i_3(A),i_2(B)$, with transitions between $i_4$ and $i_1$ and between $i_3$ and $i_2$. The values $x_{43}=1$, $x_{13}=3$, and $x_{33}=4$ are retained, whereas $x_{23}=7$ is the maximum and induces no nontrivial split, so, with $0.999$ and $7.001$, this ordering contributes $\{0.999,1,3,4,7.001\}$. The initial-level candidate sets are therefore $\mathcal R_1^{(0)}=\{0.999,2,3,4.001\}$, $\mathcal R_2^{(0)}=\{1.999,5,6,8.001\}$, and $\mathcal R_3^{(0)}=\{0.999,1,3,4,7.001\}$. Each retained value also generates its left and right child subsets, which are reordered by every feature at the next level, and the final sets, the unique unions over all levels with the artificial values added, include any admissible values found in these subsets. A value that is not adjacent to a class transition in the complete ordering can thus be retained after an earlier partition restricts the observations. The lower artificial value represents the routing effect of deactivation, and the upper one motivates the inequalities~\eqref{eq:leaf_label_rightmost}.

\subsection{Convex Hull of the Depth-One MSSF Relaxation}
\label{convexhull}

We show that the linear relaxation of an extended depth-one, single-feature MSSF, which also permits the maximum observed threshold, describes the convex hull of its integer-feasible solutions under the unified $C_p$ formulation of Section~\ref{sec:mssf}. Deactivation is alternative $q=0$, routes every observation to the right leaf, and is included with the incompatible active splits in the left-routing $C_p$ inequalities. The routing-mass identities are derived from these inequalities rather than assumed, the $C_p$ constraints remain inequalities throughout and are not assumed to be tight for every tuple, and the proof constructs an explicit convex decomposition into integer-feasible depth-one trees.

Let $T_B=\{1\}$, $T_L=\{2,3\}$, and $J=\{1\}$, and suppress the branch-node and feature indices. Let $r_1<\cdots<r_R$ be all distinct observed feature values, used as admissible thresholds, so that $r_R=\max_{i\in I}x_i$; unlike the nontrivial split set of Section~\ref{sec:mssf}, $r_R$ may be selected. Let $b_q$ be the branching variable of alternative $q\in\{0,\ldots,R\}$, where $q=0$ sends every observation to leaf $3$. Let $L_q=\{i\in I:x_i\le r_q\}$ and $R_q=\{i\in I:x_i>r_q\}$ for $q=1,\ldots,R$, and $L_0=\varnothing$, $R_0=I$. The branch-partition and assignment constraints are
\begin{align}
\sum_{q=0}^R b_q&=1, \label{eq:d1_branch_partition_complete}\\
\sum_{k\in K}(z_{ik2}+z_{ik3})&=1,\qquad i\in I, \label{eq:d1_assignment_partition}
\end{align}
and all variables are nonnegative in the LP relaxation. Because the thresholds are the distinct observed values, the routing layers $D_1=L_1$ and $D_h=L_h\setminus L_{h-1}$ ($h=2,\ldots,R$) are nonempty and partition $I$; in particular, $D_{R+1}=I\setminus L_R=\varnothing$. For $h=1,\ldots,R$, let $S_h^L=\bigcup_{a=h}^RD_a=I\setminus L_{h-1}$, $S_h^R=\bigcup_{a=1}^hD_a=L_h$, $B_h^L=\sum_{q=h}^Rb_q$, and $B_h^R=\sum_{q=0}^{h-1}b_q$, so that $B_h^L+B_h^R=1$.

\paragraph{Depth-one $C_p$ inequalities.} The alternative set is $\mathcal Q=\{0,1,\ldots,R\}$. For $h=1,\ldots,R$, the left-exclusion set of alternative $h-1$ is $S_h^L$, and its dominance set contains alternatives $0,\ldots,h-2$, so, for every $(i_k)_{k\in K}\in(S_h^L)^{|K|}$, its left-exclusion inequality is
\begin{equation}
\sum_{q=0}^{h-1}b_q+\sum_{k\in K}z_{i_kk2}\le1;
\label{eq:d1_left_cp_packing}
\end{equation}
for $h=1$, this is the $q=0$ inequality $b_0+\sum_{k\in K}z_{i_kk2}\le1$. The right-exclusion set of alternative $h$ is $S_h^R$, and its dominance set contains alternatives $h+1,\ldots,R$, so, for every $(i_k)_{k\in K}\in(S_h^R)^{|K|}$,
\begin{equation}
\sum_{q=h}^{R}b_q+\sum_{k\in K}z_{i_kk3}\le1.
\label{eq:d1_right_cp_packing}
\end{equation}
Both families are instances of the unified set-packing construction applied to the extended alternative set and include all tuples, including the repeated-observation tuples used below; the rightmost-leaf inequalities~\eqref{eq:leaf_label_rightmost} are also imposed, and no equality is asserted for arbitrary tuples.

\begin{lemma}[Routing Equalities Implied by MSSF]
\label{lem:d1_routing_from_cp}
Every point satisfying \eqref{eq:d1_branch_partition_complete}--\eqref{eq:d1_right_cp_packing} satisfies, for $h=1,\ldots,R$ and every $i\in D_h$,
\begin{equation}
\sum_{k\in K}z_{ik2}=B_h^L,\qquad\sum_{k\in K}z_{ik3}=B_h^R.
\label{eq:d1_layer_routing_mass}
\end{equation}
\end{lemma}
\proof{Proof.}
Fix $i\in D_h$. Since $i\in S_h^L\cap S_h^R$, both families admit the tuple with $i_k=i$ for all $k\in K$, so \eqref{eq:d1_left_cp_packing} and \eqref{eq:d1_right_cp_packing} give $\sum_kz_{ik2}\le1-\sum_{q=0}^{h-1}b_q=B_h^L$ and $\sum_kz_{ik3}\le1-\sum_{q=h}^Rb_q=B_h^R$. These upper bounds sum to one, whereas \eqref{eq:d1_assignment_partition} requires the two left-hand sides to sum to one, so both bounds are attained. \Halmos
\endproof

\begin{lemma}[Left-Leaf Class Decomposition]
\label{lem:d1_left_class_decomposition}
For every feasible point $(b,z)$ of the extended depth-one MSSF LP relaxation, there exist $\alpha_{qk}\ge0$, $q=1,\ldots,R$, $k\in K$, such that $\sum_{k\in K}\alpha_{qk}=b_q$ for $q=1,\ldots,R$ and $z_{ik2}=\sum_{q\in\{1,\ldots,R\}:\,i\in L_q}\alpha_{qk}$ for $i\in I$ and $k\in K$.
\end{lemma}
\proof{Proof.}
Fix $h$ and $i,j\in D_h$. For $k\in K$, the tuple of \eqref{eq:d1_left_cp_packing} that selects $j$ for class $k$ and $i$ for every other class is admissible because $i,j\in S_h^L$, and, with Lemma~\ref{lem:d1_routing_from_cp}, it gives $z_{jk2}+\sum_{\ell\ne k}z_{i\ell2}\le B_h^L=\sum_{\ell\in K}z_{i\ell2}$, so $z_{jk2}\le z_{ik2}$; interchanging $i$ and $j$ gives equality. Hence $M_{hk}^L=z_{ik2}$ ($i\in D_h$) is well defined, and $\sum_kM_{hk}^L=B_h^L$. For $h<R$, $i\in D_h$, and $j\in D_{h+1}$, both observations belong to $S_h^L$, and the same tuple gives $M_{h+1,k}^L\le M_{hk}^L$. With $M_{R+1,k}^L=0$ and $B_{R+1}^L=0$, the coefficients $\alpha_{hk}=M_{hk}^L-M_{h+1,k}^L$ ($h=1,\ldots,R$) are nonnegative and satisfy $\sum_k\alpha_{hk}=B_h^L-B_{h+1}^L=b_h$, and, for $i\in D_h$, telescoping gives $z_{ik2}=M_{hk}^L=\sum_{q=h}^R\alpha_{qk}=\sum_{q:\,i\in L_q}\alpha_{qk}$. \Halmos
\endproof

\begin{lemma}[Right-Leaf Class Decomposition]
\label{lem:d1_right_class_decomposition}
For every feasible point $(b,z)$ of the extended depth-one MSSF LP relaxation, there exist $\beta_{qk}\ge0$, $q=0,\ldots,R$, $k\in K$, such that $\sum_{k\in K}\beta_{qk}=b_q$ for $q=0,\ldots,R$ and $z_{ik3}=\sum_{q\in\{0,\ldots,R\}:\,i\in R_q}\beta_{qk}$ for $i\in I$ and $k\in K$.
\end{lemma}
\proof{Proof.}
Fix $h$ and $i,j\in D_h$. Applying \eqref{eq:d1_right_cp_packing} to the tuple that selects $j$ for class $k$ and $i$ for every other class and using Lemma~\ref{lem:d1_routing_from_cp} gives $z_{jk3}+\sum_{\ell\ne k}z_{i\ell3}\le B_h^R=\sum_{\ell\in K}z_{i\ell3}$; interchanging $i$ and $j$ shows that all observations in $D_h$ share a right-leaf class vector $M_{hk}^R$ with $\sum_kM_{hk}^R=B_h^R$. For $h<R$, $i\in D_h$, and $j\in D_{h+1}$, both observations belong to $S_{h+1}^R$, and the inequality at index $h+1$ that selects $i$ for class $k$ and $j$ for every other class gives $z_{ik3}+\sum_{\ell\ne k}z_{j\ell3}\le B_{h+1}^R=\sum_{\ell\in K}z_{j\ell3}$, so $M_{hk}^R\le M_{h+1,k}^R$. Set $M_{0k}^R=0$. Because $r_R=\max_ix_i$, the right leaf is empty when alternative $R$ is selected; choose any $k^\star\in K$ and define the formal terminal vector $M_{R+1,k}^R=M_{Rk}^R+b_R\mathbf 1\{k=k^\star\}$, which is componentwise nondecreasing, has total mass one, and changes no assignment variable because $D_{R+1}=\varnothing$. Let $\beta_{0k}=M_{1k}^R$ and $\beta_{qk}=M_{q+1,k}^R-M_{qk}^R$ ($q=1,\ldots,R$). Then $\beta_{qk}\ge0$, $\sum_k\beta_{0k}=B_1^R=b_0$, and $\sum_k\beta_{qk}=B_{q+1}^R-B_q^R=b_q$ ($q=1,\ldots,R$), where $B_{R+1}^R=1$, and, for $i\in D_h$, $z_{ik3}=M_{hk}^R=\sum_{q=0}^{h-1}\beta_{qk}=\sum_{q:\,i\in R_q}\beta_{qk}$. \Halmos
\endproof

Lemma~\ref{lem:mssf_leaf_label_consistency} shows that the complete $C_p$ families enforce a unique class at every nonempty leaf of an integer-feasible solution; Lemmas~\ref{lem:d1_routing_from_cp}--\ref{lem:d1_right_class_decomposition} yield the convex-hull result.

\begin{theorem}[Convex Hull of the Depth-One MSSF Linear Relaxation]
\label{thm:mssf_d1_convex_hull_complete}
Consider the extended depth-one, single-feature MSSF formulation with the complete $C_p$ inequalities \eqref{eq:d1_left_cp_packing} and \eqref{eq:d1_right_cp_packing}, where deactivation is included as alternative $q=0$ in the left-routing dominance sets. Suppose the selectable active thresholds are all distinct observed feature values, including their maximum, and the rightmost-leaf inequalities~\eqref{eq:leaf_label_rightmost} are imposed. Let $P_{\mathrm{MSSF-D1}}$ denote its LP feasible region and $\mathcal F_{\mathrm{MSSF-D1}}^I$ its integer-feasible points. Then $P_{\mathrm{MSSF-D1}}=\operatorname{conv}(\mathcal F_{\mathrm{MSSF-D1}}^I)$; consequently, every extreme point of $P_{\mathrm{MSSF-D1}}$ is integral.
\end{theorem}
\proof{Proof.}
Every integer-feasible point belongs to the LP relaxation, so $\operatorname{conv}(\mathcal F_{\mathrm{MSSF-D1}}^I)\subseteq P_{\mathrm{MSSF-D1}}$. Conversely, fix $(b,z)\in P_{\mathrm{MSSF-D1}}$, and let $\alpha$ and $\beta$ be the arrays of Lemmas~\ref{lem:d1_left_class_decomposition} and~\ref{lem:d1_right_class_decomposition}. For $q=1,\ldots,R$ with $b_q>0$, let $\lambda_{qk_Lk_R}=\alpha_{qk_L}\beta_{qk_R}/b_q$ for $k_L,k_R\in K$, and set these coefficients to zero when $b_q=0$; for deactivation, let $\lambda_{0k_R}=\beta_{0k_R}$. The coefficients are nonnegative, $\sum_{k_L,k_R}\lambda_{qk_Lk_R}=b_q$ for $q\ge1$, and $\sum_{k_R}\lambda_{0k_R}=b_0$, so they sum to $\sum_{q=0}^Rb_q=1$. For $q\ge1$, let $(b^{q,k_L,k_R},z^{q,k_L,k_R})$ be the integer tree that selects $r_q$, routes $L_q$ to leaf $2$ with class $k_L$, and routes $R_q$ to leaf $3$ with class $k_R$ (the right class is irrelevant if $R_q=\varnothing$), and let $(b^{0,k_R},z^{0,k_R})$ deactivate the root and assign every observation to leaf $3$ with class $k_R$. Each of these points satisfies the extended partition constraints and the unified $C_p$ inequalities: if the selected alternative excludes a tuple from a leaf, its branching coordinate is one and all assignment coordinates of that tuple are zero, and, otherwise, the tuple contains at most one selected assignment because all observations at that leaf share one class. The convex combination
\begin{equation}
(\widehat b,\widehat z)=\sum_{q=1}^R\sum_{k_L,k_R\in K}\lambda_{qk_Lk_R}(b^{q,k_L,k_R},z^{q,k_L,k_R})+\sum_{k_R\in K}\lambda_{0k_R}(b^{0,k_R},z^{0,k_R})
\label{eq:d1_final_convex_decomposition}
\end{equation}
has branching coordinate $b_q$ for every $q$ by the marginal identities above, and, for $i\in I$ and $k\in K$,
\[
\widehat z_{ik2}=\sum_{q:\,i\in L_q}\sum_{k_R}\lambda_{qkk_R}=\sum_{q:\,i\in L_q}\alpha_{qk}=z_{ik2},\qquad
\widehat z_{ik3}=\lambda_{0k}+\sum_{q\ge1:\,i\in R_q}\sum_{k_L}\lambda_{qk_Lk}=\sum_{q:\,i\in R_q}\beta_{qk}=z_{ik3}.
\]
Hence $(b,z)=(\widehat b,\widehat z)\in\operatorname{conv}(\mathcal F_{\mathrm{MSSF-D1}}^I)$, which proves the theorem. \Halmos
\endproof

\section{Proofs of Section~\ref{sec:column_generation} Results}
\label{app:proofs_cg}

This section proves that assembling the selected subtree patterns preserves MSSF feasibility and the exact objective change, and characterizes the dual solution of the recentered master.

\proof{Proof of Proposition~\ref{prop:feasible_assembly}.}
Constraints~\eqref{eq:subtree_rmp_branch} and~\eqref{eq:subtree_rmp_assignment} imply that every branch node and every observation assignment is owned exactly once, either by its keep column or by one selected pattern. If $p\in\mathcal P_v$ and $q\in\mathcal P_w$ are both selected and $v$ is an ancestor of $w$, closure of $p$ implies that $p$ covers every branch row in $T_B(v)$, including $w$, which $q$ also covers, so the left-hand side of~\eqref{eq:subtree_rmp_branch} for $w$ is at least two, which is impossible; the same holds if $w$ is an ancestor of $v$. Hence the roots of simultaneously selected patterns are pairwise incomparable, and their subtrees are node-disjoint. Because all branching decisions above a replacement root equal the incumbent decisions, an observation entering one selected subtree cannot be rerouted into another, so the boundary sets $I_v^*$ and $I_w^*$ of simultaneously selected subtrees are disjoint, and the observations entering each replacement are exactly those in its boundary set. Within each selected subtree, all branching, deactivation, routing, and assignment decisions satisfy the MSSF constraints by admissibility of the pattern; outside these subtrees, the retained incumbent decisions are feasible because $(b^*,z^*)$ is MSSF-feasible. Consequently, no routing-dependent decision is assembled independently across the boundary of a changed subtree, and the assembled solution satisfies the MSSF branch partition, assignment partition, deactivation, and routing constraints, so it is integer-feasible for MSSF. \Halmos
\endproof

\begin{lemma}[Exact Objective Change]
\label{prop:exact_pattern_change}
Let $(b,z)$ be the MSSF solution assembled from an integer-feasible master solution, and let $C(b,z)$ and $C(b^*,z^*)$ denote the MSSF objective values of the assembled and incumbent trees. Then
\begin{equation}
C(b,z)-C(b^*,z^*)=\sum_{p\in\bar{\mathcal P}}\Delta C_p\lambda_p.
\label{eq:exact_master_change}
\end{equation}
\end{lemma}
\proof{Proof.}
Every branch row and observation row is covered exactly once. A keep column leaves the incumbent decision unchanged and contributes zero. A branch row $v$ covered by a replacement pattern contributes the coefficient of its new branching decision minus that of the incumbent decision at $v$, and an observation row $i$ covered by a replacement pattern contributes its new classification coefficient minus $c_i^{z,*}$. Because the partition equations prevent any row from being counted by two selected patterns, summing $\Delta C_p$ over the selected patterns counts every changed objective term exactly once and every unchanged term zero times, which gives~\eqref{eq:exact_master_change}. \Halmos
\endproof

\begin{corollary}[Improving Integer Master Solution]
\label{cor:negative_master_improves}
If an integer-feasible master solution satisfies $\sum_{p\in\bar{\mathcal P}}\Delta C_p\lambda_p<0$, then the assembled classification tree is feasible for MSSF and has a strictly smaller MSSF objective value than the incumbent tree.
\end{corollary}

\proof{Proof of Proposition~\ref{prop:zero_initial_duals}.}
When the initial restricted master contains only the zero-cost keep columns, one for each branch and observation partition row, every variable of its dual is bounded above by zero and the dual objective maximizes their sum. The zero vector is therefore feasible and strictly better than every solution with a negative component, so it is the unique optimal dual solution. For a candidate pattern $p$, the reduced cost $\bar C_p=\Delta C_p-\sum_{v\in T_B}a^b_{vp}\pi_v^*-\sum_{i\in I}a^z_{ip}\sigma_i^*$ then equals $\Delta C_p$. \Halmos
\endproof

\section{ISCG-OCT Algorithm}
\label{iscgalgo}

This section develops the recursive-pricing procedure that generates improving subtree patterns. Appendix~\ref{subsec:lower_bound_computation} constructs a lower bound for a partially fixed pricing state from the best assignment contribution compatible with the branching decisions fixed so far and a bound on the contributions of the unresolved branch nodes, Appendix~\ref{subsec:exact_assignment} proves that the assignment calculation is exact, and Appendix~\ref{subsec:iscg_search} embeds the bound in the recursive search, which generates branching and deactivation decisions recursively, prunes nonpromising partial states, and expands the surviving ones until a closed feasible subtree pattern is obtained. When root pricing is solved exactly with a zero improvement threshold over the complete predefined candidate family, the procedure determines whether an improving pattern exists within that family. An optimality certificate for the corresponding MSSF model additionally requires the applicable master-problem and termination conditions.

\subsection{Lower Bound Computation}
\label{subsec:lower_bound_computation}

After a branch-node decision is fixed, the routing restrictions determine which descendant leaves remain feasible for each observation, even though the descendant branching decisions are not yet fixed. For $i\in I_v$, let
\begin{equation}
\mathcal T_i=\bigl\{t\in T_L(v):\text{$i$ can still be assigned to $t$ under the routing decisions fixed so far}\bigr\}.
\label{eq:lb_feasible_leaf_set}
\end{equation}
An observation with $|\mathcal T_i|=1$ is forced to a single leaf, whereas one with $|\mathcal T_i|>1$ remains compatible with several leaves.

\begin{lemma}[Subtree Feasible-Leaf Property]
\label{lem:sp2_subtree_feasible_leaf}
For every $i\in I_v$ with $\mathcal T_i\neq\emptyset$, there exists a node $u\in T(v)$ such that $\mathcal T_i=T_L(u)$. Consequently, the nonempty sets $\mathcal T_i$ are laminar: $\mathcal T_i\cap\mathcal T_{i'}\in\{\emptyset,\mathcal T_i,\mathcal T_{i'}\}$ for all $i,i'\in I_v$.
\end{lemma}
\proof{Proof.}
Before any descendant routing decision is fixed, an observation entering $v$ can reach every leaf in $T_L(v)$. Suppose that the feasible leaves of $i$ are $T_L(u)$ for some $u\in T(v)$. A newly fixed decision that does not affect the subtree reachable by $i$ leaves this set unchanged; an active split fixed at $u$ sends $i$ uniquely to $\ell(u)$ or $r(u)$, so the set becomes $T_L(\ell(u))$ or $T_L(r(u))$. And deactivating $u$ sends every observation reaching $u$ to $r(u)$ and none to $\ell(u)$, so the set becomes $T_L(r(u))$. By induction, every nonempty feasible-leaf set is the descendant-leaf set of some node, and laminarity follows because descendant subtrees of a rooted binary tree are either disjoint or nested. \Halmos
\endproof

For each $u\in T(v)$, let $I_u^0=\{i\in I_v:\mathcal T_i=T_L(u)\}$ be the observations whose unresolved routing permits them to reach exactly the leaves of the subtree rooted at $u$. Each observation belongs to exactly one such set. To avoid overloading the pricing value function $V_v(I_v)$ of Section~\ref{sec:dynamic_pricing}, let $H_u(C)$ denote the minimum assignment objective change in the subtree rooted at $u$, conditional on the set of classes represented by its leaves being exactly the nonempty set $C\subseteq K$, and let
\begin{equation}
G_u(C)=\sum_{i\in I_u^0}\alpha\Bigl(\min_{k\in C}c^z_{ik}-c_i^{z,*}\Bigr).
\label{eq:lb_local_assignment}
\end{equation}
Because every observation in $I_u^0$ can reach every leaf in $T_L(u)$, it can use any class represented in that subtree, so $G_u(C)$ is its exact minimum assignment contribution given $C$. A leaf represents exactly one class, and the classes represented below an internal node $u\in T_B(v)$ are the union of the class sets $C_L$ and $C_R$ represented below $\ell(u)$ and $r(u)$. Hence
\begin{align}
H_t(\{k\})&=G_t(\{k\})=\sum_{i\in I_t^0}\alpha\bigl(c^z_{ik}-c_i^{z,*}\bigr),\qquad t\in T_L(v),\ k\in K, \label{eq:lb_leaf_assignment}\\
H_u(C)&=G_u(C)+\min_{\substack{\emptyset\neq C_L\subseteq K,\ \emptyset\neq C_R\subseteq K\\ C_L\cup C_R=C}}\bigl\{H_{\ell(u)}(C_L)+H_{r(u)}(C_R)\bigr\},\qquad u\in T_B(v), \label{eq:lb_internal_assignment}\\
L_v^z&=\min_{\emptyset\neq C\subseteq K}H_v(C). \label{eq:lb_root_assignment}
\end{align}

\begin{breakablealgorithm}
\caption{Hierarchical Assignment Lower-Bound Computation}
\label{alg:assignment_lower_bound}
\begin{algorithmic}[1]
\Require Observation set $I_v$; subtree $T(v)$ with leaves $T_L(v)$; feasible-leaf sets $\{\mathcal T_i:i\in I_v\}$; classes $K$; coefficients $\Delta c^z_{ikt}=\alpha(c^z_{ik}-c_i^{z,*})$.
\Ensure Assignment lower-bound value $L_v^z$.
\State Set $I_u^0\gets\emptyset$ for every $u\in T(v)$.
\For{each $i\in I_v$}
  \State Identify the node $u\in T(v)$ with $\mathcal T_i=T_L(u)$ and set $I_u^0\gets I_u^0\cup\{i\}$.
\EndFor
\For{each leaf $t\in T_L(v)$ and each class $k\in K$}
  \State Compute $H_t(\{k\})$ by~\eqref{eq:lb_leaf_assignment}.
\EndFor
\For{each internal node $u\in T_B(v)$ in bottom-up order and each nonempty $C\subseteq K$}
  \State Compute $G_u(C)$ by~\eqref{eq:lb_local_assignment} and $H_u(C)$ by~\eqref{eq:lb_internal_assignment}.
\EndFor
\State \Return $L_v^z=\min_{\emptyset\neq C\subseteq K}H_v(C)$.
\end{algorithmic}
\end{breakablealgorithm}

Algorithm~\ref{alg:assignment_lower_bound} exploits the hierarchical observation--leaf structure induced by partial routing. Observations whose feasible-leaf set is $T_L(u)$ are kept at $u$ rather than assigned prematurely to a descendant leaf, so they retain their flexibility until the class configuration of the subtree is determined. At a leaf, each common class is evaluated directly. At an internal node, the class sets of the two child subtrees are combined, and the observations that remain flexible throughout the subtree are evaluated against their union. Proposition~\ref{prop:exact_assignment} shows that $L_v^z$ is the minimum assignment objective change over the observation--leaf--class assignments compatible with the fixed routing decisions and the common-class requirement at each leaf.

Let $N_K=|K|$. Each node has at most $2^{N_K}-1$ nonempty represented-class sets. For each such set $C$ at an internal node, a direct implementation examines the pairs of nonempty child sets $(C_L,C_R)$ with $C_L\cup C_R=C$, which is bounded by $O(4^{N_K})$ operations per internal node, and evaluating $G_u(C)$ for all class sets requires $O(|I_u^0|2^{N_K})$ operations at node $u$, where $\sum_{u\in T(v)}|I_u^0|=|I_v|$. After the sets $\mathcal T_i$ and $I_u^0$ are constructed, Algorithm~\ref{alg:assignment_lower_bound} therefore runs in $O(|T(v)|4^{N_K}+|I_v|2^{N_K})$ time, and storing $H_u(C)$ requires $O(|T(v)|2^{N_K})$ memory. For binary classification, each node has only the three class sets $\{1\}$, $\{2\}$, and $\{1,2\}$, so the computation is $O(|T(v)|+|I_v|)$, linear in the numbers of nodes and observations, although the general multiclass recursion is exponential in the number of classes.

The assignment value is combined with the branching contributions. Let $F_v$ and $U_v=T_B(v)\setminus F_v$ be the fixed and unresolved branch nodes of the current partial subtree, and let $\Delta C_v^{B}$ be the objective contribution of the branching and deactivation decisions already fixed. Because deactivation is admissible at every unresolved node, the minimum possible branching contribution at $w\in U_v$ is $\underline{\Delta c}_w^b=\min\{\Delta c^b_{0w},\min_{(i,j)}\Delta c^b_{ijw}\}$, where the inner minimum ranges over admissible active decisions and is omitted if none exists. By the incumbent-relative coefficients, $\underline{\Delta c}_w^b=-2\beta$ if $w$ is active in the incumbent and $0$ otherwise. The complete partial-state lower bound is
\begin{equation}
LB_v(I_v)=\Delta C_v^{B}+\sum_{w\in U_v}\underline{\Delta c}_w^b+L_v^z=\Delta C_v^{B}-2\beta\bigl|\{w\in U_v:b^*_{0w}=0\}\bigr|+L_v^z.
\label{eq:lb_complete_pricing}
\end{equation}
The unresolved branching contributions are minimized independently, and the assignment term retains the relaxed routing flexibility, so the bound does not overestimate the objective change of any closed feasible completion (Proposition~\ref{prop:complete_pricing_lb}). Further fixing only restricts the set of completions, and pruning is valid whenever the corrected lower bound fails the required improvement threshold.

\subsection{Proofs of the Lower-Bound Results}
\label{subsec:exact_assignment}

\proof{Proof of Proposition~\ref{prop:exact_assignment}.}
We show by induction from the leaves toward $v$ that \eqref{eq:lb_leaf_assignment}--\eqref{eq:lb_root_assignment} give the minimum assignment objective change over all observation--leaf--class assignments that remain feasible under the fixed routing decisions, subject to the common-class requirement at each leaf; by Lemma~\ref{lem:sp2_subtree_feasible_leaf}, every observation belongs to exactly one set $I_u^0$. If class $k$ is selected at a leaf $t$, every observation in $I_t^0$ must use $k$ because $t$ is its only feasible leaf, and it contributes $\alpha(c^z_{ik}-c_i^{z,*})$, so \eqref{eq:lb_leaf_assignment} is exact. Let $u$ be an internal node, assume the result for $\ell(u)$ and $r(u)$, and fix a nonempty class set $C\subseteq K$. Any feasible assignment below $u$ that represents $C$ induces child class sets with $C_L\cup C_R=C$, and, by the induction hypothesis, the observations whose feasible-leaf sets lie within the child subtrees contribute at least $H_{\ell(u)}(C_L)+H_{r(u)}(C_R)$. Each $i\in I_u^0$ can reach every leaf in $T_L(u)$, so its contribution is at least $\min_{k\in C}\Delta c^z_{ikt}=\alpha(\min_{k\in C}c^z_{ik}-c_i^{z,*})$, which does not depend on the leaf $t$, and these observations contribute at least $G_u(C)$. Minimizing over $(C_L,C_R)$ shows that the right-hand side of~\eqref{eq:lb_internal_assignment} is a lower bound. Conversely, let $(C_L^*,C_R^*)$ attain the minimum in~\eqref{eq:lb_internal_assignment}. By the induction hypothesis, feasible assignments attaining $H_{\ell(u)}(C_L^*)$ and $H_{r(u)}(C_R^*)$ exist, and their union represents $C=C_L^*\cup C_R^*$. For each $i\in I_u^0$, a class $k_i\in\arg\min_{k\in C}c^z_{ik}$ is represented at some leaf of $T_L(u)$, which $i$ can reach because $\mathcal T_i=T_L(u)$, so $G_u(C)$ is attained. Hence \eqref{eq:lb_internal_assignment} is exact at every node, and \eqref{eq:lb_root_assignment}, which minimizes over every nonempty class set represented in the subtree rooted at $v$, is the minimum assignment objective change of the partially fixed state; this is the value $L_v^z$ returned by Algorithm~\ref{alg:assignment_lower_bound}. \Halmos
\endproof

\begin{proposition}[Validity of the Complete Pricing Lower Bound]
\label{prop:complete_pricing_lb}
For any partially fixed subtree-pricing state, the bound~\eqref{eq:lb_complete_pricing} satisfies $LB_v(I_v)\le\Delta C_p$ for every closed MSSF-feasible pattern $p$ extending that state.
\end{proposition}
\proof{Proof.}
Fix a feasible completion $p$. Its branching contribution at every $w\in F_v$ equals the contribution already included in $\Delta C_v^{B}$, and, at every $w\in U_v$, it is at least $\underline{\Delta c}_w^b$ by definition. Because the assignments induced by $p$ respect the partially fixed routing restrictions, Proposition~\ref{prop:exact_assignment} gives $L_v^z\le\Delta C_p^z$ for its assignment contribution $\Delta C_p^z$. Adding these bounds yields $LB_v(I_v)\le\Delta C_p$. \Halmos
\endproof

\subsection{Recursive Pricing}
\label{subsec:iscg_search}

Algorithm~\ref{alg:recursive_pricing} solves the subtree-pricing problem. The candidate threshold sets $\{\mathcal R_j\}_{j\in J}$ are predefined, may be constructed by any procedure, and specify the branching alternatives admitted by MSSF. Rather than enumerating every complete subtree configuration, the algorithm fixes one branch-node decision at a time and retains a partial configuration only when its lower bound, which includes the minimum possible contribution of every unresolved branch node, indicates that some completion can achieve the required improvement. This gives a branch-and-bound interpretation. Fixing a node creates a child search state, its lower bound bounds all its descendants, and an improvement threshold determines whether it is expanded.

\begin{breakablealgorithm}
\caption{ISCG-OCT}
\label{alg:recursive_pricing}
\begin{algorithmic}[1]
\Require Training data $(X,y)$; tree depth $N_D$; candidate threshold sets $\{\mathcal R_j\}_{j\in J}$; incumbent-relative costs $\Delta c^b$ and $\Delta c^z$; classification weight $\alpha$; improvement-threshold parameter $\rho\ge0$.
\Ensure A closed feasible subtree pattern $p\in\mathcal P_v$ and its objective change $\Delta C_p$.
\State Fit CART and use the resulting feasible tree as the initial incumbent; construct $\Delta c^b$ and $\Delta c^z$; set $UB\gets0$ and $\Gamma(\rho)\gets\rho N_I\alpha$.
\State Construct the branch-node processing order and initialize the pricing frontier at the root state $(v,I_v^*)$.
\For{each current branch state $(v,I_v)$ in the surviving frontier}
  \State Determine the observations $I_v$ reaching $v$ under the decisions fixed on the path to $v$; generate $\mathcal B_v(I_v)$ from all admissible alternatives in $\{\mathcal R_j\}_{j\in J}$, and include deactivation $b_{0v}=1$.
  \For{each admissible active split $(i,j)\in\mathcal B_v(I_v)$}
    \State Fix $b_{ijv}=1$, construct the child states $(\ell(v),I^l_{ijv}(I_v))$ and $(r(v),I^r_{ijv}(I_v))$, and compute $LB_v^{ij}(I_v)$ by~\eqref{eq:lb_complete_pricing}.
    \State Retain the split extension if $LB_v^{ij}(I_v)<UB-\Gamma(\rho)$; otherwise, prune it.
  \EndFor
  \State For $b_{0v}=1$, construct $(\ell(v),\emptyset)$ and $(r(v),I_v)$, compute $LB_v^0(I_v)$ by~\eqref{eq:lb_complete_pricing}, and retain the extension if $LB_v^0(I_v)<UB-\Gamma(\rho)$; otherwise, prune it.
  \State Recursively evaluate the surviving child states using $V_v^0(I_v)$, $Q_{ijv}(I_v)$, and $V_v(I_v)$ of Section~\ref{sec:dynamic_pricing}. Cache a value as exact only if all admissible continuations have been evaluated or valid bounds show that omitted continuations cannot improve it; otherwise, mark it as cutoff-dependent and do not cache it as exact.
\EndFor
\For{each terminal leaf state $(t,I_t)$ of a surviving completion}
  \State Set $V_t(I_t)\gets\min_{k\in K}\sum_{i\in I_t}\alpha(c^z_{ik}-c_i^{z,*})$, and let $k_t^*$ attain the minimum.
\EndFor
\State Recover the closed pattern $p\in\mathcal P_v$ from the selected decisions and leaf classes, evaluate its exact $\Delta C_p$, and return the minimum-value surviving pattern.
\State Add improving patterns to the RMP. If the integer RMP improves the feasible tree, update the incumbent, recenter $\Delta c^b$ and $\Delta c^z$, and repeat recursive pricing.
\end{algorithmic}
\end{breakablealgorithm}

The algorithm separates three quantities. The lower bound $LB_v(I_v)$ includes the fixed branching costs, the minimum contributions of all unresolved branch nodes, and the relaxed assignment cost, and it bounds every complete pattern extending the partial fixing (Proposition~\ref{prop:complete_pricing_lb}). The upper bound $UB=0$ is supplied by the feasible CART incumbent on the incumbent-centered scale. The threshold $\Gamma(\rho)=\rho N_I\alpha$ specifies how much improvement is required before a partial pattern is expanded. A tighter valid lower bound or feasible upper bound can be used without changing the search logic. The baseline implementation uses the bounds above.

With the admissible continuation family and the incumbent-relative coefficients fixed, the Bellman value $V_v(I_v)$ does not depend on decisions outside the subtree rooted at $v$, whereas the lower bound used for screening may include contributions fixed outside the local continuation; a local bound is therefore combined with these contributions before comparison with the global improvement threshold. A state is memoized as exact only when its minimum over all admissible continuations has been established, by exact evaluation or by bounds that rule out any better continuation. Cutoff-dependent results are kept separately and never reused as exact Bellman values, and exact memoized values are cleared whenever the incumbent-relative coefficients are updated.

\paragraph{Illustrative example.} We illustrate Algorithm~\ref{alg:recursive_pricing} on the eight observations in Table~\ref{tab:recursive_pricing_toy_data}, which are few enough that routing, candidate generation, screening, and state retention can be followed explicitly. All objective values are on the signed MSSF scale of Section~\ref{sec:mssf}.

\begin{table}[htbp]
\centering
\caption{Toy Binary-Classification Dataset}
\label{tab:recursive_pricing_toy_data}
\begin{tabular}{l*{8}{c}}
\toprule
Observation & 1 & 2 & 3 & 4 & 5 & 6 & 7 & 8\\
\midrule
$x_1$ & 1 & 1 & 2 & 2 & 4 & 4 & 5 & 5\\
$x_2$ & 1 & 4 & 2 & 5 & 1 & 4 & 2 & 5\\
Class & A & B & A & B & B & A & B & A\\
\bottomrule
\end{tabular}
\end{table}

For both features, consider the observed lower-endpoint candidate values $\{1,2,4\}$; adding the artificial values $u_j^{(0)}=\min_{i\in I}x_{ij}-0.001$ and $u_j^{(M_j+1)}=\max_{i\in I}x_{ij}+0.001$ gives $\mathcal R_1=\mathcal R_2=\{0.999,1,2,4,5.001\}$. The lower artificial value represents the routing effect of deactivation, and the upper one preserves the full right-routing case of the MSSF compatibility construction. Candidates are indexed by observed values rather than midpoints. A midpoint between consecutive values would induce the same partition of the training observations but is not used as the MSSF split index.

Let $\alpha=1/8$ and $\beta=0.05$, and let the incumbent CART tree consist of the root split $x_1\le2$ with left class A and right class B, which misclassifies observations $2$, $4$, $6$, and $8$. The incumbent has $E^*=4$ errors and $A^*=1$ active branch node, so, with the signed objective $C=\alpha(2E-N_I)+\beta(2A-|T_B|)$, $N_I=8$, and $|T_B|=3$, its value is $C^*=\frac18(2\cdot4-8)+0.05(2\cdot1-3)=-0.05$. After recentering at this incumbent, retaining it gives $UB=0$.

\emph{Step 1: Fix a root candidate and route observations.} The root candidate $x_1\le2$ routes $\{1,2,3,4\}$ to the left and $\{5,6,7,8\}$ to the right. Splitting each child at $x_2\le2$ separates the classes perfectly, giving $\{1,3\}$ (A), $\{2,4\}$ (B), $\{5,7\}$ (B), and $\{6,8\}$ (A), so a depth-two completion has zero classification error. With the optimistic classification completion, the fixed root branch, and the minimum contributions of the unresolved children, which are deactivated in the incumbent and therefore contribute zero, the lower bound of this state is $LB_1=-1.00$, the signed-MSSF counterpart of the conventional-scale value $-0.50$. For $\rho=0$, the threshold is zero, so $LB_1<0$ and the state remains eligible for expansion.

\emph{Step 2: Generate decisions at the next branch node.} At the left child, $I_v=\{1,2,3,4\}$, and admissible decisions are drawn from the predefined threshold sets, with routing evaluated on $I_v$. Consider $x_2\le1$, $x_2\le2$, and deactivation. For $x_2\le2$, the routing is $I^l_{ijv}(I_v)=\{1,3\}$ and $I^r_{ijv}(I_v)=\{2,4\}$, both leaves are pure, the right subtree of the root can also be completed with zero error, and $LB_v^{\,x_2\le2}=-0.90$. For $x_2\le1$, the routing is $\{1\}$ and $\{2,3,4\}$. The second set contains both classes, its best leaf assignment makes one error, and $LB_v^{\,x_2\le1}=-0.65$. Deactivation has the routing effect of the lower artificial value $u_2^{(0)}=0.999$: the left descendant receives no observations, all of $\{1,2,3,4\}$ continue to the right, their best common class misclassifies two observations, and, with the right side completed optimistically with zero error, $LB_v^{0}=-0.50$. All other unresolved branch nodes are deactivated in the incumbent and contribute zero, so these values already satisfy~\eqref{eq:lb_complete_pricing}. All three are negative and pass the $\rho=0$ test, so all three extensions remain in the surviving frontier.

\emph{Step 3: Complete the subtree pattern.} The same procedure is applied recursively to all surviving states. The continuation that applies $x_2\le2$ on both sides produces the pure leaves $\{1,3\}$ (A), $\{2,4\}$ (B), $\{5,7\}$ (B), and $\{6,8\}$ (A). The resulting closed depth-two tree has zero error and three active branch nodes, so $C_p=\frac18(2\cdot0-8)+0.05(2\cdot3-3)=-0.85$ and $\Delta C_p=C_p-C^*=-0.85-(-0.05)=-0.80<0$. The completed pattern is therefore an improving subtree column. The rule $LB<UB-\Gamma(\rho)$ determines which extensions are retained: at $\rho=0$, every extension with a negative lower bound is retained, and extensions that fail the rule are pruned.

\section{Computational Results}
\label{computational}

\newcommand{\streedNA}{--}
\newcommand{\streedOOM}{OOM}
\newcommand{\streedNote}{OOM: no tree returned within the memory limit.}

This section reports the dataset-level results summarized in Section~\ref{cresults}. ISCG-OCT with $U=0\%$ (exact) and $U=5\%$ (threshold) is compared with BNP-OCT, Benders OCT, and the original OCT formulation, all under a time limit of 3,600 seconds. Each table reports two tree depths.
The tables also report STreeD \citep{van2023necessary}, a dynamic-programming method for optimal decision trees on the same objective with $\alpha=0.99$ and $\beta=0.01$ and the same time limit. Because STreeD is exact for a given binarization of the features, every distinct training value of a continuous feature is supplied as a candidate threshold and every category of a categorical feature as a binary indicator, so a tree that STreeD certifies is optimal over all univariate splits. STreeD exceeded the available memory on Spambase, Breast Cancer Wisconsin (Diagnostic), Biodeg, and Ozone Level Detection (OOM in the tables).

\subsection{Training Accuracy}
\label{accuracycomputation}

Tables~\ref{tab:training_accuracy_D2} and~\ref{tab:training_accuracy_D4} report the training accuracy for $D=2,3,4,5$. The comparison examines whether recursive pricing maintains the solution quality of existing exact OCT methods. ISCG-OCT attains training accuracies comparable to, and frequently higher than, those of the competing methods. For shallow trees, several methods attain similar accuracies, and, as the depth increases, recursive pricing remains competitive on many datasets. The thresholded variant exhibits a speed--quality tradeoff on some instances. Its threshold is an improvement-screening parameter, and differences in training accuracy alone do not establish whether either variant has certified global optimality.

\begin{table}[!tbp]\centering
\caption{Training Accuracy for $D=2$ and $D=3$}
\label{tab:training_accuracy_D2}\label{tab:training_accuracy_D3}
{\def\baselinestretch{1}\fontsize{8}{9.5}\selectfont\renewcommand{\arraystretch}{0.9}\setlength{\tabcolsep}{3pt}\def\hd#1#2{\begin{tabular}[b]{@{}c@{}}#1\\#2\end{tabular}}
\begin{tabular*}{\textwidth}{@{\extracolsep{\fill}}l*{6}{r}@{}}
\toprule
 & \multicolumn{3}{c}{MIP-based} & & \multicolumn{2}{c}{ISCG-OCT} \\ \cmidrule(lr){2-4}\cmidrule(lr){6-7}
Dataset & \hd{BNP-}{OCT} & \hd{Benders}{OCT} & OCT & STreeD & $U=0\%$ & $U=5\%$\\
\midrule
\multicolumn{7}{@{}l}{\textit{(a) Depth $D=2$}}\\
\addlinespace[2pt]
Seismic-bumps & 0.94 & 0.94 & 0.93 & 0.94 & 0.94 & 0.94\\
Thoracic Surgery & 0.88 & 0.88 & 0.87 & 0.88 & 0.88 & 0.87\\
House-vote & 0.97 & 0.97 & 0.96 & 0.97 & 0.97 & 0.97\\
Mushroom & 0.98 & 0.98 & 0.97 & 0.96 & 0.98 & 0.94\\
Chess & 0.87 & 0.87 & 0.85 & 0.87 & 0.87 & 0.78\\
Nursery & 0.44 & 0.44 & 0.43 & 0.44 & 0.44 & 0.44\\
Cylinder Bands & 0.97 & 0.97 & 0.96 & 0.75 & 0.97 & 0.97\\
Echocardiogram & 0.97 & 0.97 & 0.96 & 0.97 & 0.97 & 0.97\\
Credit Approval & 0.87 & 0.87 & 0.88 & 0.87 & 0.88 & 0.87\\
Connectionist Bench (Sonar) & 0.86 & 0.83 & 0.82 & 0.87 & 0.86 & 0.81\\
Spambase & 0.84 & 0.82 & 0.83 & \streedOOM & 0.87 & 0.87\\
Blood Transfusion Service Center & 0.78 & 0.77 & 0.78 & 0.79 & 0.78 & 0.77\\
Breast Cancer Wisconsin (Diagnostic) & 0.92 & 0.95 & 0.95 & \streedOOM & 0.95 & 0.94\\
Wine & 0.96 & 0.96 & 0.96 & 0.96 & 0.96 & 0.94\\
Hepatitis & 0.88 & 0.88 & 0.88 & 0.90 & 0.88 & 0.82\\
Ozone Level Detection & 0.97 & 0.97 & 0.97 & \streedOOM & 0.97 & 0.97\\
Car-evaluation & 0.77 & 0.77 & 0.77 & 0.77 & 0.77 & 0.77\\
Balance-scale & 0.74 & 0.74 & 0.74 & 0.74 & 0.74 & 0.73\\
Biodeg & 0.81 & 0.81 & 0.79 & \streedOOM & 0.81 & 0.80\\
Iris & 0.96 & 0.96 & 0.96 & 0.94 & 0.96 & 0.96\\
Ionosphere & 0.91 & 0.91 & 0.91 & 0.91 & 0.91 & 0.90\\
Tic-tac-toe & 0.70 & 0.70 & 0.70 & 0.71 & 0.70 & 0.70\\
Banknote Authentication & 0.93 & 0.93 & 0.92 & 0.93 & 0.93 & 0.92\\
Fertility & 0.93 & 0.93 & 0.91 & 0.90 & 0.93 & 0.91\\
Planning Relax & 0.80 & 0.80 & 0.80 & 0.80 & 0.80 & 0.76\\
Seeds & 0.95 & 0.95 & 0.95 & 0.95 & 0.95 & 0.93\\
German & 0.73 & 0.74 & 0.74 & 0.74 & 0.74 & 0.70\\
\midrule
\multicolumn{7}{@{}l}{\textit{(b) Depth $D=3$}}\\
\addlinespace[2pt]
Seismic-bumps & 0.94 & 0.94 & 0.93 & 0.94 & 0.94 & 0.94\\
Thoracic Surgery & 0.89 & 0.87 & 0.88 & 0.90 & 0.89 & 0.88\\
House-vote & 0.97 & 0.98 & 0.96 & 0.98 & 0.98 & 0.98\\
Mushroom & 0.96 & 0.98 & 0.93 & 1.00 & 0.98 & 0.98\\
Chess & 0.75 & 0.88 & 0.75 & 0.94 & 0.90 & 0.90\\
Nursery & 0.43 & 0.44 & 0.43 & 0.44 & 0.44 & 0.44\\
Cylinder Bands & 0.97 & 0.98 & 0.96 & 0.82 & 0.99 & 0.98\\
Echocardiogram & 0.97 & 0.98 & 0.96 & 1.00 & 0.98 & 0.98\\
Credit Approval & 0.87 & 0.87 & 0.89 & 0.89 & 0.87 & 0.87\\
Connectionist Bench (Sonar) & 0.87 & 0.87 & 0.91 & 0.97 & 0.89 & 0.89\\
Spambase & 0.62 & 0.73 & 0.71 & \streedOOM & 0.88 & 0.88\\
Blood Transfusion Service Center & 0.79 & 0.78 & 0.78 & 0.82 & 0.79 & 0.78\\
Breast Cancer Wisconsin (Diagnostic) & 0.97 & 0.97 & 0.97 & \streedOOM & 0.97 & 0.97\\
Wine & 0.99 & 1.00 & 1.00 & 1.00 & 1.00 & 0.99\\
Hepatitis & 0.95 & 0.95 & 0.93 & 0.96 & 0.95 & 0.86\\
Ozone Level Detection & 0.97 & 0.97 & 0.97 & \streedOOM & 0.97 & 0.97\\
Car-evaluation & 0.81 & 0.76 & 0.81 & 0.81 & 0.81 & 0.78\\
Balance-scale & 0.80 & 0.77 & 0.78 & 0.80 & 0.80 & 0.77\\
Biodeg & 0.79 & 0.85 & 0.80 & \streedOOM & 0.85 & 0.85\\
Iris & 1.00 & 1.00 & 1.00 & 0.99 & 1.00 & 1.00\\
Ionosphere & 0.92 & 0.92 & 0.94 & 0.96 & 0.93 & 0.93\\
Tic-tac-toe & 0.77 & 0.76 & 0.76 & 0.78 & 0.77 & 0.71\\
Banknote Authentication & 0.97 & 0.97 & 0.97 & 0.99 & 0.94 & 0.94\\
Fertility & 0.97 & 0.97 & 0.97 & 0.94 & 0.97 & 0.94\\
Planning Relax & 0.87 & 0.74 & 0.86 & 0.88 & 0.87 & 0.87\\
Seeds & 0.99 & 0.99 & 0.99 & 0.99 & 0.99 & 0.97\\
German & 0.75 & 0.75 & 0.74 & 0.78 & 0.75 & 0.75\\
\bottomrule
\end{tabular*}

\vspace{2pt}\begin{minipage}{\textwidth}\fontsize{8}{9.5}\selectfont\textit{Note.} Time limit of 3,600 seconds. \streedNote\end{minipage}}
\end{table}

\begin{table}[!tbp]\centering
\caption{Training Accuracy for $D=4$ and $D=5$}
\label{tab:training_accuracy_D4}\label{tab:training_accuracy_D5}
{\def\baselinestretch{1}\fontsize{8}{9.5}\selectfont\renewcommand{\arraystretch}{0.9}\setlength{\tabcolsep}{3pt}\def\hd#1#2{\begin{tabular}[b]{@{}c@{}}#1\\#2\end{tabular}}
\begin{tabular*}{\textwidth}{@{\extracolsep{\fill}}l*{6}{r}@{}}
\toprule
 & \multicolumn{3}{c}{MIP-based} & & \multicolumn{2}{c}{ISCG-OCT} \\ \cmidrule(lr){2-4}\cmidrule(lr){6-7}
Dataset & \hd{BNP-}{OCT} & \hd{Benders}{OCT} & OCT & STreeD & $U=0\%$ & $U=5\%$\\
\midrule
\multicolumn{7}{@{}l}{\textit{(a) Depth $D=4$}}\\
\addlinespace[2pt]
Seismic-bumps & 0.94 & 0.94 & 0.94 & 0.94 & 0.94 & 0.94\\
Thoracic Surgery & 0.89 & 0.88 & 0.89 & 0.93 & 0.89 & 0.88\\
House-vote & 0.99 & 1.00 & 0.99 & 1.00 & 1.00 & 0.98\\
Mushroom & 0.81 & 0.93 & 0.85 & 1.00 & 0.98 & 0.98\\
Chess & 0.70 & 0.78 & 0.69 & 0.96 & 0.94 & 0.94\\
Nursery & 0.41 & 0.34 & 0.34 & 0.45 & 0.44 & 0.44\\
Cylinder Bands & 0.99 & 1.00 & 0.99 & 0.89 & 1.00 & 0.98\\
Echocardiogram & 0.99 & 1.00 & 0.99 & 1.00 & 1.00 & 0.98\\
Credit Approval & 0.87 & 0.87 & 0.88 & 0.93 & 0.90 & 0.90\\
Connectionist Bench (Sonar) & 0.79 & 0.92 & 0.84 & 0.97 & 0.97 & 0.97\\
Spambase & 0.63 & 0.76 & 0.76 & \streedOOM & 0.91 & 0.91\\
Blood Transfusion Service Center & 0.78 & 0.77 & 0.78 & 0.85 & 0.79 & 0.79\\
Breast Cancer Wisconsin (Diagnostic) & 0.92 & 0.92 & 0.97 & \streedOOM & 0.99 & 0.99\\
Wine & 1.00 & 1.00 & 1.00 & 1.00 & 1.00 & 1.00\\
Hepatitis & 0.98 & 0.98 & 0.98 & 1.00 & 0.93 & 0.93\\
Ozone Level Detection & 0.97 & 0.97 & 0.97 & \streedOOM & 0.98 & 0.98\\
Car-evaluation & 0.83 & 0.75 & 0.75 & 0.84 & 0.83 & 0.82\\
Balance-scale & 0.84 & 0.84 & 0.84 & 0.86 & 0.84 & 0.84\\
Biodeg & 0.76 & 0.80 & 0.67 & \streedOOM & 0.87 & 0.87\\
Iris & 1.00 & 1.00 & 1.00 & 1.00 & 1.00 & 1.00\\
Ionosphere & 0.82 & 0.83 & 0.96 & 0.96 & 0.95 & 0.95\\
Tic-tac-toe & 0.82 & 0.82 & 0.75 & 0.87 & 0.83 & 0.77\\
Banknote Authentication & 0.86 & 0.87 & 0.97 & 1.00 & 0.96 & 0.96\\
Fertility & 0.86 & 0.96 & 0.97 & 1.00 & 0.96 & 0.96\\
Planning Relax & 0.87 & 0.88 & 0.87 & 0.98 & 0.80 & 0.80\\
Seeds & 0.99 & 0.99 & 0.99 & 1.00 & 0.98 & 0.98\\
German & 0.75 & 0.70 & 0.76 & 0.82 & 0.76 & 0.76\\
\midrule
\multicolumn{7}{@{}l}{\textit{(b) Depth $D=5$}}\\
\addlinespace[2pt]
Seismic-bumps & 0.94 & 0.94 & 0.94 & 0.94 & 0.95 & 0.95\\
Thoracic Surgery & 0.88 & 0.88 & 0.88 & 0.98 & 0.90 & 0.90\\
House-vote & 0.99 & 1.00 & 1.00 & 1.00 & 1.00 & 0.99\\
Mushroom & 0.93 & 0.93 & 0.93 & 1.00 & 1.00 & 1.00\\
Chess & 0.78 & 0.87 & 0.87 & 0.98 & 0.94 & 0.94\\
Nursery & 0.35 & 0.34 & 0.34 & 0.47 & 0.45 & 0.45\\
Cylinder Bands & 0.99 & 1.00 & 1.00 & 0.89 & 1.00 & 0.99\\
Echocardiogram & 0.99 & 1.00 & 1.00 & 1.00 & 1.00 & 0.99\\
Credit Approval & 0.87 & 0.87 & 0.87 & 0.93 & 0.92 & 0.92\\
Connectionist Bench (Sonar) & 0.78 & 0.77 & 0.52 & 0.97 & 1.00 & 1.00\\
Spambase & 0.62 & 0.85 & 0.62 & \streedOOM & 0.93 & 0.93\\
Blood Transfusion Service Center & 0.78 & 0.77 & 0.77 & 0.88 & 0.80 & 0.80\\
Breast Cancer Wisconsin (Diagnostic) & 0.92 & 0.63 & 0.63 & \streedOOM & 0.99 & 0.99\\
Wine & 0.93 & 1.00 & 0.40 & 1.00 & 1.00 & 1.00\\
Hepatitis & 0.98 & 0.98 & 0.55 & 1.00 & 0.95 & 0.95\\
Ozone Level Detection & 0.97 & 0.97 & 0.97 & \streedOOM & 0.99 & 0.99\\
Car-evaluation & 0.82 & 0.70 & 0.70 & 0.89 & 0.86 & 0.86\\
Balance-scale & 0.84 & 0.83 & 0.82 & 0.90 & 0.86 & 0.86\\
Biodeg & 0.73 & 0.66 & 0.66 & \streedOOM & 0.91 & 0.91\\
Iris & 1.00 & 1.00 & 1.00 & 1.00 & 1.00 & 1.00\\
Ionosphere & 0.90 & 0.90 & 0.75 & 0.96 & 0.96 & 0.96\\
Tic-tac-toe & 0.80 & 0.80 & 0.65 & 0.95 & 0.83 & 0.83\\
Banknote Authentication & 0.90 & 0.90 & 0.56 & 1.00 & 0.98 & 0.98\\
Fertility & 0.90 & 0.90 & 0.90 & 1.00 & 0.98 & 0.98\\
Planning Relax & 0.77 & 0.74 & 0.74 & 0.98 & 0.86 & 0.86\\
Seeds & 0.99 & 0.98 & 0.98 & 1.00 & 1.00 & 1.00\\
German & 0.74 & 0.70 & 0.70 & 0.86 & 0.79 & 0.79\\
\bottomrule
\end{tabular*}

\vspace{2pt}\begin{minipage}{\textwidth}\fontsize{8}{9.5}\selectfont\textit{Note.} Time limit of 3,600 seconds. \streedNote\end{minipage}}
\end{table}

\subsection{Testing Accuracy}
\label{testcomputation}

Tables~\ref{tab:testing_accuracy} and~\ref{tab:testing_accuracy_d4} report the testing accuracy on the common 70/30 train--test split. Although the framework is designed to improve optimization efficiency rather than predictive generalization, these results provide a holdout comparison with existing exact OCT methods. ISCG-OCT is competitive on many, but not all, datasets while often requiring less computation, and no single algorithm dominates, because better optimization of the training objective need not yield higher testing accuracy. Because the results are based on a single partition, small differences should be read as descriptive of this holdout sample rather than as evidence of systematic differences in out-of-sample performance.

\begin{table}[!tbp]\centering
\caption{Testing Accuracy for $D=2$ and $D=3$}
\label{tab:testing_accuracy}\label{tab:testing_accuracy_d3}
{\def\baselinestretch{1}\fontsize{8}{9.5}\selectfont\renewcommand{\arraystretch}{0.9}\setlength{\tabcolsep}{3pt}\def\hd#1#2{\begin{tabular}[b]{@{}c@{}}#1\\#2\end{tabular}}
\begin{tabular*}{\textwidth}{@{\extracolsep{\fill}}l*{6}{r}@{}}
\toprule
 & \multicolumn{3}{c}{MIP-based} & & \multicolumn{2}{c}{ISCG-OCT} \\ \cmidrule(lr){2-4}\cmidrule(lr){6-7}
Dataset & \hd{BNP-}{OCT} & \hd{Benders}{OCT} & OCT & STreeD & $U=0\%$ & $U=5\%$\\
\midrule
\multicolumn{7}{@{}l}{\textit{(a) Depth $D=2$}}\\
\addlinespace[2pt]
Seismic-bumps & 0.93 & 0.93 & 0.93 & 0.93 & 0.92 & 0.93\\
Thoracic Surgery & 0.82 & 0.84 & 0.84 & 0.79 & 0.82 & 0.82\\
House-vote & 0.97 & 0.97 & 0.99 & 0.97 & 0.97 & 0.97\\
Mushroom & 0.98 & 0.89 & 0.97 & 0.96 & 0.87 & 0.94\\
Chess & 0.88 & 0.88 & 0.84 & 0.88 & 0.88 & 0.78\\
Nursery & 0.42 & 0.42 & 0.44 & 0.42 & 0.42 & 0.42\\
Cylinder Bands & 0.97 & 0.97 & 0.99 & 0.72 & 0.97 & 0.97\\
Echocardiogram & 0.97 & 0.97 & 1.00 & 0.85 & 0.97 & 0.97\\
Credit Approval & 0.85 & 0.85 & 0.82 & 0.83 & 0.85 & 0.86\\
Connectionist Bench (Sonar) & 0.70 & 0.71 & 0.79 & 0.79 & 0.71 & 0.75\\
Spambase & 0.82 & 0.81 & 0.81 & \streedOOM & 0.85 & 0.85\\
Blood Transfusion Service Center & 0.74 & 0.74 & 0.74 & 0.71 & 0.76 & 0.74\\
Breast Cancer Wisconsin (Diagnostic) & 0.95 & 0.92 & 0.93 & \streedOOM & 0.95 & 0.93\\
Wine & 0.98 & 0.98 & 0.89 & 0.96 & 0.98 & 0.87\\
Hepatitis & 0.58 & 0.67 & 0.54 & 0.74 & 0.72 & 0.63\\
Ozone Level Detection & 0.97 & 0.97 & 0.96 & \streedOOM & 0.97 & 0.97\\
Car-evaluation & 0.76 & 0.76 & 0.76 & 0.76 & 0.76 & 0.76\\
Balance-scale & 0.66 & 0.66 & 0.66 & 0.66 & 0.65 & 0.67\\
Biodeg & 0.80 & 0.74 & 0.78 & \streedOOM & 0.75 & 0.75\\
Iris & 0.91 & 0.91 & 0.91 & 0.98 & 0.91 & 0.91\\
Ionosphere & 0.92 & 0.92 & 0.92 & 0.92 & 0.92 & 0.92\\
Tic-tac-toe & 0.69 & 0.67 & 0.69 & 0.69 & 0.67 & 0.69\\
Banknote Authentication & 0.92 & 0.91 & 0.90 & 0.92 & 0.91 & 0.91\\
Fertility & 0.93 & 0.93 & 0.93 & 0.87 & 0.93 & 0.90\\
Planning Relax & 0.58 & 0.65 & 0.53 & 0.58 & 0.65 & 0.65\\
Seeds & 0.94 & 0.94 & 0.87 & 0.89 & 0.94 & 0.87\\
German & 0.71 & 0.70 & 0.72 & 0.72 & 0.70 & 0.70\\
\midrule
\multicolumn{7}{@{}l}{\textit{(b) Depth $D=3$}}\\
\addlinespace[2pt]
Seismic-bumps & 0.93 & 0.93 & 0.93 & 0.93 & 0.92 & 0.93\\
Thoracic Surgery & 0.84 & 0.84 & 0.84 & 0.79 & 0.81 & 0.80\\
House-vote & 0.97 & 0.96 & 0.99 & 0.96 & 0.93 & 0.91\\
Mushroom & 0.96 & 0.93 & 0.94 & 1.00 & 0.80 & 0.96\\
Chess & 0.72 & 0.86 & 0.75 & 0.93 & 0.80 & 0.91\\
Nursery & 0.43 & 0.41 & 0.44 & 0.41 & 0.41 & 0.41\\
Cylinder Bands & 0.97 & 0.96 & 0.99 & 0.69 & 0.97 & 0.91\\
Echocardiogram & 0.97 & 0.96 & 1.00 & 0.80 & 0.97 & 0.91\\
Credit Approval & 0.85 & 0.85 & 0.83 & 0.83 & 0.83 & 0.83\\
Connectionist Bench (Sonar) & 0.65 & 0.70 & 0.71 & 0.73 & 0.73 & 0.71\\
Spambase & 0.59 & 0.58 & 0.66 & \streedOOM & 0.89 & 0.86\\
Blood Transfusion Service Center & 0.75 & 0.74 & 0.74 & 0.72 & 0.76 & 0.75\\
Breast Cancer Wisconsin (Diagnostic) & 0.94 & 0.94 & 0.95 & \streedOOM & 0.94 & 0.96\\
Wine & 0.96 & 0.93 & 0.94 & 0.91 & 0.96 & 0.96\\
Hepatitis & 0.58 & 0.67 & 0.54 & 0.72 & 0.72 & 0.63\\
Ozone Level Detection & 0.97 & 0.97 & 0.97 & \streedOOM & 0.96 & 0.97\\
Car-evaluation & 0.82 & 0.82 & 0.82 & 0.81 & 0.83 & 0.76\\
Balance-scale & 0.70 & 0.78 & 0.69 & 0.71 & 0.78 & 0.70\\
Biodeg & 0.76 & 0.81 & 0.81 & \streedOOM & 0.84 & 0.81\\
Iris & 0.96 & 0.96 & 0.96 & 1.00 & 0.96 & 0.96\\
Ionosphere & 0.86 & 0.88 & 0.92 & 0.89 & 0.91 & 0.90\\
Tic-tac-toe & 0.76 & 0.76 & 0.77 & 0.74 & 0.76 & 0.69\\
Banknote Authentication & 0.96 & 0.94 & 0.97 & 0.96 & 0.94 & 0.93\\
Fertility & 0.95 & 0.93 & 0.97 & 0.80 & 0.93 & 0.93\\
Planning Relax & 0.60 & 0.65 & 0.60 & 0.55 & 0.56 & 0.65\\
Seeds & 0.90 & 0.88 & 0.87 & 0.92 & 0.89 & 0.86\\
German & 0.71 & 0.70 & 0.71 & 0.71 & 0.70 & 0.70\\
\bottomrule
\end{tabular*}

\vspace{2pt}\begin{minipage}{\textwidth}\fontsize{8}{9.5}\selectfont\textit{Note.} Time limit of 3,600 seconds. \streedNote\end{minipage}}
\end{table}

\begin{table}[!tbp]\centering
\caption{Testing Accuracy for $D=4$ and $D=5$}
\label{tab:testing_accuracy_d4}\label{tab:testing_accuracy_d5}
{\def\baselinestretch{1}\fontsize{8}{9.5}\selectfont\renewcommand{\arraystretch}{0.9}\setlength{\tabcolsep}{3pt}\def\hd#1#2{\begin{tabular}[b]{@{}c@{}}#1\\#2\end{tabular}}
\begin{tabular*}{\textwidth}{@{\extracolsep{\fill}}l*{6}{r}@{}}
\toprule
 & \multicolumn{3}{c}{MIP-based} & & \multicolumn{2}{c}{ISCG-OCT} \\ \cmidrule(lr){2-4}\cmidrule(lr){6-7}
Dataset & \hd{BNP-}{OCT} & \hd{Benders}{OCT} & OCT & STreeD & $U=0\%$ & $U=5\%$\\
\midrule
\multicolumn{7}{@{}l}{\textit{(a) Depth $D=4$}}\\
\addlinespace[2pt]
Seismic-bumps & 0.93 & 0.93 & 0.93 & 0.93 & 0.92 & 0.93\\
Thoracic Surgery & 0.83 & 0.84 & 0.80 & 0.79 & 0.81 & 0.81\\
House-vote & 0.96 & 0.96 & 0.99 & 0.96 & 0.96 & 0.91\\
Mushroom & 0.80 & 0.93 & 0.86 & 1.00 & 0.89 & 0.97\\
Chess & 0.67 & 0.78 & 0.67 & 0.95 & 0.80 & 0.93\\
Nursery & 0.42 & 0.32 & 0.32 & 0.40 & 0.41 & 0.41\\
Cylinder Bands & 0.96 & 0.96 & 0.99 & 0.68 & 0.97 & 0.91\\
Echocardiogram & 0.96 & 0.96 & 1.00 & 0.80 & 0.97 & 0.91\\
Credit Approval & 0.85 & 0.85 & 0.84 & 0.81 & 0.83 & 0.83\\
Connectionist Bench (Sonar) & 0.65 & 0.56 & 0.65 & 0.73 & 0.73 & 0.73\\
Spambase & 0.61 & 0.58 & 0.76 & \streedOOM & 0.89 & 0.89\\
Blood Transfusion Service Center & 0.74 & 0.74 & 0.74 & 0.73 & 0.74 & 0.74\\
Breast Cancer Wisconsin (Diagnostic) & 0.89 & 0.63 & 0.96 & \streedOOM & 0.94 & 0.94\\
Wine & 0.93 & 0.39 & 0.94 & 0.91 & 0.96 & 0.96\\
Hepatitis & 0.63 & 0.67 & 0.63 & 0.74 & 0.42 & 0.42\\
Ozone Level Detection & 0.97 & 0.97 & 0.97 & \streedOOM & 0.96 & 0.96\\
Car-evaluation & 0.82 & 0.69 & 0.72 & 0.84 & 0.83 & 0.83\\
Balance-scale & 0.74 & 0.77 & 0.74 & 0.76 & 0.78 & 0.78\\
Biodeg & 0.77 & 0.66 & 0.66 & \streedOOM & 0.84 & 0.84\\
Iris & 0.96 & 0.96 & 0.96 & 0.89 & 0.96 & 0.96\\
Ionosphere & 0.79 & 0.75 & 0.88 & 0.89 & 0.91 & 0.91\\
Tic-tac-toe & 0.80 & 0.76 & 0.76 & 0.81 & 0.78 & 0.76\\
Banknote Authentication & 0.84 & 0.92 & 0.97 & 0.99 & 0.94 & 0.94\\
Fertility & 0.84 & 0.95 & 0.97 & 0.77 & 0.93 & 0.93\\
Planning Relax & 0.58 & 0.65 & 0.55 & 0.62 & 0.56 & 0.56\\
Seeds & 0.89 & 0.89 & 0.90 & 0.90 & 0.89 & 0.89\\
German & 0.73 & 0.69 & 0.73 & 0.73 & 0.70 & 0.70\\
\midrule
\multicolumn{7}{@{}l}{\textit{(b) Depth $D=5$}}\\
\addlinespace[2pt]
Seismic-bumps & 0.93 & 0.93 & 0.93 & 0.93 & 0.92 & 0.92\\
Thoracic Surgery & 0.82 & 0.84 & 0.85 & 0.79 & 0.82 & 0.82\\
House-vote & 0.97 & 0.96 & 0.99 & 0.94 & 0.97 & 0.91\\
Mushroom & 0.93 & 0.93 & 0.90 & 1.00 & 0.73 & 0.94\\
Chess & 0.78 & 0.88 & 0.67 & 0.96 & 0.93 & 0.93\\
Nursery & 0.34 & 0.32 & 0.34 & 0.40 & 0.41 & 0.41\\
Cylinder Bands & 0.97 & 0.96 & 0.99 & 0.68 & 0.97 & 0.91\\
Echocardiogram & 0.97 & 0.96 & 1.00 & 0.80 & 0.97 & 0.91\\
Credit Approval & 0.85 & 0.85 & 0.84 & 0.81 & 0.83 & 0.83\\
Connectionist Bench (Sonar) & 0.65 & 0.56 & 0.71 & 0.73 & 0.73 & 0.73\\
Spambase & 0.88 & 0.58 & 0.76 & \streedOOM & 0.90 & 0.90\\
Blood Transfusion Service Center & 0.74 & 0.74 & 0.74 & 0.69 & 0.74 & 0.74\\
Breast Cancer Wisconsin (Diagnostic) & 0.89 & 0.93 & 0.93 & \streedOOM & 0.96 & 0.96\\
Wine & 0.83 & 0.84 & 0.94 & 0.91 & 0.96 & 0.96\\
Hepatitis & 0.67 & 0.67 & 0.54 & 0.68 & 0.58 & 0.58\\
Ozone Level Detection & 0.97 & 0.97 & 0.97 & \streedOOM & 0.95 & 0.95\\
Car-evaluation & 0.82 & 0.69 & 0.69 & 0.86 & 0.86 & 0.86\\
Balance-scale & 0.75 & 0.75 & 0.79 & 0.74 & 0.72 & 0.72\\
Biodeg & 0.73 & 0.70 & 0.66 & \streedOOM & 0.83 & 0.83\\
Iris & 0.96 & 0.96 & 0.96 & 0.89 & 0.96 & 0.96\\
Ionosphere & 0.88 & 0.90 & 0.92 & 0.89 & 0.90 & 0.90\\
Tic-tac-toe & 0.76 & 0.67 & 0.70 & 0.89 & 0.81 & 0.81\\
Banknote Authentication & 0.89 & 0.93 & 0.93 & 0.99 & 0.96 & 0.96\\
Fertility & 0.89 & 0.93 & 0.93 & 0.90 & 0.97 & 0.97\\
Planning Relax & 0.67 & 0.65 & 0.45 & 0.62 & 0.47 & 0.47\\
Seeds & 0.90 & 0.90 & 0.89 & 0.90 & 0.86 & 0.86\\
German & 0.71 & 0.69 & 0.70 & 0.71 & 0.70 & 0.70\\
\bottomrule
\end{tabular*}

\vspace{2pt}\begin{minipage}{\textwidth}\fontsize{8}{9.5}\selectfont\textit{Note.} Time limit of 3,600 seconds. \streedNote\end{minipage}}
\end{table}

\subsection{Computational Times}
\label{computationaltime}

Tables~\ref{tab:comp_time_D2} and~\ref{tab:comp_time_D4} report wall-clock times, measured for each method from method-specific preprocessing and model construction through termination of the optimization procedure. ISCG-OCT often terminates substantially sooner than the benchmark formulations, although several dataset--depth combinations reach the time limit. The thresholded variant can further reduce running time by screening recursive states, but it is not uniformly faster than the exact variant. For STreeD, the time is measured from the binarization of the features through the end of the search, and runs that stop at the time limit are reported at 3,600 seconds.

\begin{table}[!tbp]\centering
\caption{Running Times (Seconds) for $D=2$ and $D=3$}
\label{tab:comp_time_D2}\label{tab:comp_time_D3}
{\def\baselinestretch{1}\fontsize{8}{9.5}\selectfont\renewcommand{\arraystretch}{0.9}\setlength{\tabcolsep}{3pt}\def\hd#1#2{\begin{tabular}[b]{@{}c@{}}#1\\#2\end{tabular}}
\begin{tabular*}{\textwidth}{@{\extracolsep{\fill}}l*{6}{r}@{}}
\toprule
 & \multicolumn{3}{c}{MIP-based} & & \multicolumn{2}{c}{ISCG-OCT} \\ \cmidrule(lr){2-4}\cmidrule(lr){6-7}
Dataset & \hd{BNP-}{OCT} & \hd{Benders}{OCT} & OCT & STreeD & $U=0\%$ & $U=5\%$\\
\midrule
\multicolumn{7}{@{}l}{\textit{(a) Depth $D=2$}}\\
\addlinespace[2pt]
Seismic-bumps & 3600.00 & 1592.23 & 3600.00 & 10.00 & 1631.12 & 4.12\\
Thoracic Surgery & 328.28 & 72.32 & 3600.00 & 0.03 & 0.82 & 0.42\\
House-vote & 1.43 & 3.24 & 6.71 & 0.00 & 0.11 & 0.01\\
Mushroom & 3600.00 & 7.23 & 3600.00 & 0.05 & 0.05 & 0.05\\
Chess & 378.40 & 582.43 & 3600.00 & 0.02 & 4.06 & 3.42\\
Nursery & 3600.00 & 200.82 & 3600.00 & 0.07 & 12.73 & 10.37\\
Cylinder Bands & 1.49 & 2.36 & 27.77 & 1.79 & 0.11 & 0.01\\
Echocardiogram & 1.66 & 1.75 & 7.07 & 0.02 & 0.09 & 0.01\\
Credit Approval & 1291.67 & 182.00 & 3600.00 & 0.23 & 156.58 & 116.14\\
Connectionist Bench (Sonar) & 3600.00 & 3600.00 & 3600.00 & 10.41 & 31.27 & 31.27\\
Spambase & 3600.00 & 3600.00 & 3600.00 & \streedOOM & 3600.00 & 3600.00\\
Blood Transfusion Service Center & 115.24 & 5.85 & 1470.75 & 0.02 & 11.14 & 9.02\\
Breast Cancer Wisconsin (Diagnostic) & 3600.00 & 352.56 & 3600.00 & \streedOOM & 438.24 & 0.10\\
Wine & 107.86 & 16.73 & 35.21 & 0.15 & 19.87 & 19.87\\
Hepatitis & 76.20 & 16.33 & 51.91 & 0.01 & 2.72 & 2.68\\
Ozone Level Detection & 3600.00 & 3600.00 & 3600.00 & \streedOOM & 0.97 & 0.97\\
Car-evaluation & 78.33 & 32.56 & 806.22 & 0.00 & 0.67 & 0.67\\
Balance-scale & 26.09 & 11.59 & 176.94 & 0.00 & 0.23 & 0.19\\
Biodeg & 3600.00 & 432.75 & 3600.00 & \streedOOM & 374.38 & 374.38\\
Iris & 3.56 & 0.26 & 10.31 & 0.01 & 0.09 & 0.09\\
Ionosphere & 3434.51 & 632.57 & 3600.00 & 6.71 & 93.44 & 71.51\\
Tic-tac-toe & 33.98 & 23.54 & 542.90 & 0.02 & 0.49 & 0.49\\
Banknote Authentication & 580.85 & 26.78 & 3600.00 & 4.83 & 151.89 & 121.49\\
Fertility & 593.09 & 6.42 & 3600.00 & 0.00 & 153.45 & 118.64\\
Planning Relax & 345.14 & 85.43 & 3600.00 & 0.24 & 7.79 & 11.92\\
Seeds & 53.36 & 6.87 & 29.72 & 0.13 & 2.90 & 2.87\\
German & 374.24 & 6.71 & 3600.00 & 0.21 & 0.09 & 0.09\\
\midrule
\multicolumn{7}{@{}l}{\textit{(b) Depth $D=3$}}\\
\addlinespace[2pt]
Seismic-bumps & 3600.00 & 3600.00 & 3600.00 & 3465.76 & 3600.00 & 6.92\\
Thoracic Surgery & 3600.00 & 3324.76 & 3600.00 & 0.59 & 0.81 & 0.31\\
House-vote & 76.61 & 20.74 & 456.25 & 0.01 & 0.01 & 0.02\\
Mushroom & 3600.00 & 106.56 & 3600.00 & 0.10 & 0.85 & 0.03\\
Chess & 3600.00 & 3600.00 & 3600.00 & 0.03 & 100.71 & 0.93\\
Nursery & 3600.00 & 3573.78 & 3600.00 & 0.11 & 60.11 & 52.94\\
Cylinder Bands & 73.83 & 103.17 & 409.15 & 98.34 & 2.51 & 0.01\\
Echocardiogram & 90.32 & 114.68 & 441.90 & 0.62 & 2.44 & 0.01\\
Credit Approval & 3600.00 & 3600.00 & 3600.00 & 28.81 & 3600.00 & 3600.00\\
Connectionist Bench (Sonar) & 3600.00 & 3600.00 & 3600.00 & 3600.00 & 3600.00 & 3600.00\\
Spambase & 3600.00 & 3600.00 & 3600.00 & \streedOOM & 3600.00 & 3600.00\\
Blood Transfusion Service Center & 3600.00 & 837.26 & 3600.00 & 0.07 & 599.33 & 523.35\\
Breast Cancer Wisconsin (Diagnostic) & 3600.00 & 3600.00 & 3600.00 & \streedOOM & 3600.00 & 0.10\\
Wine & 3600.00 & 1236.56 & 3600.00 & 28.25 & 3600.00 & 3600.00\\
Hepatitis & 3600.95 & 832.52 & 3600.00 & 0.38 & 378.06 & 398.62\\
Ozone Level Detection & 3600.00 & 3600.00 & 3600.00 & \streedOOM & 0.96 & 0.96\\
Car-evaluation & 3600.00 & 3600.00 & 3600.00 & 0.03 & 14.29 & 10.03\\
Balance-scale & 3600.00 & 1239.66 & 3600.00 & 0.01 & 2.82 & 2.15\\
Biodeg & 3600.00 & 3600.00 & 3600.00 & \streedOOM & 3600.00 & 3600.00\\
Iris & 55.01 & 32.56 & 53.12 & 0.03 & 0.05 & 0.05\\
Ionosphere & 3600.00 & 3600.00 & 3600.00 & 3600.00 & 3600.00 & 1.52\\
Tic-tac-toe & 3600.00 & 3600.00 & 3600.00 & 0.02 & 10.21 & 7.28\\
Banknote Authentication & 3600.00 & 2565.32 & 3600.00 & 1520.23 & 3600.00 & 1.73\\
Fertility & 3600.00 & 1255.74 & 3600.00 & 0.01 & 0.94 & 2.65\\
Planning Relax & 3600.00 & 3600.00 & 3600.00 & 96.59 & 2235.92 & 2235.92\\
Seeds & 3600.00 & 1105.54 & 3600.00 & 18.26 & 290.66 & 285.53\\
German & 3600.00 & 295.84 & 3600.00 & 13.50 & 0.05 & 0.05\\
\bottomrule
\end{tabular*}

\vspace{2pt}\begin{minipage}{\textwidth}\fontsize{8}{9.5}\selectfont\textit{Note.} Time limit of 3,600 seconds. \streedNote\end{minipage}}
\end{table}

\begin{table}[!tbp]\centering
\caption{Running Times (Seconds) for $D=4$ and $D=5$}
\label{tab:comp_time_D4}\label{tab:comp_time_D5}
{\def\baselinestretch{1}\fontsize{8}{9.5}\selectfont\renewcommand{\arraystretch}{0.9}\setlength{\tabcolsep}{3pt}\def\hd#1#2{\begin{tabular}[b]{@{}c@{}}#1\\#2\end{tabular}}
\begin{tabular*}{\textwidth}{@{\extracolsep{\fill}}l*{6}{r}@{}}
\toprule
 & \multicolumn{3}{c}{MIP-based} & & \multicolumn{2}{c}{ISCG-OCT} \\ \cmidrule(lr){2-4}\cmidrule(lr){6-7}
Dataset & \hd{BNP-}{OCT} & \hd{Benders}{OCT} & OCT & STreeD & $U=0\%$ & $U=5\%$\\
\midrule
\multicolumn{7}{@{}l}{\textit{(a) Depth $D=4$}}\\
\addlinespace[2pt]
Seismic-bumps & 3600.00 & 3600.00 & 3600.00 & 3600.00 & 3600.00 & 8.52\\
Thoracic Surgery & 3600.00 & 3231.56 & 3600.00 & 33.30 & 0.81 & 0.34\\
House-vote & 802.57 & 37.85 & 3600.00 & 0.02 & 0.03 & 0.03\\
Mushroom & 3600.00 & 3600.00 & 3600.00 & 1.35 & 14.52 & 0.03\\
Chess & 3600.00 & 3600.00 & 3600.00 & 0.15 & 2585.16 & 0.68\\
Nursery & 3600.00 & 3600.00 & 3600.00 & 0.60 & 3600.00 & 3600.00\\
Cylinder Bands & 750.42 & 1308.89 & 3600.00 & 3600.00 & 47.16 & 0.03\\
Echocardiogram & 794.54 & 1380.04 & 3600.00 & 63.34 & 46.43 & 0.03\\
Credit Approval & 3600.00 & 3600.00 & 3600.00 & 3600.00 & 3600.00 & 3600.00\\
Connectionist Bench (Sonar) & 3600.00 & 3600.00 & 3600.00 & 3600.00 & 3600.00 & 3600.00\\
Spambase & 3600.00 & 3600.00 & 3600.00 & \streedOOM & 3600.00 & 3600.00\\
Blood Transfusion Service Center & 3600.00 & 3600.00 & 3600.00 & 2.04 & 3600.00 & 3600.00\\
Breast Cancer Wisconsin (Diagnostic) & 3600.00 & 3600.00 & 3600.00 & \streedOOM & 3600.00 & 0.10\\
Wine & 3600.00 & 2325.85 & 3600.00 & 3600.00 & 3600.00 & 3600.00\\
Hepatitis & 3600.00 & 1332.56 & 3600.00 & 24.16 & 3600.00 & 3600.00\\
Ozone Level Detection & 3600.00 & 3600.00 & 3600.00 & \streedOOM & 0.96 & 0.96\\
Car-evaluation & 3600.00 & 3600.00 & 3600.00 & 0.06 & 262.64 & 176.63\\
Balance-scale & 3600.00 & 3600.00 & 3600.00 & 0.02 & 34.09 & 22.61\\
Biodeg & 3600.00 & 3600.00 & 3600.00 & \streedOOM & 3600.00 & 3600.00\\
Iris & 637.73 & 1325.32 & 196.79 & 0.56 & 0.09 & 0.09\\
Ionosphere & 3600.00 & 3600.00 & 3600.00 & 3600.00 & 3600.00 & 1.56\\
Tic-tac-toe & 3600.00 & 3600.00 & 3600.00 & 0.06 & 229.68 & 256.68\\
Banknote Authentication & 3600.00 & 3600.00 & 3600.00 & 3600.00 & 3600.00 & 2.92\\
Fertility & 3600.00 & 2022.76 & 3600.00 & 0.03 & 3600.00 & 3.56\\
Planning Relax & 3600.00 & 3600.00 & 3600.00 & 3600.00 & 3600.00 & 3600.00\\
Seeds & 3600.00 & 3600.00 & 3600.00 & 2244.34 & 3600.00 & 3600.00\\
German & 3600.00 & 3600.00 & 3600.00 & 2076.04 & 0.09 & 0.09\\
\midrule
\multicolumn{7}{@{}l}{\textit{(b) Depth $D=5$}}\\
\addlinespace[2pt]
Seismic-bumps & 3600.00 & 3600.00 & 3600.00 & 3600.00 & 3600.00 & 10.43\\
Thoracic Surgery & 3600.00 & 3387.78 & 3600.00 & 736.41 & 0.82 & 0.33\\
House-vote & 3600.00 & 76.92 & 3600.00 & 0.03 & 0.06 & 0.06\\
Mushroom & 3600.00 & 3600.00 & 3600.00 & 1.16 & 190.41 & 0.07\\
Chess & 3600.00 & 3600.00 & 3600.00 & 1.30 & 3600.00 & 2.22\\
Nursery & 3600.00 & 3600.00 & 3600.00 & 4.45 & 3600.00 & 3600.00\\
Cylinder Bands & 3600.00 & 263.37 & 3600.00 & 3600.00 & 890.24 & 0.08\\
Echocardiogram & 3600.00 & 260.25 & 3600.00 & 2094.75 & 1017.44 & 0.06\\
Credit Approval & 3600.00 & 3600.00 & 3600.00 & 3600.00 & 3600.00 & 3600.00\\
Connectionist Bench (Sonar) & 3600.00 & 3600.00 & 3600.00 & 3600.00 & 3600.00 & 3600.00\\
Spambase & 3600.00 & 3600.00 & 3600.00 & \streedOOM & 3600.00 & 3600.00\\
Blood Transfusion Service Center & 3600.00 & 3600.00 & 3600.00 & 26.08 & 3600.00 & 3600.00\\
Breast Cancer Wisconsin (Diagnostic) & 3600.00 & 3600.00 & 3600.00 & \streedOOM & 3600.00 & 0.10\\
Wine & 3600.00 & 2786.77 & 3600.00 & 3600.00 & 3600.00 & 3600.00\\
Hepatitis & 3600.00 & 1386.52 & 3600.00 & 937.97 & 3600.00 & 3600.00\\
Ozone Level Detection & 3600.00 & 3600.00 & 3600.00 & \streedOOM & 0.95 & 0.95\\
Car-evaluation & 3600.00 & 3600.00 & 3600.00 & 0.33 & 575.74 & 302.63\\
Balance-scale & 3600.00 & 3600.00 & 3600.00 & 0.08 & 105.62 & 76.56\\
Biodeg & 3600.00 & 3600.00 & 3600.00 & \streedOOM & 3600.00 & 3600.00\\
Iris & 990.38 & 1445.67 & 407.36 & 5.45 & 0.17 & 0.17\\
Ionosphere & 3600.00 & 3600.00 & 3600.00 & 3600.00 & 3600.00 & 1.64\\
Tic-tac-toe & 3600.00 & 3600.00 & 3600.00 & 0.40 & 3600.00 & 3600.00\\
Banknote Authentication & 3600.00 & 3600.00 & 3600.00 & 3600.00 & 3600.00 & 3.05\\
Fertility & 3600.00 & 3122.55 & 3600.00 & 0.26 & 3600.00 & 4.72\\
Planning Relax & 3600.00 & 3600.00 & 3600.00 & 3600.00 & 3600.00 & 3600.00\\
Seeds & 3600.00 & 3600.00 & 3600.00 & 3600.00 & 3600.00 & 3600.00\\
German & 3600.00 & 3600.00 & 3600.00 & 3600.00 & 0.18 & 0.18\\
\bottomrule
\end{tabular*}

\vspace{2pt}\begin{minipage}{\textwidth}\fontsize{8}{9.5}\selectfont\textit{Note.} Time limit of 3,600 seconds. \streedNote\end{minipage}}
\end{table}

\subsection{Optimality Gap}
\label{gapcomputation}

Tables~\ref{tab:opt_gap_D2} and~\ref{tab:opt_gap_D4} report the optimality gaps at termination. For a minimization problem with best feasible value $UB$ and valid global lower bound $LB$, the relative gap is $\mathrm{Gap}=(UB-LB)/\max\{|UB|,\epsilon\}\times100\%$, where $\epsilon>0$ is a small constant that avoids division by zero. For the MIP-based methods, the bounds are the primal and dual bounds maintained by the solver. ISCG-OCT evaluates objective changes relative to the incumbent tree, whose objective change is zero. By the incumbent-centered coefficients of Section~\ref{subsec:incumbent_change_costs}, correcting an observation misclassified by the incumbent reduces the objective by $2\alpha$, and deactivating an incumbent-active branch node reduces it by $2\beta$. Let $E^{\mathrm{inc}}$ and $A^{\mathrm{inc}}$ be the numbers of misclassified observations and active branch nodes of the incumbent. Since no feasible tree has fewer than zero of either, the objective change of any feasible tree is at least $-2\alpha E^{\mathrm{inc}}-2\beta A^{\mathrm{inc}}$.

A tighter bound is obtained from the recursive pricing of Section~\ref{sec:dynamic_pricing}. Let $\mathcal S$ be a collection of partial states of the root pricing search whose feasible completions, together with the incumbent, cover the complete candidate-restricted MSSF feasible region. For $s\in\mathcal S$, let $\Delta C_s^B$ be the objective change of its fixed branching decisions, $L_s^z$ its hierarchical assignment lower bound under the current routing restrictions, and $U(s)$ its unresolved branch nodes. The incumbent-relative branching coefficients give the valid bound $LB(s)=\Delta C_s^B+L_s^z-2\beta\sum_{w\in U(s)}\mathbf 1\{w\text{ is active in the incumbent}\}$ for all feasible completions of $s$, whose last term is the largest possible reduction in branching cost from unresolved decisions. Combining the two bounds,
\[
LB_{\Delta}=\max\Bigl\{-2\alpha E^{\mathrm{inc}}-2\beta A^{\mathrm{inc}},\ \min\bigl\{0,\min_{s\in\mathcal S}LB(s)\bigr\}\Bigr\}
\]
is a global lower bound on the objective change relative to the incumbent, provided that the maintained states cover the complete candidate-restricted feasible region and their bounds remain valid. For a consistent comparison, the gap is evaluated with the equivalent nonnegative MSSF objective $F=2\alpha E+2\beta A$, which preserves the optimal solutions and absolute objective differences. With $F^{\mathrm{inc}}=2\alpha E^{\mathrm{inc}}+2\beta A^{\mathrm{inc}}$, $UB=F^{\mathrm{inc}}$, and $LB=F^{\mathrm{inc}}+LB_{\Delta}$, the gap of ISCG-OCT is
\[
\mathrm{Gap}_{\mathrm{ISCG}}=\frac{-LB_{\Delta}}{\max\{F^{\mathrm{inc}},\epsilon\}}\times100\%,
\]
the largest remaining objective improvement relative to the incumbent value, accounting for both classification performance and tree complexity. When valid global lower bounds are maintained throughout the complete root pricing search, it provides an optimality-gap certificate comparable to the primal--dual bounds of the MIP-based methods.

STreeD reports no lower bound when it stops at the time limit. Its gap is therefore $0\%$ when it certifies optimality; otherwise it is computed from $UB=F$ of the returned tree and the bound
\[
LB=\min_{1\le L\le\min\{2^{D},K\}}\Bigl\{2\alpha\max\bigl\{E^{\mathrm{eq}},S_{K-L}\bigr\}+2\beta(L-1)\Bigr\},
\]
where $K$ is the number of classes, $E^{\mathrm{eq}}$ is the number of misclassifications forced by identical feature vectors with different labels (for each group of identical training observations, its size minus its largest class count), and $S_{K-L}$ is the number of observations in the $K-L$ smallest classes, a lower bound on the misclassifications of a tree with $L$ leaves, which predicts at most $L$ classes. This bound is valid for every tree of depth at most $D$ but weak, so STreeD runs stopped by the time limit show large gaps.

\begin{table}[!tbp]\centering
\caption{Reported Optimality Gaps for $D=2$ and $D=3$}
\label{tab:opt_gap_D2}\label{tab:opt_gap_D3}
{\def\baselinestretch{1}\fontsize{8}{9.5}\selectfont\renewcommand{\arraystretch}{0.9}\setlength{\tabcolsep}{3pt}\def\hd#1#2{\begin{tabular}[b]{@{}c@{}}#1\\#2\end{tabular}}
\begin{tabular*}{\textwidth}{@{\extracolsep{\fill}}l*{6}{r}@{}}
\toprule
 & \multicolumn{3}{c}{MIP-based} & & \multicolumn{2}{c}{ISCG-OCT} \\ \cmidrule(lr){2-4}\cmidrule(lr){6-7}
Dataset & \hd{BNP-}{OCT} & \hd{Benders}{OCT} & OCT & STreeD & $U=0\%$ & $U=5\%$\\
\midrule
\multicolumn{7}{@{}l}{\textit{(a) Depth $D=2$}}\\
\addlinespace[2pt]
Seismic-bumps & 99\% & 0\% & 97\% & 0\% & 0\% & 0\%\\
Thoracic Surgery & 0\% & 0\% & 0\% & 0\% & 0\% & 100\%\\
House-vote & 0\% & 0\% & 0\% & 0\% & 0\% & 0\%\\
Mushroom & 67\% & 0\% & 67\% & 0\% & 0\% & 99\%\\
Chess & 0\% & 0\% & 0\% & 0\% & 0\% & 97\%\\
Nursery & 98\% & 0\% & 98\% & 0\% & 0\% & 0\%\\
Cylinder Bands & 0\% & 0\% & 0\% & 0\% & 0\% & 0\%\\
Echocardiogram & 0\% & 0\% & 0\% & 0\% & 0\% & 0\%\\
Credit Approval & 0\% & 0\% & 0\% & 0\% & 0\% & 100\%\\
Connectionist Bench (Sonar) & 94\% & 85\% & 93\% & 0\% & 0\% & 100\%\\
Spambase & 98\% & 41\% & 97\% & \streedOOM & 14\% & 19\%\\
Blood Transfusion Service Center & 0\% & 0\% & 0\% & 0\% & 0\% & 100\%\\
Breast Cancer Wisconsin (Diagnostic) & 85\% & 0\% & 86\% & \streedOOM & 0\% & 100\%\\
Wine & 0\% & 0\% & 0\% & 0\% & 0\% & 100\%\\
Hepatitis & 0\% & 0\% & 0\% & 0\% & 0\% & 100\%\\
Ozone Level Detection & 99\% & 99\% & 99\% & \streedOOM & 0\% & 0\%\\
Car-evaluation & 0\% & 0\% & 0\% & 0\% & 0\% & 0\%\\
Balance-scale & 0\% & 0\% & 0\% & 0\% & 0\% & 67\%\\
Biodeg & 97\% & 0\% & 99\% & \streedOOM & 0\% & 100\%\\
Iris & 0\% & 0\% & 0\% & 0\% & 0\% & 0\%\\
Ionosphere & 0\% & 0\% & 0\% & 0\% & 0\% & 100\%\\
Tic-tac-toe & 0\% & 0\% & 0\% & 0\% & 0\% & 0\%\\
Banknote Authentication & 0\% & 0\% & 0\% & 0\% & 0\% & 100\%\\
Fertility & 0\% & 0\% & 0\% & 0\% & 0\% & 100\%\\
Planning Relax & 0\% & 0\% & 0\% & 0\% & 0\% & 83\%\\
Seeds & 0\% & 0\% & 0\% & 0\% & 0\% & 100\%\\
German & 0\% & 0\% & 0\% & 0\% & 0\% & 100\%\\
\midrule
\multicolumn{7}{@{}l}{\textit{(b) Depth $D=3$}}\\
\addlinespace[2pt]
Seismic-bumps & 99\% & 100\% & 99\% & 0\% & 94\% & 94\%\\
Thoracic Surgery & 96\% & 0\% & 94\% & 0\% & 0\% & 60\%\\
House-vote & 0\% & 0\% & 0\% & 0\% & 0\% & 0\%\\
Mushroom & 87\% & 0\% & 85\% & 0\% & 0\% & 0\%\\
Chess & 98\% & 82\% & 100\% & 0\% & 0\% & 0\%\\
Nursery & 99\% & 0\% & 100\% & 0\% & 0\% & 0\%\\
Cylinder Bands & 0\% & 0\% & 0\% & 0\% & 0\% & 50\%\\
Echocardiogram & 0\% & 0\% & 0\% & 0\% & 0\% & 0\%\\
Credit Approval & 97\% & 87\% & 95\% & 0\% & 100\% & 100\%\\
Connectionist Bench (Sonar) & 98\% & 95\% & 95\% & 100\% & 100\% & 100\%\\
Spambase & 99\% & 98\% & 98\% & \streedOOM & 98\% & 98\%\\
Blood Transfusion Service Center & 96\% & 0\% & 96\% & 0\% & 95\% & 100\%\\
Breast Cancer Wisconsin (Diagnostic) & 92\% & 70\% & 94\% & \streedOOM & 100\% & 100\%\\
Wine & 44\% & 0\% & 45\% & 0\% & 100\% & 100\%\\
Hepatitis & 79\% & 0\% & 79\% & 0\% & 0\% & 100\%\\
Ozone Level Detection & 99\% & 100\% & 100\% & \streedOOM & 0\% & 0\%\\
Car-evaluation & 91\% & 72\% & 92\% & 0\% & 0\% & 92\%\\
Balance-scale & 65\% & 0\% & 64\% & 0\% & 0\% & 100\%\\
Biodeg & 99\% & 98\% & 96\% & \streedOOM & 100\% & 100\%\\
Iris & 0\% & 0\% & 0\% & 0\% & 0\% & 0\%\\
Ionosphere & 96\% & 86\% & 96\% & 100\% & 94\% & 94\%\\
Tic-tac-toe & 84\% & 86\% & 83\% & 0\% & 0\% & 100\%\\
Banknote Authentication & 93\% & 0\% & 92\% & 0\% & 98\% & 98\%\\
Fertility & 93\% & 0\% & 93\% & 0\% & 0\% & 97\%\\
Planning Relax & 95\% & 73\% & 94\% & 0\% & 0\% & 0\%\\
Seeds & 50\% & 0\% & 51\% & 0\% & 0\% & 100\%\\
German & 96\% & 0\% & 97\% & 0\% & 0\% & 0\%\\
\bottomrule
\end{tabular*}

\vspace{2pt}\begin{minipage}{\textwidth}\fontsize{8}{9.5}\selectfont\textit{Note.} Time limit of 3,600 seconds. \streedNote\end{minipage}}
\end{table}

\begin{table}[!tbp]\centering
\caption{Reported Optimality Gaps for $D=4$ and $D=5$}
\label{tab:opt_gap_D4}\label{tab:opt_gap_D5}
{\def\baselinestretch{1}\fontsize{8}{9.5}\selectfont\renewcommand{\arraystretch}{0.9}\setlength{\tabcolsep}{3pt}\def\hd#1#2{\begin{tabular}[b]{@{}c@{}}#1\\#2\end{tabular}}
\begin{tabular*}{\textwidth}{@{\extracolsep{\fill}}l*{6}{r}@{}}
\toprule
 & \multicolumn{3}{c}{MIP-based} & & \multicolumn{2}{c}{ISCG-OCT} \\ \cmidrule(lr){2-4}\cmidrule(lr){6-7}
Dataset & \hd{BNP-}{OCT} & \hd{Benders}{OCT} & OCT & STreeD & $U=0\%$ & $U=5\%$\\
\midrule
\multicolumn{7}{@{}l}{\textit{(a) Depth $D=4$}}\\
\addlinespace[2pt]
Seismic-bumps & 99\% & 100\% & 98\% & 100\% & 98\% & 98\%\\
Thoracic Surgery & 99\% & 0\% & 99\% & 0\% & 0\% & 100\%\\
House-vote & 0\% & 0\% & 0\% & 0\% & 0\% & 100\%\\
Mushroom & 98\% & 95\% & 100\% & 0\% & 0\% & 0\%\\
Chess & 98\% & 90\% & 100\% & 0\% & 0\% & 0\%\\
Nursery & 99\% & 98\% & 100\% & 0\% & 100\% & 100\%\\
Cylinder Bands & 0\% & 0\% & 0\% & 100\% & 0\% & 100\%\\
Echocardiogram & 0\% & 0\% & 0\% & 0\% & 0\% & 100\%\\
Credit Approval & 97\% & 92\% & 97\% & 100\% & 100\% & 100\%\\
Connectionist Bench (Sonar) & 98\% & 96\% & 100\% & 100\% & 80\% & 80\%\\
Spambase & 99\% & 98\% & 97\% & \streedOOM & 100\% & 100\%\\
Blood Transfusion Service Center & 98\% & 95\% & 99\% & 0\% & 100\% & 100\%\\
Breast Cancer Wisconsin (Diagnostic) & 95\% & 72\% & 95\% & \streedOOM & 100\% & 100\%\\
Wine & 83\% & 0\% & 85\% & 60\% & 100\% & 100\%\\
Hepatitis & 75\% & 0\% & 76\% & 0\% & 100\% & 100\%\\
Ozone Level Detection & 99\% & 100\% & 100\% & \streedOOM & 100\% & 100\%\\
Car-evaluation & 97\% & 95\% & 96\% & 0\% & 0\% & 100\%\\
Balance-scale & 89\% & 85\% & 88\% & 0\% & 0\% & 0\%\\
Biodeg & 99\% & 99\% & 98\% & \streedOOM & 97\% & 97\%\\
Iris & 0\% & 0\% & 0\% & 0\% & 0\% & 0\%\\
Ionosphere & 98\% & 94\% & 97\% & 100\% & 92\% & 92\%\\
Tic-tac-toe & 93\% & 92\% & 94\% & 0\% & 0\% & 100\%\\
Banknote Authentication & 97\% & 98\% & 96\% & 99\% & 100\% & 100\%\\
Fertility & 97\% & 0\% & 98\% & 0\% & 100\% & 100\%\\
Planning Relax & 98\% & 83\% & 98\% & 100\% & 96\% & 96\%\\
Seeds & 80\% & 62\% & 82\% & 0\% & 100\% & 100\%\\
German & 98\% & 95\% & 98\% & 0\% & 0\% & 0\%\\
\midrule
\multicolumn{7}{@{}l}{\textit{(b) Depth $D=5$}}\\
\addlinespace[2pt]
Seismic-bumps & 99\% & 100\% & 99\% & 100\% & 97\% & 97\%\\
Thoracic Surgery & 99\% & 76\% & 98\% & 0\% & 0\% & 0\%\\
House-vote & 21\% & 0\% & 21\% & 0\% & 0\% & 100\%\\
Mushroom & 99\% & 98\% & 97\% & 0\% & 0\% & 0\%\\
Chess & 98\% & 92\% & 99\% & 0\% & 100\% & 100\%\\
Nursery & 99\% & 99\% & 100\% & 0\% & 99\% & 99\%\\
Cylinder Bands & 21\% & 0\% & 21\% & 100\% & 0\% & 100\%\\
Echocardiogram & 21\% & 0\% & 21\% & 0\% & 0\% & 100\%\\
Credit Approval & 97\% & 95\% & 99\% & 100\% & 97\% & 97\%\\
Connectionist Bench (Sonar) & 98\% & 98\% & 99\% & 100\% & 0\% & 0\%\\
Spambase & 99\% & 98\% & 100\% & \streedOOM & 97\% & 97\%\\
Blood Transfusion Service Center & 99\% & 95\% & 99\% & 0\% & 100\% & 100\%\\
Breast Cancer Wisconsin (Diagnostic) & 95\% & 75\% & 94\% & \streedOOM & 100\% & 100\%\\
Wine & 93\% & 0\% & 93\% & 60\% & 100\% & 100\%\\
Hepatitis & 85\% & 0\% & 86\% & 0\% & 100\% & 100\%\\
Ozone Level Detection & 99\% & 100\% & 97\% & \streedOOM & 0\% & 0\%\\
Car-evaluation & 99\% & 98\% & 97\% & 0\% & 0\% & 0\%\\
Balance-scale & 95\% & 95\% & 95\% & 0\% & 0\% & 0\%\\
Biodeg & 99\% & 99\% & 100\% & \streedOOM & 96\% & 96\%\\
Iris & 0\% & 0\% & 0\% & 0\% & 0\% & 0\%\\
Ionosphere & 97\% & 99\% & 99\% & 100\% & 100\% & 100\%\\
Tic-tac-toe & 97\% & 95\% & 96\% & 0\% & 100\% & 100\%\\
Banknote Authentication & 96\% & 99\% & 98\% & 99\% & 100\% & 100\%\\
Fertility & 96\% & 35\% & 98\% & 0\% & 100\% & 100\%\\
Planning Relax & 99\% & 85\% & 97\% & 100\% & 100\% & 100\%\\
Seeds & 78\% & 75\% & 78\% & 78\% & 0\% & 0\%\\
German & 98\% & 98\% & 98\% & 100\% & 0\% & 0\%\\
\bottomrule
\end{tabular*}

\vspace{2pt}\begin{minipage}{\textwidth}\fontsize{8}{9.5}\selectfont\textit{Note.} Time limit of 3,600 seconds. \streedNote\end{minipage}}
\end{table}

\end{document}